\documentclass{amsart}

\usepackage{graphicx}
\usepackage{amsmath}
\usepackage{algorithm}
\usepackage{blindtext}
\usepackage{enumitem}
\usepackage{framed}
\usepackage[end]{algpseudocode}
\usepackage{multirow}
\usepackage{booktabs}
\usepackage{rotating}

\theoremstyle{definition}

\theoremstyle{remark}

\numberwithin{equation}{section}

\begin{document}

\title{Heavy Ball GMRESR method for Nonsymmetric Linear Systems}

\author{Mei Yang }
\address{Department of Mathematics, The University of Texas at Arlington, Arlington, TX, 76019}
\email{mei.yang@uta.edu}



\date{\today
}



\begin{abstract}
The heavy ball GMRES (HBGMRES) method is one of Krylov subspace methods for linear systems combined with the restarted GMRES and the heavy ball method which is applied in optimization.
HBGMRES not only keeps benefit of the restarted GMRES in limiting memory usage and controlling orthogonalization cost,
 but also is able to cover up the slow convergence problem in the restarted GMRES.
 Another type of Krylov subspace methods is GMRESR which consists of GCR as the outer algorithm and certain steps of GMRES as an inner method.
 Compared with HBGMRES, GMRESR gives the approximately optimal solution over the new search vector gained from GMRES and all previously kept
 search vectors in GCR. Even though GMRESR performs better than HBGMRES, it has slow convergence which is similar to GMRES and the restarted GMRES.
 Inspired by HBGMRES, we present the heavy ball GMRESR method (HBGMRESR) by using HBGMRES as the inner loop
 instead of using GMRES to salvage the lost convergence speed while still keeping the benefit of GMRESR
 that in the sense of a global minimization over some specific part of the Krylov subspace is done.
 Numerical tests on real data are presented to demonstratee the superiority of the new methods over GMRESR and HBGMRES.
\end{abstract}

\maketitle



\section{\bf{Introduction}}

For the solution of systems of linear equations the so-called Krylov subspace methods are widely used, such as GMRES \cite{A}. Given an initial guess $x_0$, the $k$-step approximate solution is $x_k=x_0+z$, and the vector $z$ is sought in the Krylov subspace
\begin{equation}\label{1}
 \mathcal{K}_k(A, r_0)=\text{span}\{r_0, Ar_0,...,A^{k-1}r_0\},
\end{equation}
where $r_0=b-Ax_0$.

 \par For large scale linear systems, GMRES can be very expensive for large $k$ due to the need to store all the basis vectors in the Arnoldi process and quadratically growth with $k$ of orthogonalization cost of basis vectors. Then, the so-called $k$-step restarted GMRES is proposed, which usually is: fix $k$ and repeatedly iterate the $k$-step GMRES with the current initial guess being the very previous $k$-step GMRES solution. We often denote it by GMRES($k$) \cite{A}. In this paper we use REGMRES($k$) or simply REGMRES as the notation when $k$ is either clear from the context or not particularly important for the discussion. Although REGMRES successfully solve the problem of heavy burden memory and orthogonalization cost by limiting the number of Arnoldi process in every GMRES cycle, it possibly has a tendency to stagnate in convergence with relative small $k$. In REGMRES($k$), based on the previous cycle's solution, each current GMRES cylce minimizes the residual by fixing $k$, which can not extend the searching space as big as in the usual GMRES by increasing $k$. Also the information of each GMRES cycle is not well used in REGMRES($k$). Recently, A.Imakura, R.C.Li, and S.L.Zhang introduced a new type of method, Heavy Ball GMRES (HBGMRES) \cite{E}, which combined the REGMRES with the heavy ball method \cite{D}, which is a optimization technique. HBGMRES is able to salvage the lost convergence since it utilizes information of previous cycles.

 The generalized conjugate gradient method (GCR)\cite{N} is another type of Krylov subspace method, which minimizes the residual over $$\text{span}\{U_k\}=\text{span}\{u_1,u_2,...,u_k\}=\mathcal{K}_k(A, r_0),$$
 where $U_k$ satisfies $AU_k=C_k$, and $C_k$ is column orthonormal. In mathematics, CGR and GMRES are equivalent to each other. The $k+1$-step searching direction $u_{k+1}$ is set to be $r_k$ which can be replaced by a better choice in order to accelerate the convergence. Vuik and Van der Vorst proposed the GMRESR \cite{O} method to consist of GCR as the outer alogrithem and $k$-step GMRES method as the inner loop to approximately solve for searching direction $u_{k+1}$ from $Ae_{k}=r_k$. This strategy significantly speeds up the convergence. However, it remains the disadvantage of GMRES because the inner loop is essentially a restarted GMRES with the same operator A. It therefore also displays the severe degradation in the convergence behavior. Inspired by HBGMRES, we use HBGMRES instead of GMRES as the inner loop, which dramatically salvage the cost lost convergence speed in GMRESR. We refer to this algorithm as the heavy ball GMRESR method (HBGMRESR).

 The rest of this paper is organized as follows. Section 2 reviews HBGMRES and GMRESR. Section 3 presents the algorithmic framework of HBGMRESR. Section 4 gives our numerical results and the concluding remarks is showed in section 5.


\section{\bf {HBGMRES and GMRESR}}

\subsection{HBGMRES}
Before talking about HBGMRES, we first review GMRES given in Algorithm 1. Given initial guess $x_0$, write $x=x_0+z$, where $z$ is the solution of the problem
\begin{equation}
\begin{aligned}
 \|r_k\|_2 &=\min_{z\in {\mathcal{K}_k}(A, r_0)}\|b-A(x_0+z)\|_2 \\
   &=\min_{y\in \mathbb{C}^k}\|r_0-V_ky\|_2 \\
   &=\min_{y\in \mathbb{C}^k}\|\beta e_1-\bar{H}_ky\|_2,
\end{aligned}
\end{equation}
where $\text{span}\{V_k\}=\mathcal{K}_k(A, r_0), V_k^TV_k=I_k$, and $\bar{H}_k \in \mathbb{C}^{(k+1)\times k}$ is upper-Hessenberg. In REGMRES($k$), $x_k^{(\ell)}$ is searched over the subspace $\mathcal{K}_k(A, r_0^{(\ell)})$ at each $\ell$-stage. By adding a difference vector of $x^{(\ell)}-x^{(\ell-1)}$, HBGMRES brings sufficient history information before $x^{(\ell)}$ into the current search, usually enough to make up the lost of previous search spaces. The algorithm is given by Algorithm 2.
\begin{algorithm}
\caption{$k$-step GMRES}\label{GMRES}
\begin{flushleft}
Given any initial guess $x_0\in \mathbb{C}^n$ and an integer $k\geq 1$, this algorithm computes
a generalized minimal residual solution to the linear system $Ax=b$.
\end{flushleft}
\hrule
\begin{algorithmic}[1]
\State {$r_0=b-Ax_0, \beta=\|r_0\|_2$;}
\State {$V_{(:,1)}=r_0/\beta, \hat{H}=0_{(k+1)\times k}$}
\For{ $k=1,2,\ldots, $}:
  \State {$f=AV_{(:,i)}$;}
  \State {$\hat{H}_{(1:i,i)}=V^H_{(:,1:i)}f,f=f-V_{(:,1:i)}\hat{H}_{(1,:i,i)}$}
  \State {$\hat{H}_{(i+1,i)}=\|f\|_2$}
 \If {$\hat{H}_{(i+1,i)}>0$}
     \State{$V_{(:,i+1)}=f/\hat{H}_{(i+1,i)}$;}
 \Else
   \State{reset $k=i, \hat{H}=\hat{H}_{(1:i,1:i)}$}
   \State {\bf {break};}
  \EndIf
\EndFor
\State{$y=\arg \min_y\|\hat{H}y-\beta e_1\|_2$}
\State{\Return{$x_{k}=x_0+Vy$ as an approximate solution to $Ax=b$.}}
\end{algorithmic}
\end{algorithm}

\begin{algorithm}
\caption{HBGMRES(k)}\label{HBGMRES}
\begin{flushleft}
Given any initial guess $x_0^{(1)}\in \mathbb{C}^n$ and an integer $k\geq 1$, this algorithm computes
an approximate to the linear system $Ax=b$ via the heavy ball GMRESR.
\end{flushleft}
\hrule
\begin{algorithmic}[1]
\State{$r_0^{(1)}=b-Ax_0^{(1)}, x_0^{(-1)}=0$;}
\For{ $\ell=1,2,\ldots, $}:
 \If {$\|r_0^{(\ell)}\|_2<\text{tol}\times (\|A\|_1\|x_0^{(\ell)}\|_2+\|b\|_2)$}
     \State{\bf{break};}
 \Else
   \State{compute
	\[
	x^{\left(\ell\right)}_k=x_0^{(\ell)} + \arg \min_{z\in\mathcal{K}_k\left(A,r_0^{(\ell)}\right)+\text{span}\left\{x^{\left(\ell\right)}_0-x^{\left(\ell-1\right)}_0\right\}}\left\|b-A(x_0^{(\ell)}+z)\right\|_2\]}
   \State {$x^{\left(\ell+1\right)}_0=x^{\left(\ell\right)}_k$;}
   \State {$r^{\left(\ell+1\right)}_0=b-Ax^{\left(\ell+1\right)}_0$;}
  \EndIf
\EndFor
\State{\Return{$x^{\left(\ell\right)}_0$ as the computed solution to $Ax=b$.}}
\end{algorithmic}
\end{algorithm}
\newpage
\subsection{GMRESR}
GMRESR was proposed by Vuik and Van der Vorst \cite{O}.  This algorithm uses GCR as the outer algorithm and $k$-steps of GMRES as an inner method for computation of $\mathcal{P}_{k,\ell}(A)r_{\ell}$, where $\mathcal{P}_{k,\ell}(A)$ is an $k$-th order polynomial function to approximate $A^{-1}$. The inner GMRES method computes a new search vector by solving the residual equation $Ae_k=r_k$, where $e_{\ell}=x^{\star}-x_{\ell}$ is the error between the exact solution $x^{\star}$ and the $\ell$-th approximate solution $x_{\ell}$. The approximate solution of $Ae_{\ell}=r_\ell$ is denoted by $u_{\ell+1}$. The outer GCR algorithm minimizes the residual over the new search vector and all previously kept search vectors $u_i$. Therefore the approximate solution $x_\ell$ is sought in the subspace $\text{span}\{U_\ell\}$, where $U_\ell=(u_1,u_2,...,u_\ell)$. The algorithm is given by Algorithm 3.
\begin{algorithm}
\caption{GMRESR}\label{GMRESR}
\begin{flushleft}
Given any initial guess $x_0\in \mathbb{C}^n$ and an integer $\ell \geq 1$, this algorithm computes
an approximate to the linear system $Ax=b$ via GMRESR.
\end{flushleft}
\hrule
\begin{algorithmic}[1]
\State{$r_0=b-Ax_0, \ell=0$;}
\For{ $\ell=1,2,\ldots, $}:
     \State{$u_\ell=P_{k,\ell}(A)r_{\ell-1}$;}
     \State{$c_\ell=Au_\ell$}
     \For{$i=1,2,\ldots,\ell-1 $;}
        \State{$\alpha_i=c_i^T c_\ell$;}
        \State{$c_\ell=c_\ell-\alpha_i c_i$;}
        \State{$u_\ell=u_\ell-\alpha_i u_i$;}
     \EndFor
     \State{$u_\ell=u_\ell/\|c_\ell\|_2$;}
     \State{$c_\ell=c_\ell/\|c_\ell\|_2$;}
     \State{$x_\ell=x_{\ell-1}+(c_\ell^T r_{\ell-1})u_\ell$;}
     \State{$r_\ell=r_{\ell-1}-(c_\ell^T r_{\ell-1})c_\ell$;}
     \If {$\|r_\ell\|_2<\text{tol}\times (\|A\|_1\|x_\ell\|_2+\|b\|_2)$}
        \State{\bf{break};}
     \EndIf
\EndFor
\State{\Return{$x_\ell$ as the approximate solution to $Ax=b$.}}
\end{algorithmic}
\end{algorithm}
At line 3, $\mathcal{P}_{k,\ell}$ indicates the GMRES polynomial that is implicitly constructed in $k$ steps of GMRES.

 After $\ell$ steps GMRESR, we have two matrices $U_\ell=(u_1,u_2,...,u_\ell)$ and $C_\ell=AU_\ell$ with the property that $C_\ell^TC_\ell=I_\ell$. GMRESR minimizes the following problem:
\[\min_{x\in \text{span}(U_\ell)}\|b-Ax\|_2.\]
over $\text{span}(U_\ell)$. In GCR, $U_\ell$ satisfies
\begin{equation}\label{2}
\text{span}(U_\ell)=\mathcal{K}_k(A,r_0)
\end{equation}
because of using $u_{\ell+1}=r_\ell$ at line 3. However, in GMRESR, (\ref{2}) no longer holds. At each stage,
\[ u_{\ell+1}\in \mathcal{K}_k(A,r_{\ell}).
\]
Thus
\[x_{\ell+1}\in x_\ell+\mathcal{K}_k(A,r_{\ell}).\]
It shows that the inner loop of GMRESR is essentially a restarted GMRES with the same operator $A$. It therefore would display the tendency of stagnation in GMRES.
\section{\bf{Heavy Ball GMRESR}}
  In GMRESR, the approximate solution $x_\ell$ is in $\text{span}\{U_\ell\}$. Each $u_i$ is an approximate solution to the error $e_{i-1}$. The better approximation of $e_{\ell}$, the better approximation solution of the linear system we would have. By extending the searching subspace in each cycle by adding an difference vector of previous two solutions in order to give an better approximation, HBGMRES effectively salvages the stagnation in GMRES. Therefore, we choose HBGMRES instead of GMRES for the inner loop. At each stage, we merely add the difference vector $u_\ell-u_{\ell-1}$ in order to solve for $u_{\ell+1}$ in $k$ steps HBGMRES. We call this method heavy ball GMRESR (HBGMRESR). In each inner loop, we have
  \[u_{\ell+1}\in \mathcal{K}_k(A,r_\ell)+\text{span}\{u_\ell-u_{\ell-1}\} .\]
  Therefore,
  \[x_{\ell+1}\in x_\ell+\mathcal{K}_k(A,r_{\ell})+\text{span}\{u_\ell-u_{\ell-1}\}.\]
  HBGMRESR is given by Algorithm 4.
\begin{algorithm}
\caption{HBGMRESR}\label{HBGMRESR}
\begin{flushleft}
Given any initial guess $x_0\in \mathbb{C}^n$ and an integer $\ell \geq 1$, this algorithm computes
an approximate to the linear system $Ax=b$ via HBGMRESR.
\end{flushleft}
\hrule
\begin{algorithmic}[1]
\State{$r_0=b-Ax_0, \ell=0$;}
\For{ $\ell=1,2,\ldots, $}:
  \If {$\ell\leq 2$}
     \State{Compute $u_\ell=\arg \min\limits_{z \in \mathcal{K}_k(A,r_{\ell-1})}\|r_{\ell-1}-Az\|_2$ by $k$-step GMRES;}
  \Else
     \State{Compute $u_\ell=\arg \min\limits_{z \in \mathcal{K}_k(A,r_{\ell-1})+\text{span}\{u_\ell-u_{\ell-1}\}}\|r_{\ell-1}-Az\|_2$ by HBGMRES;}
  \EndIf
     \State{$c_\ell=Au_\ell$;}
     \For{$i=1,2,\ldots,\ell-1 $;}
        \State{$\alpha_i=c_i^T c_\ell$;}
        \State{$c_\ell=c_\ell-\alpha_i c_i$;}
        \State{$u_\ell=u_\ell-\alpha_i u_i$;}
     \EndFor
     \State{$u_\ell=u_\ell/\|c_\ell\|_2$;}
     \State{$c_\ell=c_\ell/\|c_\ell\|_2$;}
     \State{$x_\ell=x_{\ell-1}+(c_\ell^T r_{\ell-1})u_\ell$;}
     \State{$r_\ell=r_{\ell-1}-(c_\ell^T r_{\ell-1})c_\ell$;}
     \If {$\|r_\ell\|_2<\text{tol}\times (\|A\|_1\|x_\ell\|_2+\|b\|_2)$}
        \State{\bf{break};}
     \EndIf
\EndFor
\State{\Return{$x_\ell$ as the approximate solution to $Ax=b$.}}
\end{algorithmic}
\end{algorithm}
\section{\bf{Numerical Examples}}
Table 1 lists the numbers of flops for GMRES variants, where (MV) is the number of flops by one matrix-vector multiplication with A. We take it to be twice the number of nonzero entries in A. For simplicity, we only keep the leading terms in flops by three major actions within each inner cycle: matrix-vector multiplication, orthogonalizations, and solutions of the reduced least squares problems. In the $\ell$th outer cycle of GMRESR, we just keep the only leading terms of orthogonalizations.
\begin{table}
\caption{Flops for of GMRES variants}
\begin{tabular}{|c|c|}
  \hline
  per cycle of REGMRES ($k$) & $(k+1)$(MV)+$2k^2n+4k^2$  \\
  \hline
  per cycle of HBGMRES ($k$) & $(k+2)$(MV)+$2(k+2)^2n+4(k+1)^2$  \\
  \hline
  $\ell$th loop of GMRESR($k$)&$(k+2)$(MV)+$2k^2n+4k^2+6n\ell$\\
  \hline
  $\ell$th loop of HBGMRESR($k$)&$(k+3)$(MV)+$2(k+2)^2n+4(k+1)^2+6n\ell$\\
  \hline
\end{tabular}
\end{table}
\vspace{2mm}

In what follows, we will present several numerical tests to compare REGMRES, HBGMRES, GMRESR and HBGMRESR. Comparisons will be done in two aspects:
\begin{itemize}
  \item Normalized resdicual (NRes)
       \[
        \text{NRes}=\frac{\|Ax-b\|_2}{\|A\|_1\|x\|_2+\|b\|_2}
       \]
       against the number of cycles for each methods in Table 1, where using $\|A\|_1$ is for its easiness in computation.
  \item NRes against the numbers of flops for all four methods.
\end{itemize}

 In the selective reorthogonalization case, we firstly use modified Gram-Schmidt (MGC) in stead of classical Gram-Schmidt (CGS) showed at line 5 of Algorithm 1. The MGC is as followed:

\small{
\begin{table}[H]
\caption{Testing Matrices}
\begin{tabular}{|c|c|c|c|c|c|}
  \hline
 ID & matrix & $n$ & \texttt{nnz} & sparsity (\%) & application \\
  \hline
 1 & \texttt{cavity05}     & 1182   & 32747    & 2.3439    & computational fluid dynamics \\
 2 & \texttt{cavity10}     & 2597   & 76367    & 1.1323    & computational fluid dynamics \\
 3 & \texttt{cavity16}     & 4562   & 137887   & 0.6625    & computational fluid dynamics \\
 4 & \texttt{chipcool0}    & 20082  & 281150   & 0.0697    & model reduction problem \\
 5 & \texttt{comsol}       & 1500   & 97645    & 4.3398    & structural problem \\
 6 & \texttt{flowmeter5}   & 9669   & 67391    & 0.0721    & model reduction problem \\
 7 & \texttt{fpga\_trans\_02}& 1220 & 7382     & 0.4960    & circuit simulation \\
 8 & \texttt{memplus}      & 17758  & 126150   & 0.0400    & circuit simulation \\
 9 & \texttt{wang3}        & 26064  & 177168   & 0.0261    & semiconductor device problem \\
 10 & \texttt{raefsky1}    & 3242   & 294276   & 2.7998    & computational fluid dynamics \\
 11 & \texttt{raefsky2}    & 3242   & 294276   & 2.7998    & computational fluid dynamics \\
 12 & \texttt{ns3da}       & 20414  & 1679599  & 0.4000    & computational fluid dynamics \\
 13 & \texttt{atmosmodd}   & 1270432& 8814880  & 0.0005    & computational fluid dynamics \\
 14 & \texttt{atmosmodj}   & 1270432& 8814880  & 0.0005    & computational fluid dynamics \\
 15 & \texttt{atmosmodl}   & 1489752& 10319760 & 0.0005    & computational fluid dynamics \\
 16 & \texttt{atmosmodm}   & 1489752& 10319760 & 0.0005    & computational fluid dynamics \\
 \hline
\end{tabular}
\end{table}
}
\vspace{2mm}

All testing matrices are taken from the University of Florida sparse matrix collection \cite{Q} with no particular preference in their selections. Each comes with their right-hand sides $b$. Different from No. 1-11, examples No. 12-14 are all double linear systems. For simplicity, we select the first right-hand sides $b$ as our test vectors.  Table 2 lists 15 testing examples and their characteristics, where $n$ is the size of the matrix, \texttt{nnz} is the number of nonzero entries, and sparsity is $\texttt{nnz}/n^2$.

In running our numerical tests, in order to be fair for comparing the algorithms, we use the following two scenarios
\begin{itemize}
  \item We will run REGMRES($k+1$), HBGMRES($k$), GMRESR($k+1$), and HBGMRESR($k$) in order to illustrate convergence against cycle indices. One consideration is to make sure that all approximation solutions at a cycle of REGMRES and HBGMRES computed from a subspace of dimension $k+2$: $\text{span}\{x_0\}+\mathcal{K}_{k+1}(A,r_0)$ for REGMRES($k+1$), and $\text{span}\{x_0\}+\mathcal{K}_{k}(A,r_0)+\text{span}\{d\}$ for HBGMRES($k$). The same idea is used to make GMRESR and HBGMRESR have the same dimension subspace in an inner cycle of GMRESR($k+1$) and HBGMRESR($k$).
  \item We will experiment two different reorthogonalization strategies in the Arnoldi processes for computing orthonormal bases of Krylov subspaces and in the additional orthogonalizations in each method.  They are selective reorthogonalization and always reorthogonalization. The strategy is as same as in \cite{E}.
  \item The stopping criteria is either NRes is less than or equal to $10^{-12}$ or the number of cycles exceeds 1000.
\end{itemize}

We now report our numerical results for the testing examples listed in Table 2. For each matrix, we ran REGMRES(31), HBGMRES(30), GMRESR(31), and HBGMRESR(30) on the vector b that comes with the matrix and a randomly generated $b$,  with selective and always reorthogonalization. Table 3 lists the number of cycles needed by the algorithms to achieve NRes less or equal $10^{-12}$. In the table, RE, HB, GR and HBR stand for REGMRES, HBGMRES, GMRESR and HBGMRESR. The table clearly demonstrates huge savings achieved by HBGMRESR over the other three.

To get a better impression as how each algorithm behaves,  we plot Figures 1 - 6 to show how NRes in the computed solutions moves against the cycle index. We make the following observations from Figures 1 - 8 and Table 3.
\begin{itemize}
  \item HBGMRESR is overwhelmingly faster than REGMRES, HBGMRES and GMRESE on all examples. We see that the number of iterative cycles of all the examples are much less than that of the other three. This advantage could make up the extra flops cost in every loop of HBGMRESR. In addition, we can tell that GMRESR is faster than REGMRES and HBGMRES for most of examples except for \texttt{memplus},\texttt{fpga\_trans\_02} and \texttt{atmosmodm}.
  \item For HBGMRESR and GMRESR, there is almost no differences in the number of iterative cycles for with always reorthogonalization or selective reorthogonalization. This makes selective reorthogonalization a better choice for cost consideration.
  \item HBGMRESR successfully salvaged the stagnation problem in GMRESR.
\end{itemize}
\newpage
\begin{table}
  \caption{Number of Iterative Cycles}\label{tbl:cycles}
  \centering
  \begin{tabular}{|c|c|c|c|c|c|c|c|c|c|}
  \hline
  \multirow{2}{*}{matrix} & \multicolumn{4}{|c|}{always reorth} &\multicolumn{4}{|c|}{selective reorth} &\multirow{2}{*}{$b$} \\
  \cline{2-9}
     & RE & HB & GR & HBR & RE & HB & GR & HBR &  \\ \hline
  \hline
  \texttt{cavity05}   & 510 & 77   & 30    & 5  & 275  & 77  & 30  & 5 & \multirow{15}{*}{\rotatebox{270}{original}}\\
  \texttt{cavity10}     &1878 & 172   & 52    & 6  &1587 & 172  & 52  & 6 &\\
  \texttt{cavity16}     & -   & 206   & 65    & 6  & -   & 206  & 65  & 6 &\\
  \texttt{chipcool0}    & -   & 517   & 225   & 17 & -   & 518  & 225 & 17& \\
  \texttt{comsol}       &2910 & 284   & 25    & 7  & -   & 308  & 25  & 7 & \\
  \texttt{flowmeter5}   & -   &1482   & 364   &41  & -   & 1414 & 365 & 42& \\
  \texttt{fpga\_trans\_02}&156& 39    &40     & 4  &156  & 39   & 40  & 4 & \\
  \texttt{memplus}      & 84  & 39    & 35    & 4  & 84  & 39   & 35  & 4 & \\
  \texttt{wang3}        & 27  & 24    & 14    & 3  & 27  & 24   & 14  & 3 & \\
  \texttt{raefsky1}     & 208 &39     & 15    & 4  & 204 &39    & 15  & 4 & \\
  \texttt{raefsky2}     & 189 & 137   & 24    & 6  & 189 & 142  & 24  & 6 & \\
  \texttt{ns3da}        & 77  & 72    & 43    & 5  & 77  & 72   & 43  & 5 & \\
  \texttt{atmosmodd}    & 31  & 30    & 15    & 3  & 31  & 30   & 15  & 3 & \\
  \texttt{atmosmodj}    & 56  & 25    & 12    & 3  & 56  & 25   & 12  & 3 & \\
  \texttt{atmosmodl}    & 11  & 11    & 8     & 3  & 11  & 11   & 8   & 3 & \\
  \texttt{atmosmodm}    & 8   & 7     & 7     & 3  & 8   & 7    & 7   & 3 & \\
 \hline
 \hline
 \texttt{cavity05}   & 363 & 97   & 28    & 5  & 411  & 97  & 28  & 5 & \multirow{15}{*}{\rotatebox{270}{random}}\\
  \texttt{cavity10}     &1736 & 198  & 43    & 6  &1618& 198  & 43  & 6 &\\
  \texttt{cavity16}     & -   & 535  & 52    & 7  & -  & 604  & 52  & 7 &\\
  \texttt{chipcool0}    & -   & 491  & 218   & 16 & -  & 501  & 218 &16 &\\
  \texttt{comsol}       & -   & 264  & 42    & 7  & -  & 285  & 42  & 7 &\\
  \texttt{flowmeter5}   & -   & 1586 & 343   & 43 & -  & 1584 & 348 &43 &\\
  \texttt{fpga\_trans\_02}&132& 41   & 45    & 4  & 132& 41   & 45  & 4 & \\
  \texttt{memplus}      & 206 & 53   & 50    & 4  & 206& 53   & 50  & 4 & \\
  \texttt{wang3}        & 31  & 21   & 14    & 3  & 31 & 21   & 14  & 3 & \\
  \texttt{raefsky1}     & 174 &52    & 15    & 4  & 195& 52   & 15  & 4 & \\
  \texttt{raefsky2}     & 191 & 135  & 26    & 6  & 192& 114  & 26  & 6 & \\
  \texttt{ns3da}        & 74  & 76   & 45    & 5  & 74 & 76   & 45  & 5 & \\
  \texttt{atmosmodd}    & 29  & 31   & 18    & 3  & 29 & 31   & 18  & 3 & \\
  \texttt{atmosmodj}    & 69  & 25   & 17    & 3  & 69 & 25   & 17  & 3 & \\
  \texttt{atmosmodl}    & 16  & 13   & 12    & 3  & 16 & 13   & 12  & 3 & \\
  \texttt{atmosmodm}    & 10  & 10   & 9     & 3  & 10 & 10   & 9   & 3 & \\
 \hline

\end{tabular}
\end{table}

\vspace{2mm}
\newpage
\begin{figure}
{\centering
\begin{tabular}{cc}
\hspace{-0.3 cm}
\resizebox*{0.48\textwidth}{0.240\textheight}{\includegraphics{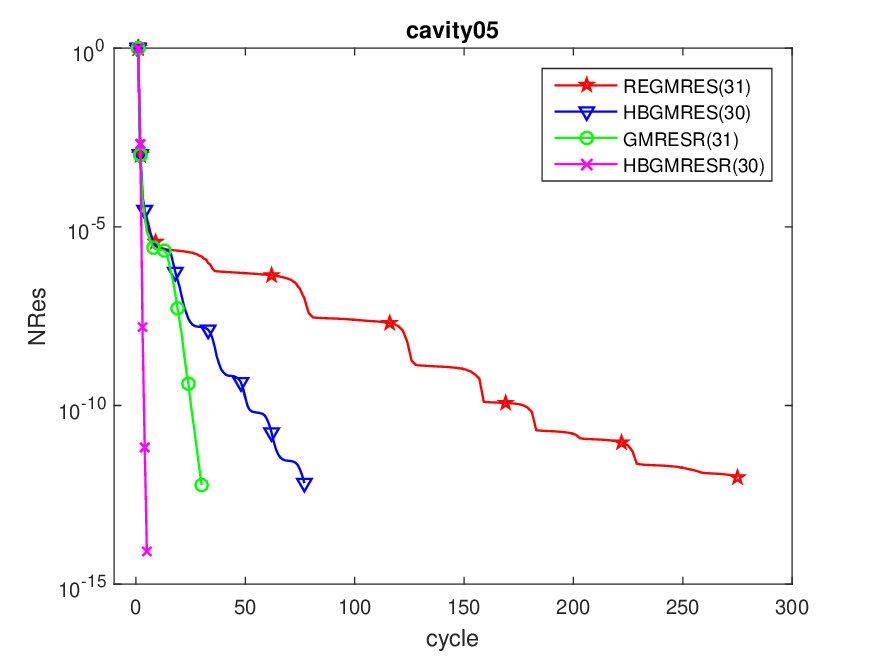}}
&  \hspace{-0.5 cm}
\resizebox*{0.48\textwidth}{0.240\textheight}{\includegraphics{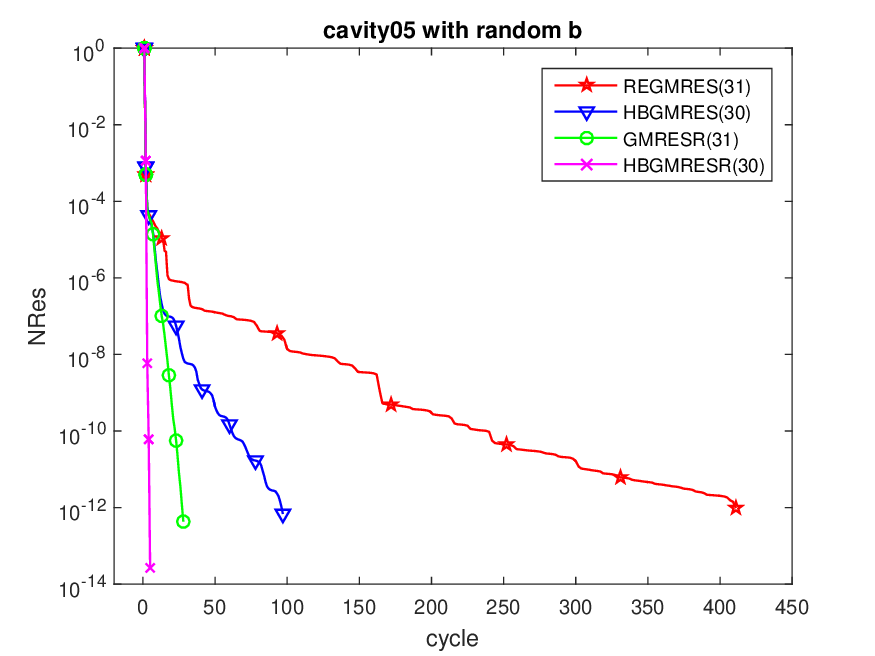}} \\
\hspace{-0.3 cm}
\resizebox*{0.48\textwidth}{0.240\textheight}{\includegraphics{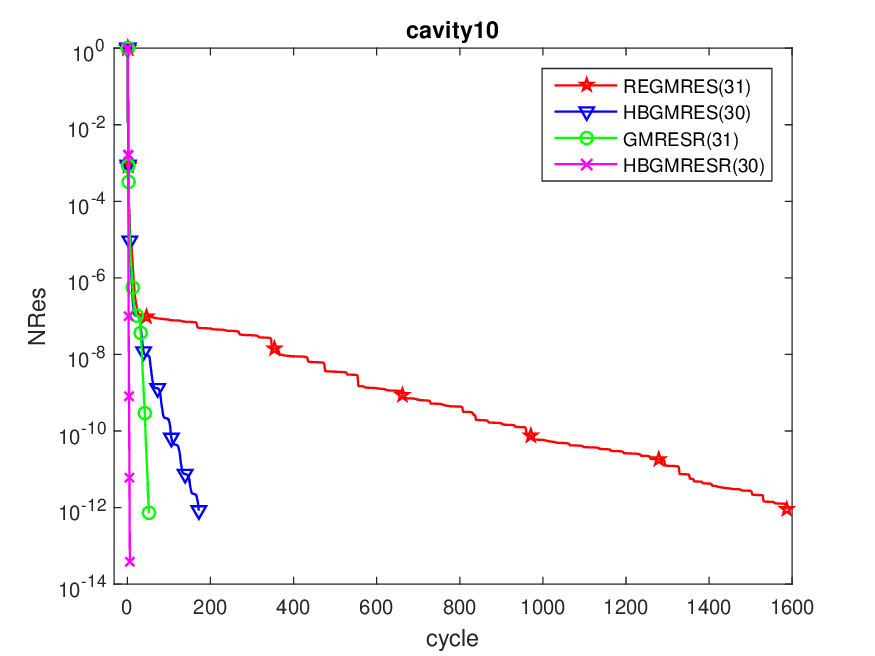}}
&  \hspace{-0.5 cm}
\resizebox*{0.48\textwidth}{0.240\textheight}{\includegraphics{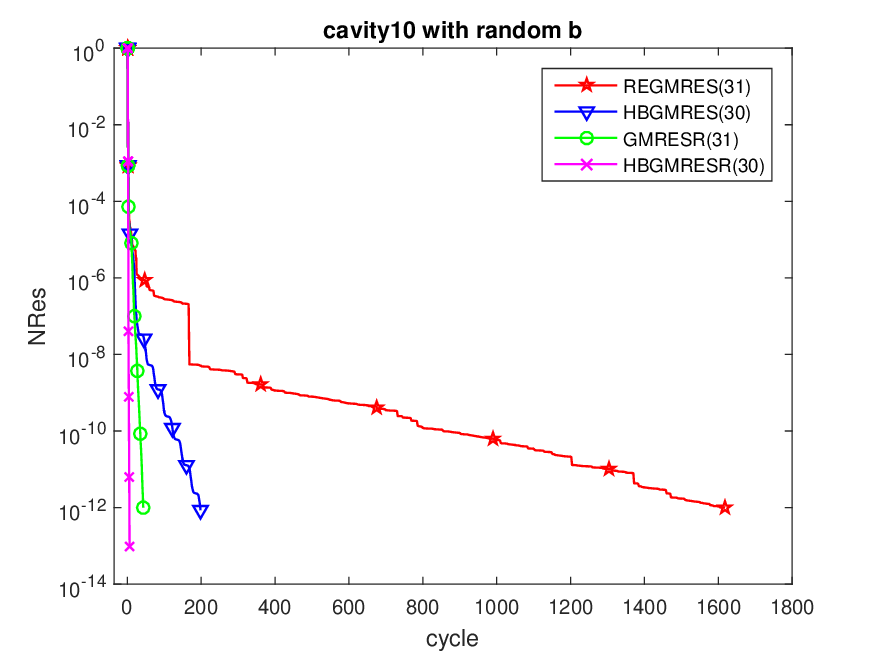}} \\
\hspace{-0.3 cm}
\resizebox*{0.48\textwidth}{0.240\textheight}{\includegraphics{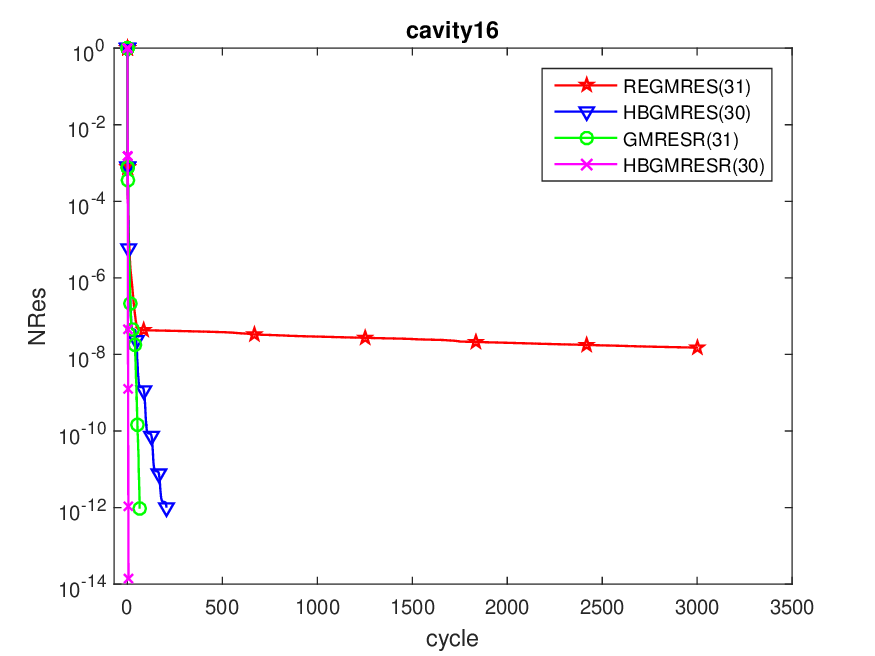}}
&  \hspace{-0.5 cm}
\resizebox*{0.48\textwidth}{0.240\textheight}{\includegraphics{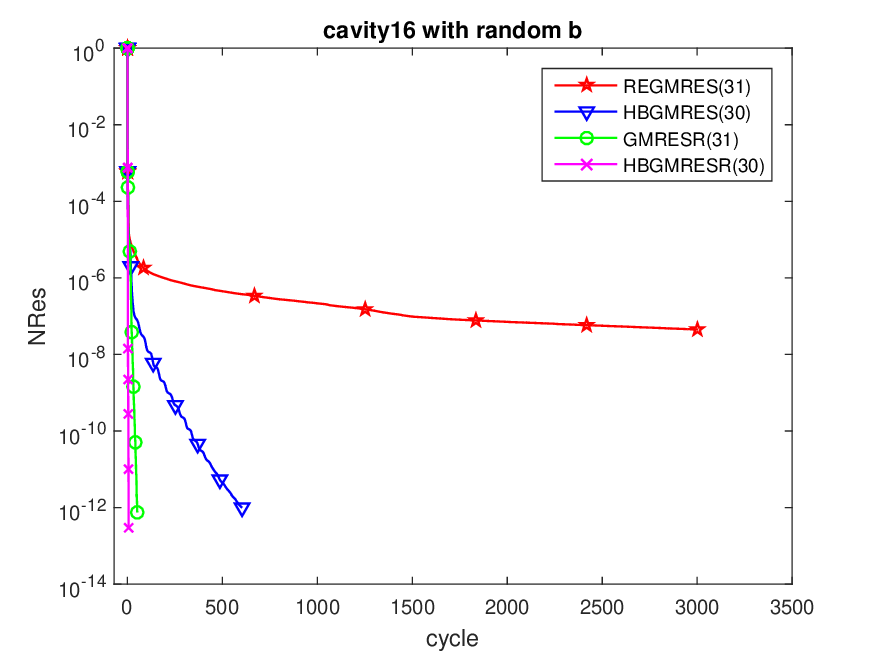}} \\
\hspace{-0.3 cm}
\resizebox*{0.48\textwidth}{0.240\textheight}{\includegraphics{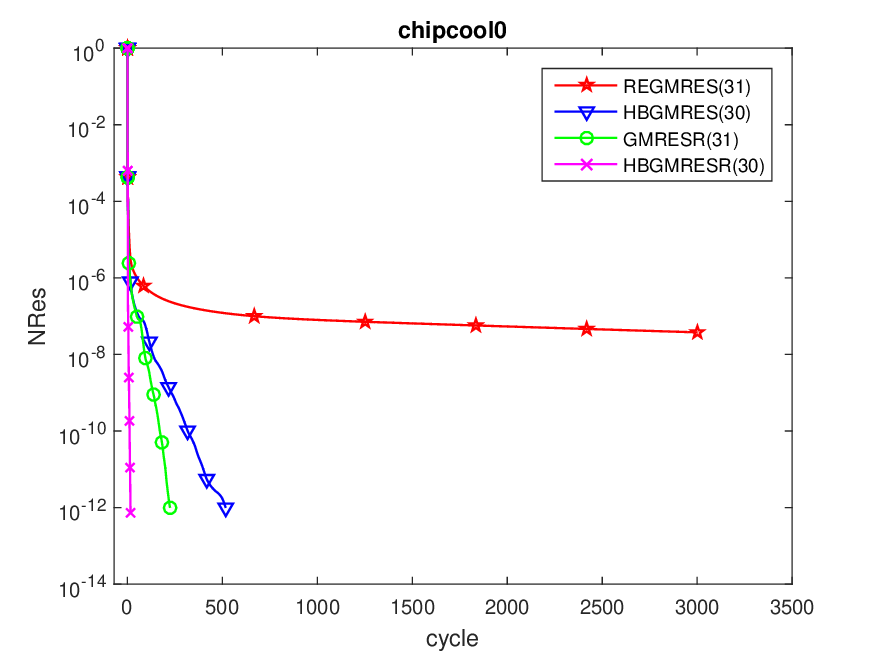}}
&  \hspace{-0.5 cm}
\resizebox*{0.48\textwidth}{0.240\textheight}{\includegraphics{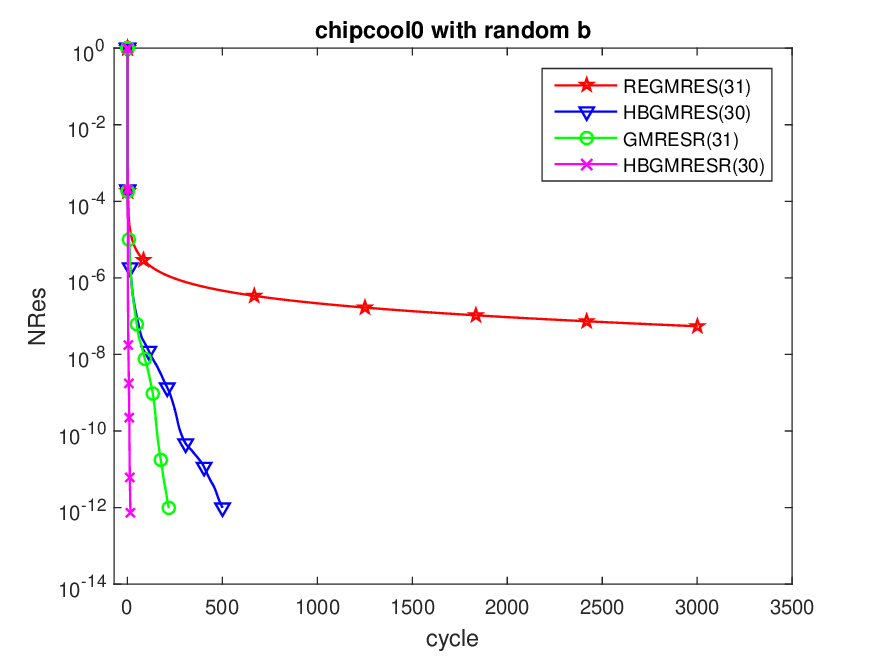}}
\end{tabular}\par
}\vspace{-0.15 cm}
\caption{\small 
    NRes {\em vs.} cycle for
   {\tt cavity05}, {\tt cavity10}, {\tt cavity16}, and {\tt chipcool0}
   with selective reorthogonalization. {\em Left:\/} original $b$; {\em Right:\/} random $b$.
   }
\label{fig:1st4}
\end{figure}
\vspace{2mm}
\newpage
\begin{figure}
{\centering
\begin{tabular}{cc}
\hspace{-0.3 cm}
\resizebox*{0.48\textwidth}{0.240\textheight}{\includegraphics{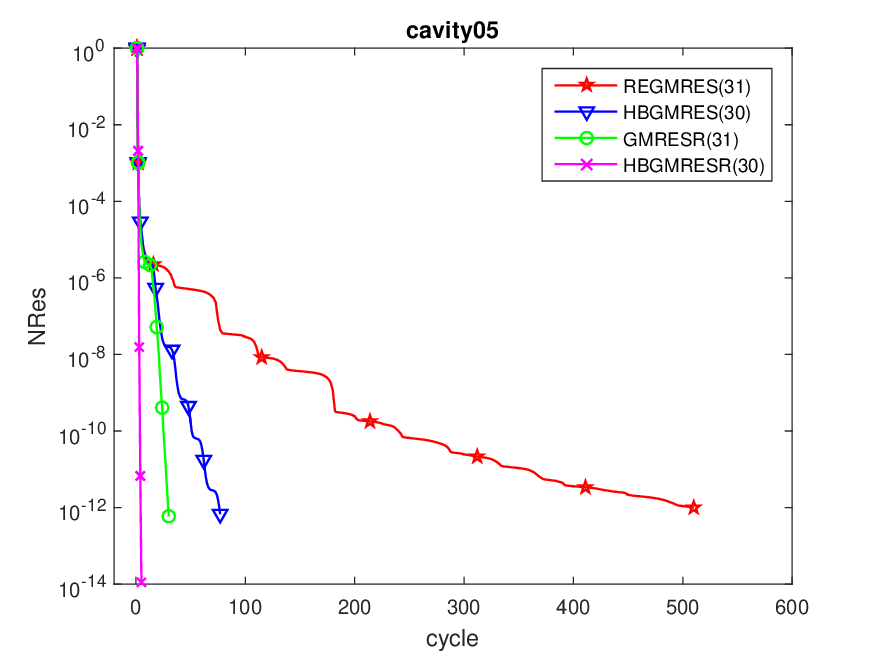}}
&  \hspace{-0.5 cm}
\resizebox*{0.48\textwidth}{0.240\textheight}{\includegraphics{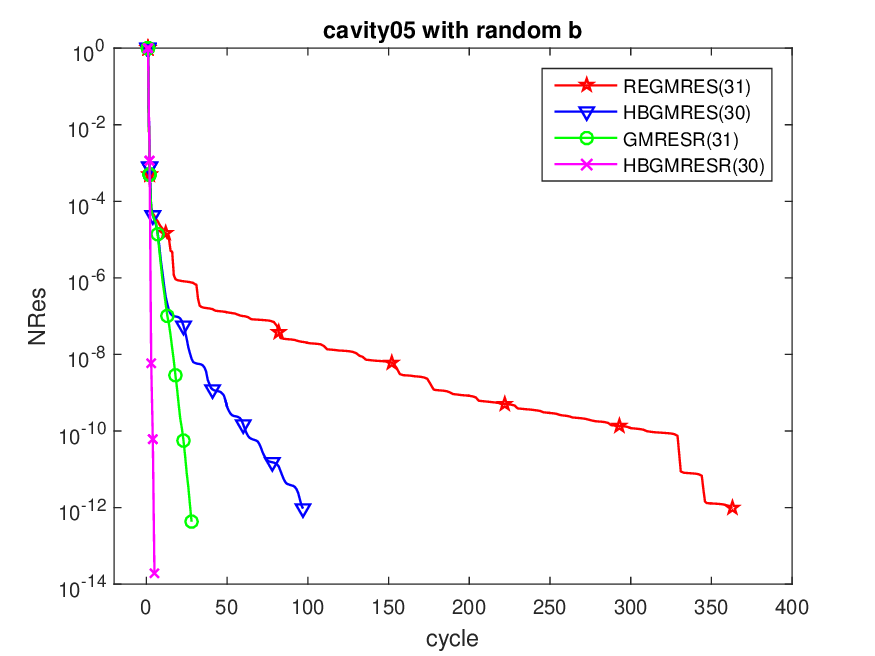}} \\
\hspace{-0.3 cm}
\resizebox*{0.48\textwidth}{0.240\textheight}{\includegraphics{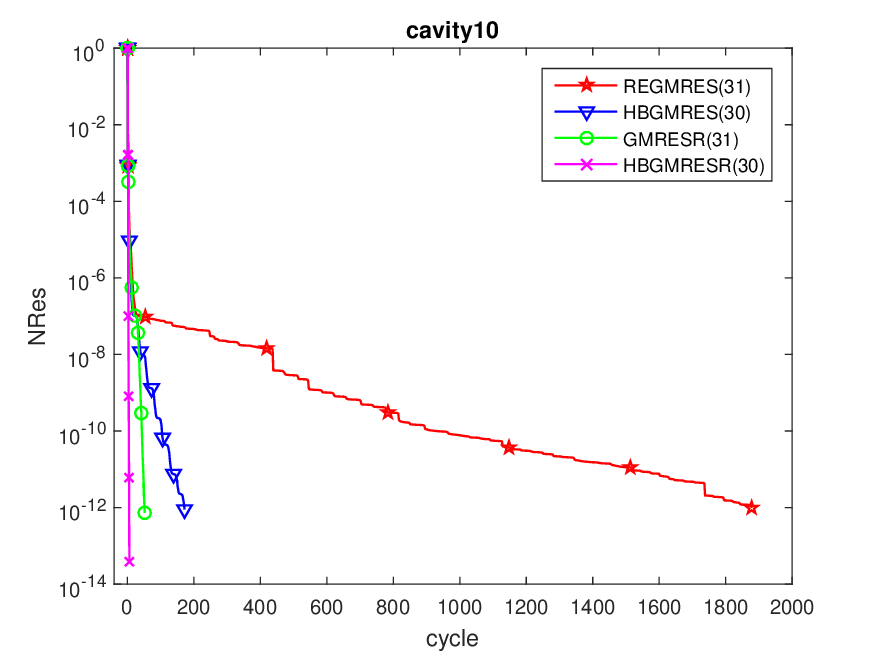}}
&  \hspace{-0.5 cm}
\resizebox*{0.48\textwidth}{0.240\textheight}{\includegraphics{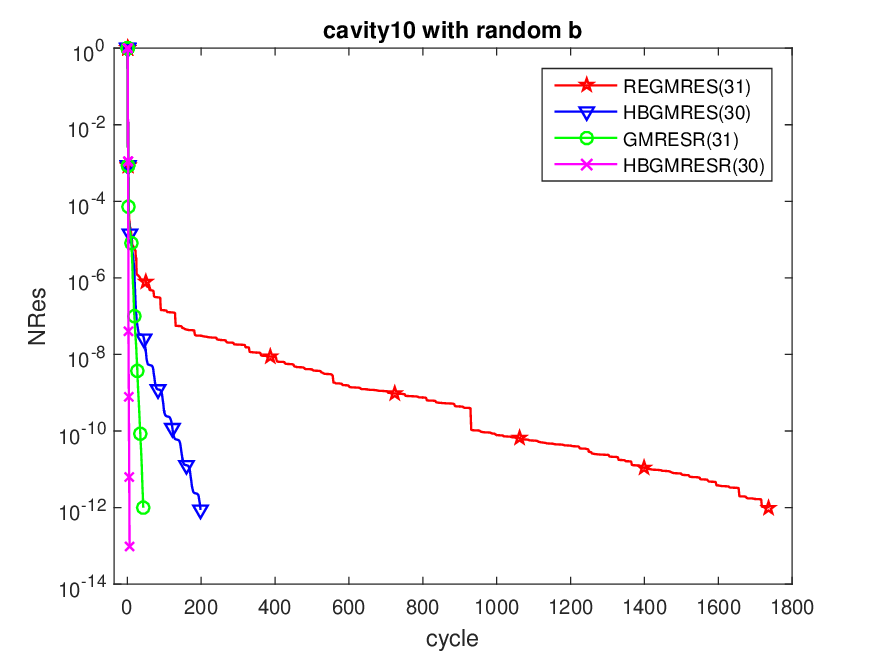}} \\
\hspace{-0.3 cm}
\resizebox*{0.48\textwidth}{0.240\textheight}{\includegraphics{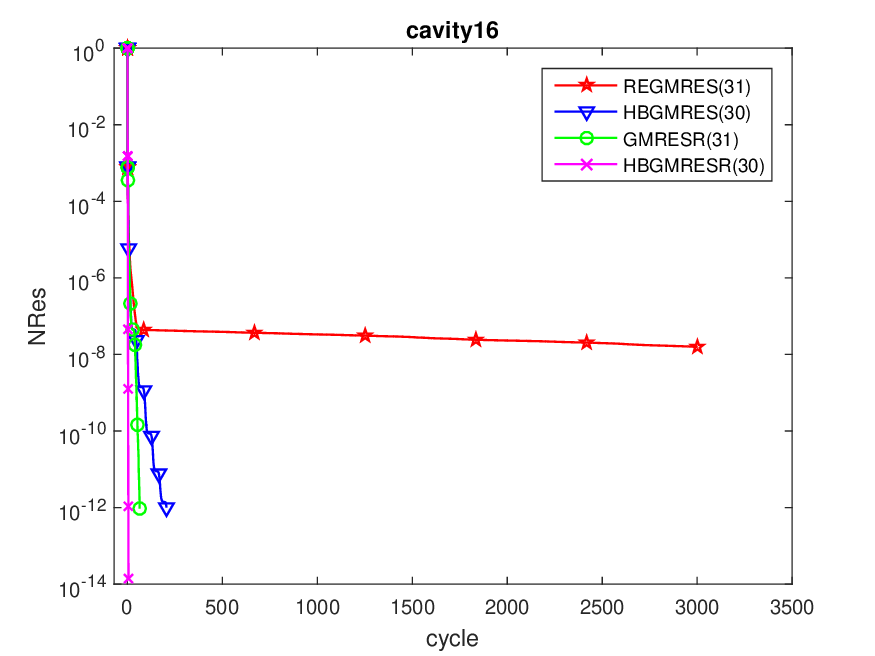}}
&  \hspace{-0.5 cm}
\resizebox*{0.48\textwidth}{0.240\textheight}{\includegraphics{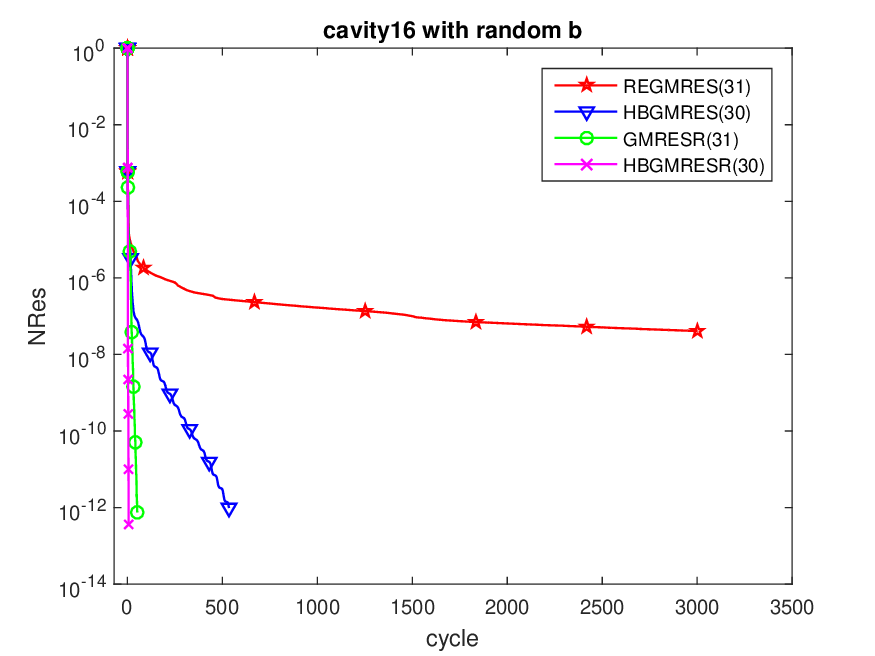}} \\
\hspace{-0.3 cm}
\resizebox*{0.48\textwidth}{0.240\textheight}{\includegraphics{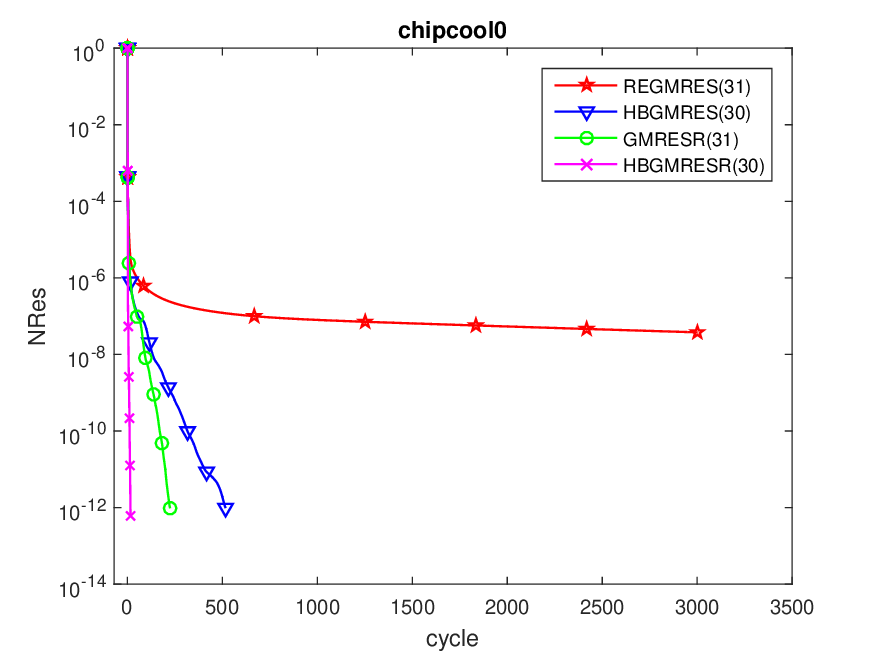}}
&  \hspace{-0.5 cm}
\resizebox*{0.48\textwidth}{0.240\textheight}{\includegraphics{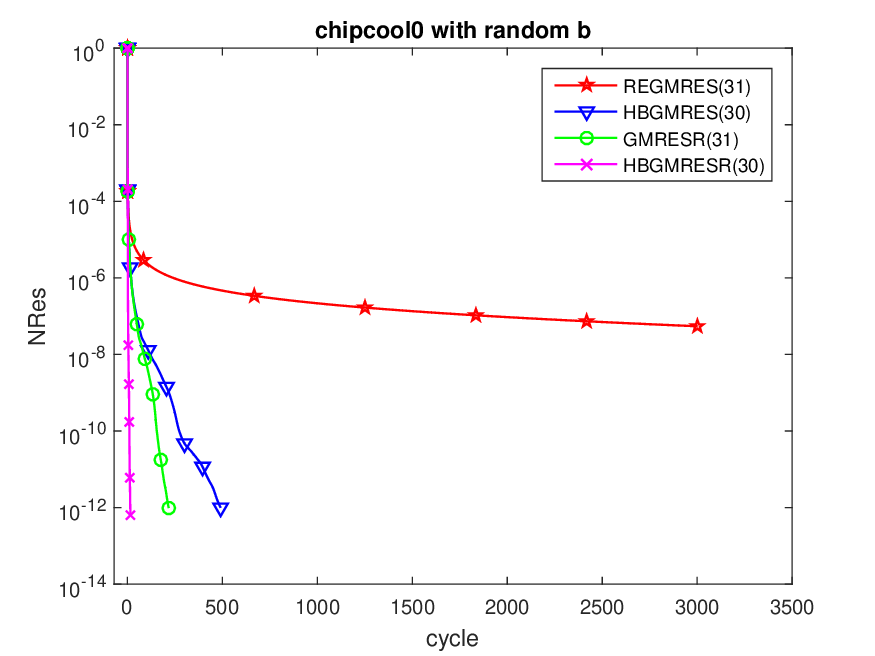}}
\end{tabular}\par
}\vspace{-0.15 cm}
\caption{\small 
    NRes {\em vs.} cycle for
   {\tt cavity05}, {\tt cavity10}, {\tt cavity16}, and {\tt chipcool0}
   with always reorthogonalization. {\em Left:\/} original $b$; {\em Right:\/} random $b$.
   }
\label{fig:1st4w}
\end{figure}
\vspace{2mm}

\newpage
\begin{figure}
{\centering
\begin{tabular}{cc}
\hspace{-0.3 cm}
\resizebox*{0.48\textwidth}{0.240\textheight}{\includegraphics{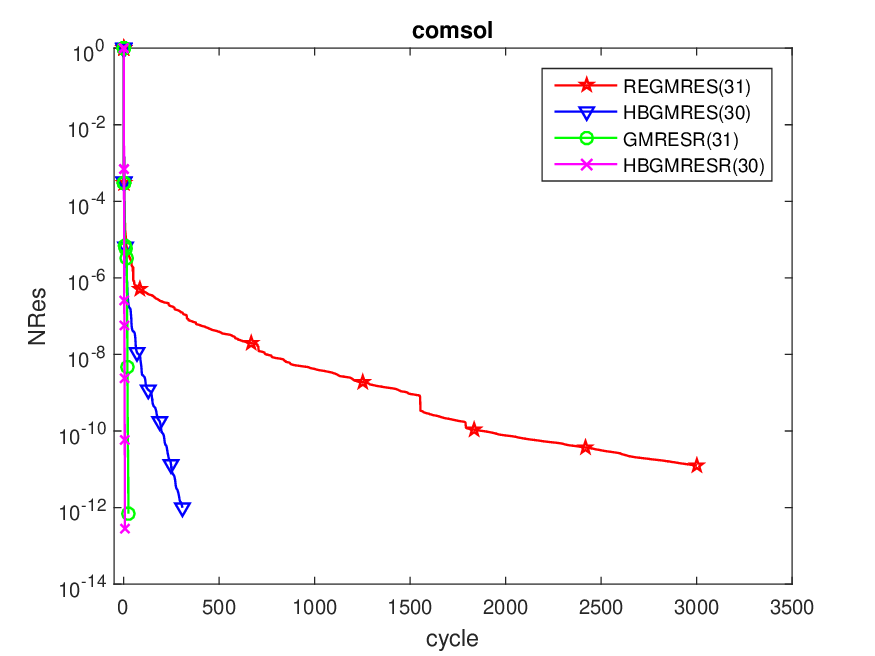}}
&  \hspace{-0.5 cm}
\resizebox*{0.48\textwidth}{0.240\textheight}{\includegraphics{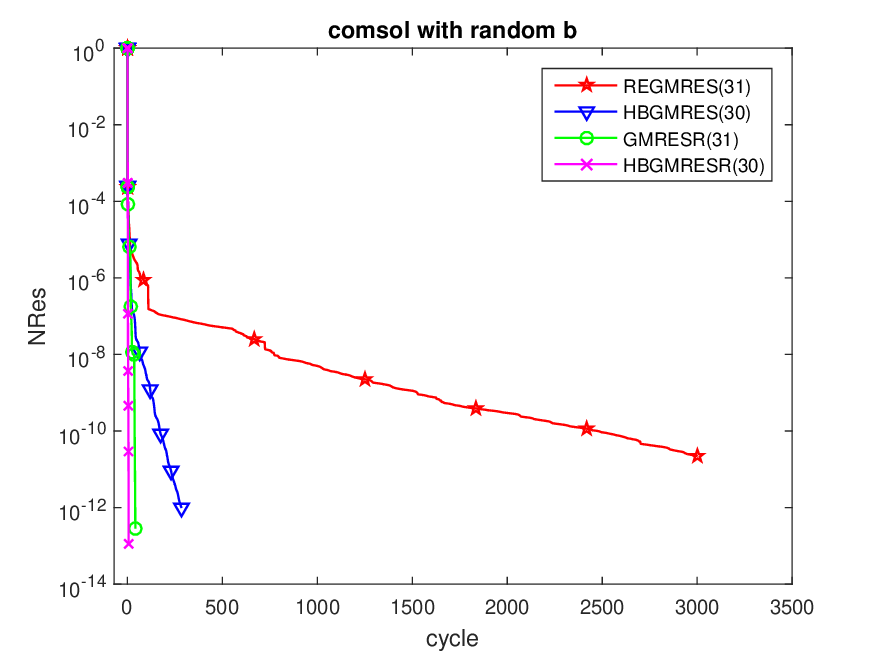}} \\
\hspace{-0.3 cm}
\resizebox*{0.48\textwidth}{0.240\textheight}{\includegraphics{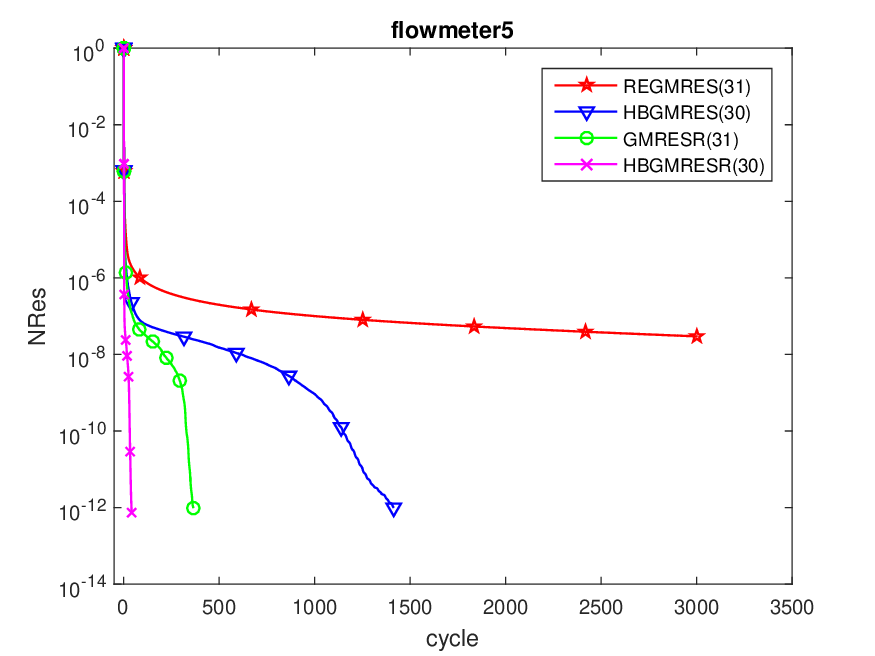}}
&  \hspace{-0.5 cm}
\resizebox*{0.48\textwidth}{0.240\textheight}{\includegraphics{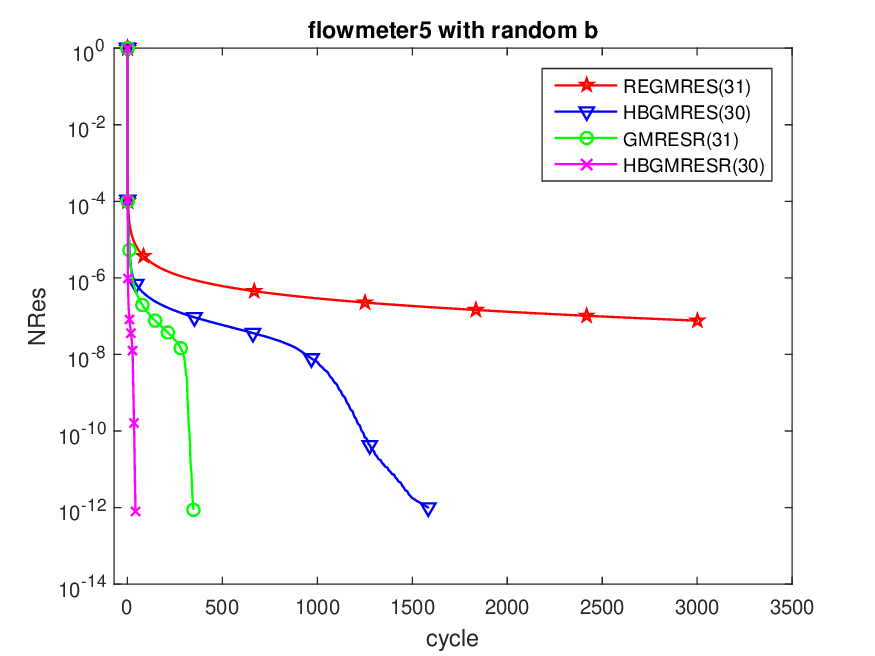}} \\
\hspace{-0.3 cm}
\resizebox*{0.48\textwidth}{0.240\textheight}{\includegraphics{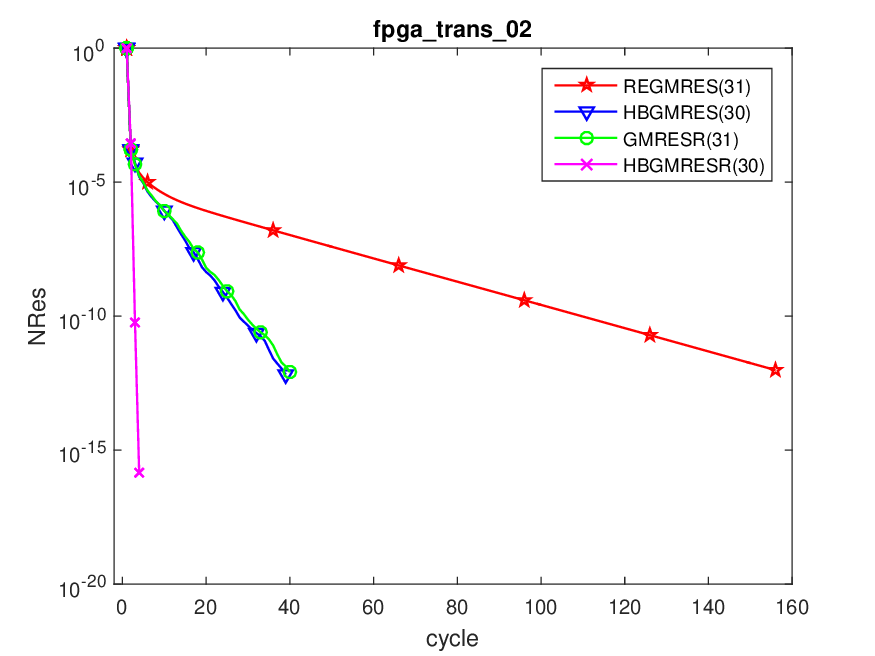}}
&  \hspace{-0.5 cm}
\resizebox*{0.48\textwidth}{0.240\textheight}{\includegraphics{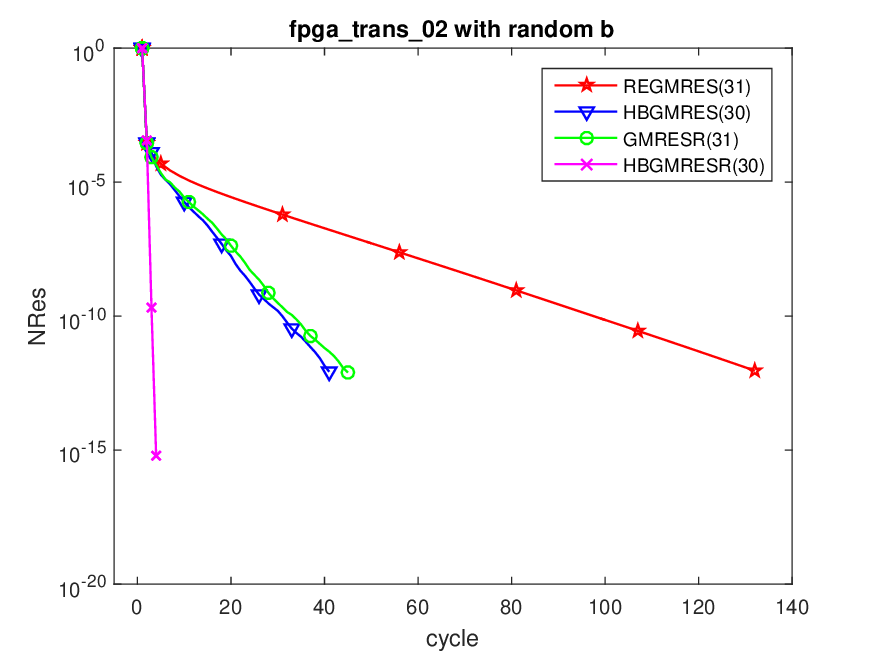}} \\
\hspace{-0.3 cm}
\resizebox*{0.48\textwidth}{0.240\textheight}{\includegraphics{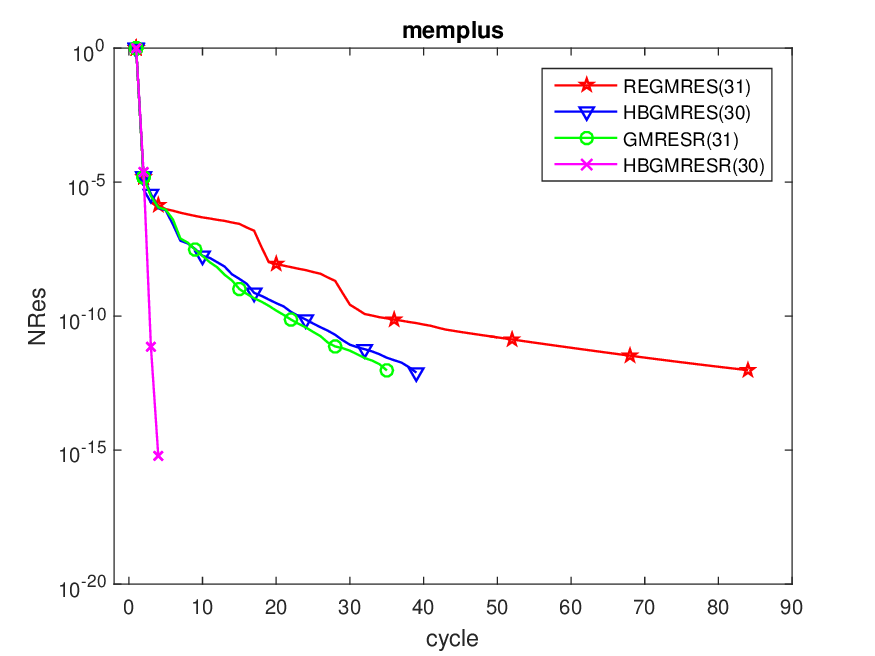}}
&  \hspace{-0.5 cm}
\resizebox*{0.48\textwidth}{0.240\textheight}{\includegraphics{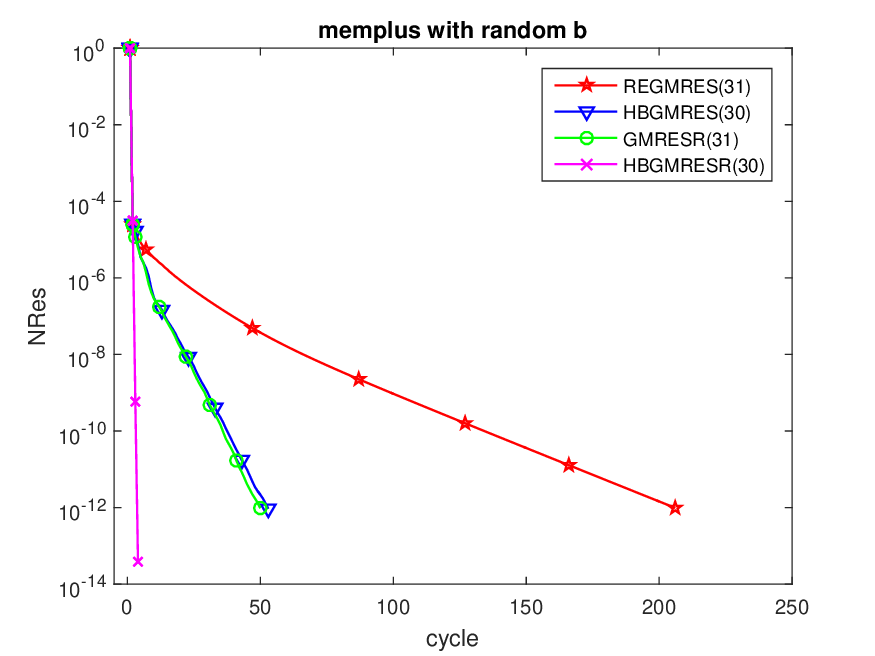}}
\end{tabular}\par
}\vspace{-0.15 cm}
\caption{\small 
    NRes {\em vs.} cycle for
   {\tt comsol}, {\tt flowmeter5}, {\tt fpga\_trans\_02}, and {\tt memplus}
   with selective reorthogonalization. {\em Left:\/} original $b$; {\em Right:\/} random $b$.
   }
\label{fig:2nd4}
\end{figure}
\vspace{2mm}

\newpage
\begin{figure}
{\centering
\begin{tabular}{cc}
\hspace{-0.3 cm}
\resizebox*{0.48\textwidth}{0.240\textheight}{\includegraphics{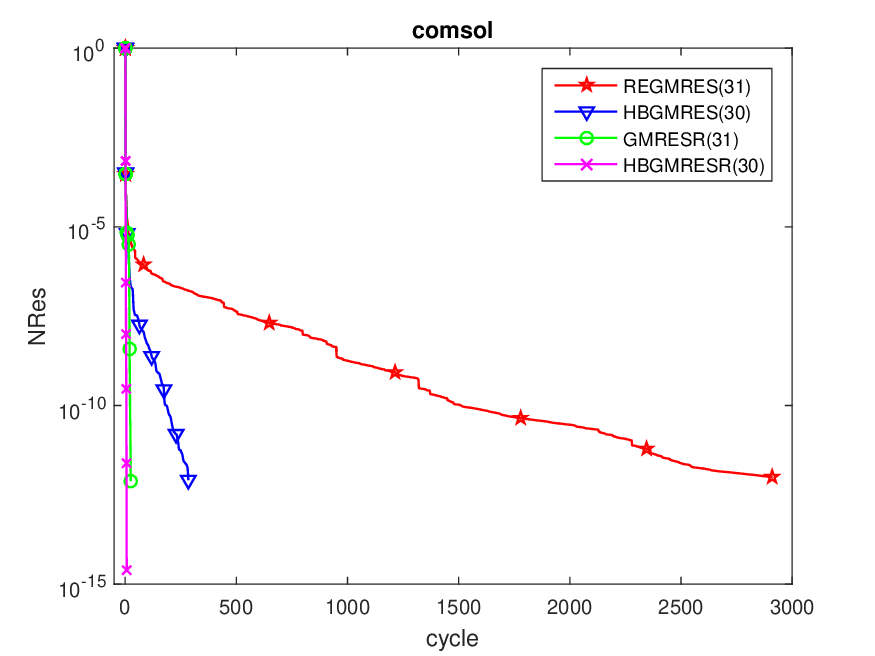}}
&  \hspace{-0.5 cm}
\resizebox*{0.48\textwidth}{0.240\textheight}{\includegraphics{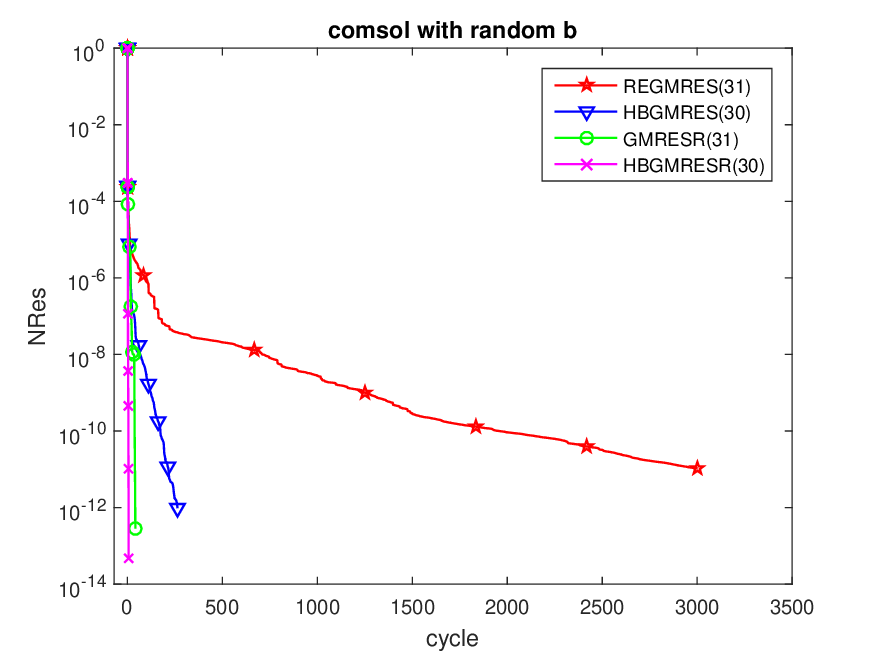}} \\
\hspace{-0.3 cm}
\resizebox*{0.48\textwidth}{0.240\textheight}{\includegraphics{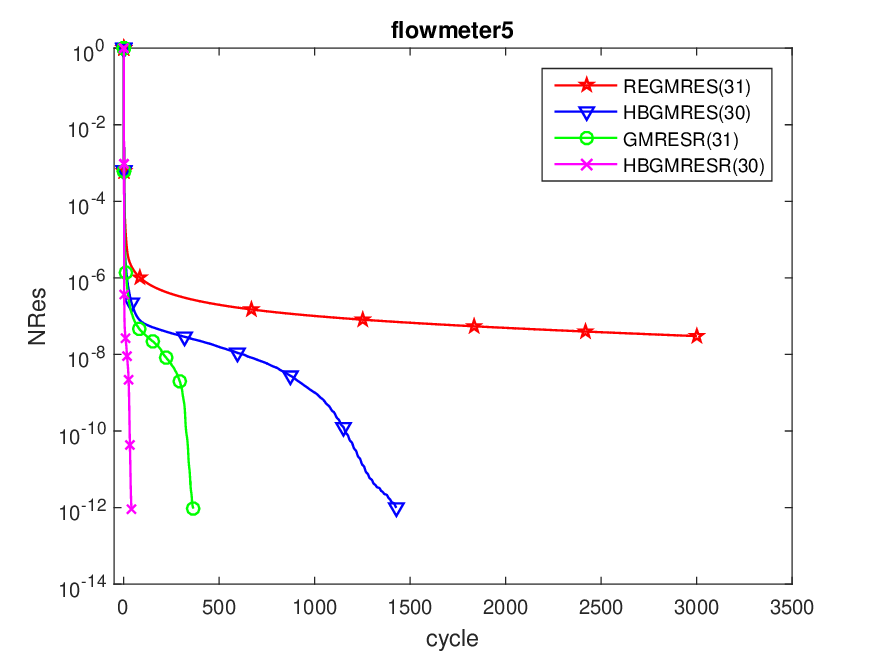}}
&  \hspace{-0.5 cm}
\resizebox*{0.48\textwidth}{0.240\textheight}{\includegraphics{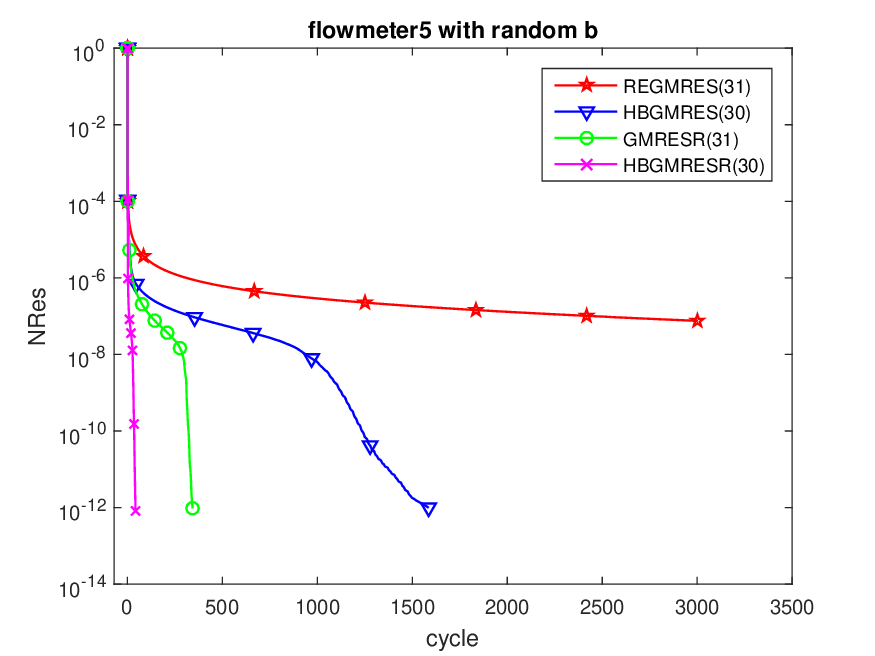}} \\
\hspace{-0.3 cm}
\resizebox*{0.48\textwidth}{0.240\textheight}{\includegraphics{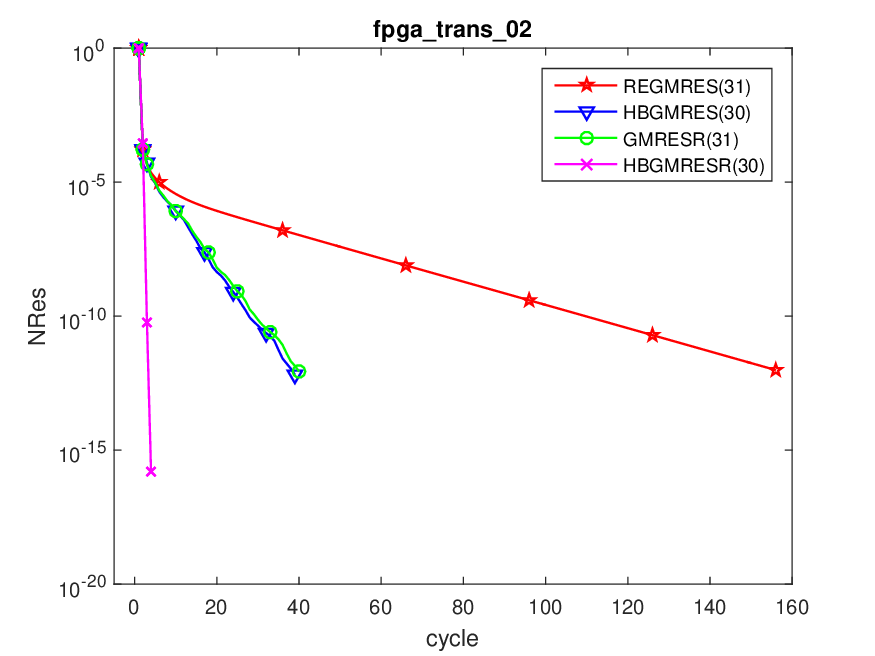}}
&  \hspace{-0.5 cm}
\resizebox*{0.48\textwidth}{0.240\textheight}{\includegraphics{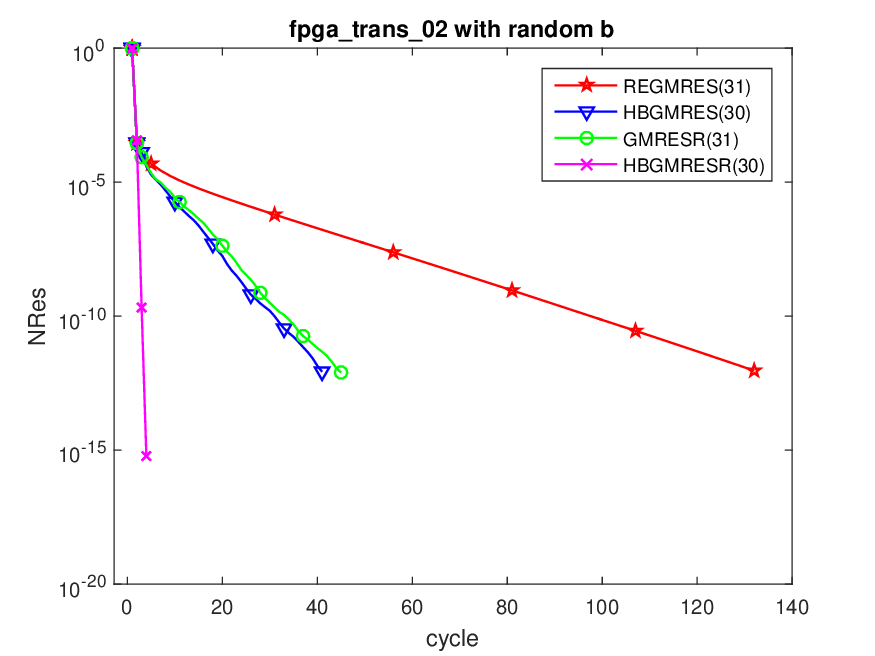}} \\
\hspace{-0.3 cm}
\resizebox*{0.48\textwidth}{0.240\textheight}{\includegraphics{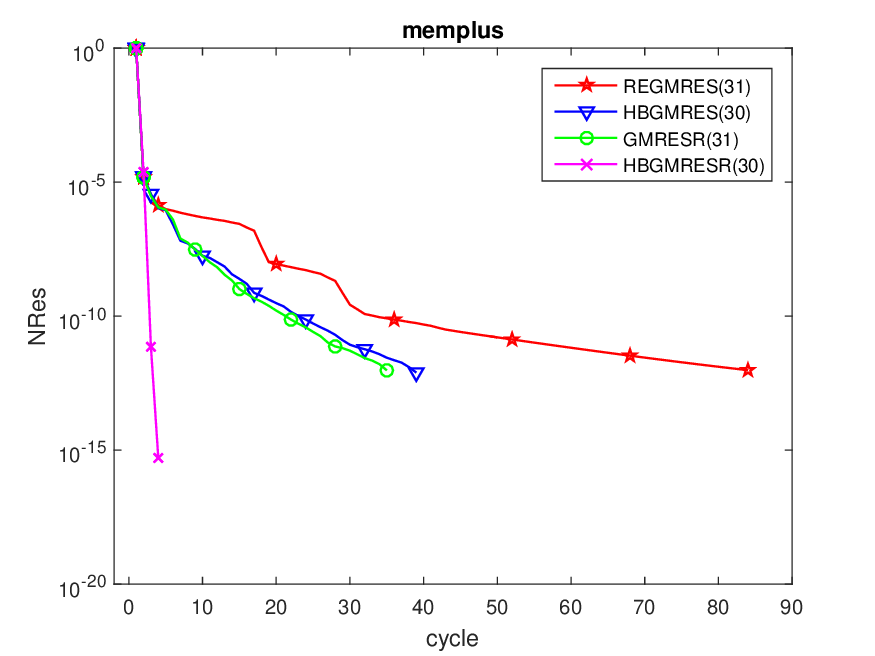}}
&  \hspace{-0.5 cm}
\resizebox*{0.48\textwidth}{0.240\textheight}{\includegraphics{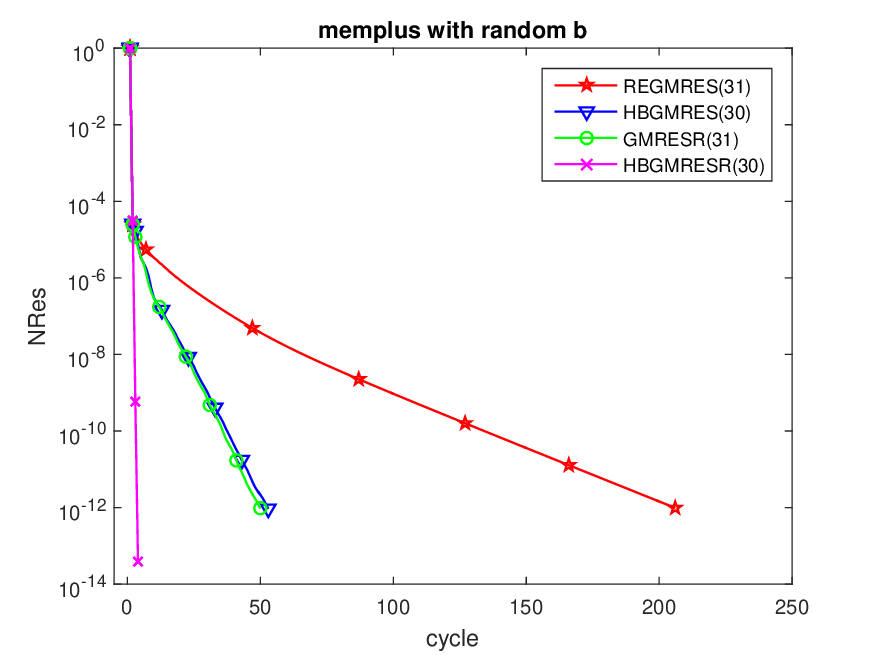}}
\end{tabular}\par
}\vspace{-0.15 cm}
\caption{\small 
    NRes {\em vs.} cycle for
   {\tt comsol}, {\tt flowmeter5}, {\tt fpga\_trans\_02}, and {\tt memplus}
   with always reorthogonalization. {\em Left:\/} original $b$; {\em Right:\/} random $b$.
   }
\label{fig:2nd4w}
\end{figure}

\newpage
\begin{figure}
{\centering
\begin{tabular}{cc}
\hspace{-0.3 cm}
\resizebox*{0.48\textwidth}{0.240\textheight}{\includegraphics{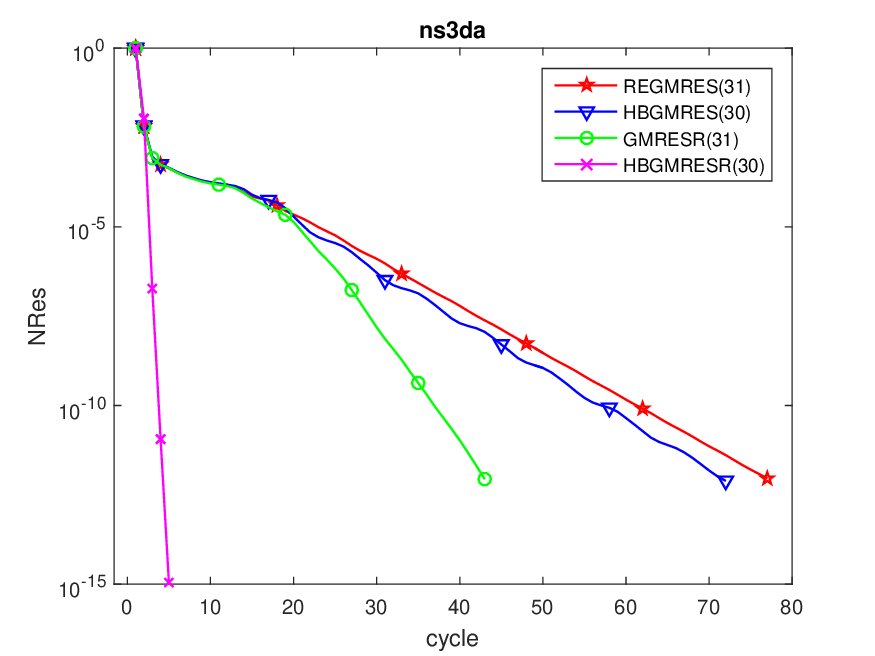}}
&  \hspace{-0.5 cm}
\resizebox*{0.48\textwidth}{0.240\textheight}{\includegraphics{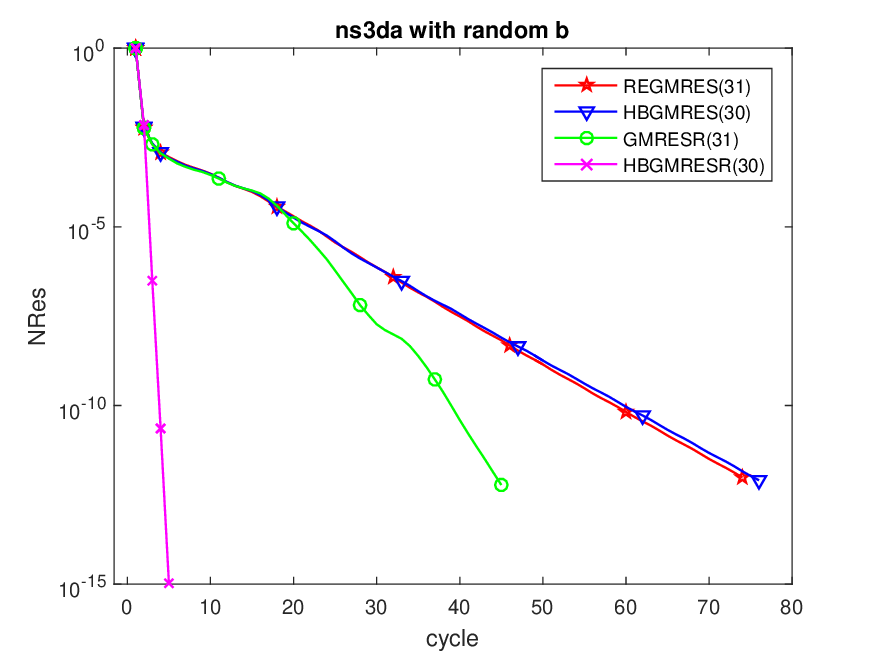}} \\
\hspace{-0.3 cm}
\resizebox*{0.48\textwidth}{0.240\textheight}{\includegraphics{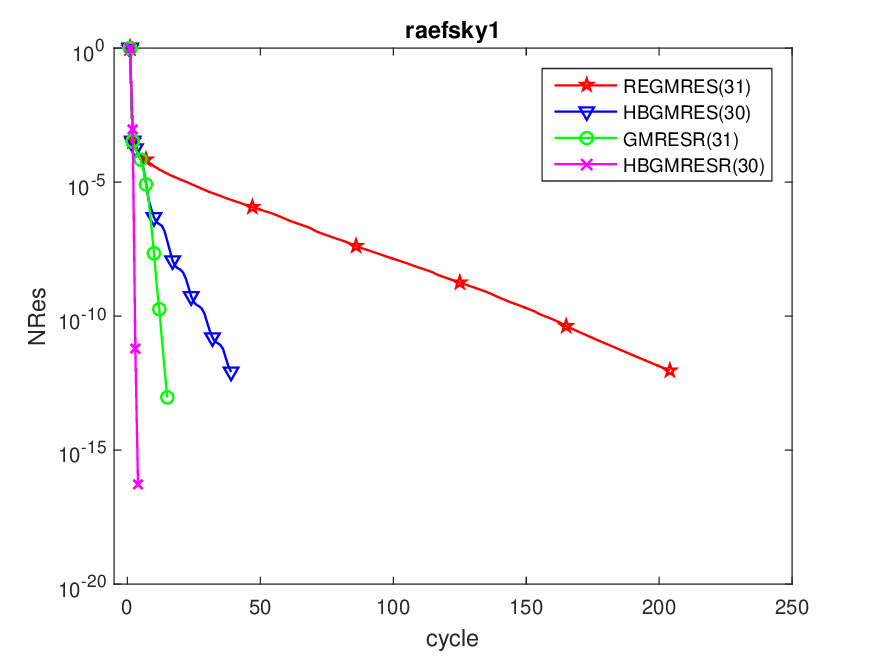}}
&  \hspace{-0.5 cm}
\resizebox*{0.48\textwidth}{0.240\textheight}{\includegraphics{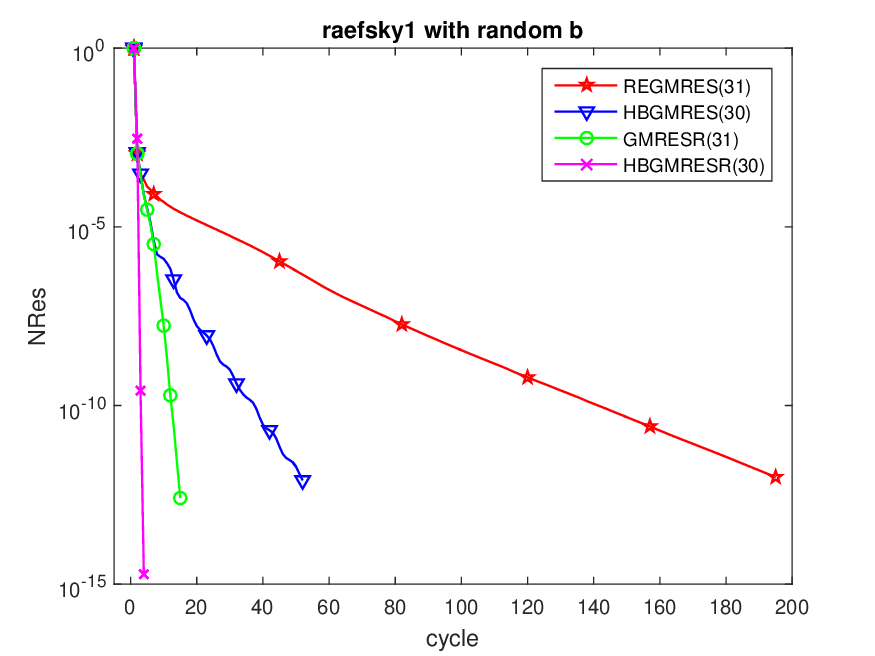}} \\
\hspace{-0.3 cm}
\resizebox*{0.48\textwidth}{0.240\textheight}{\includegraphics{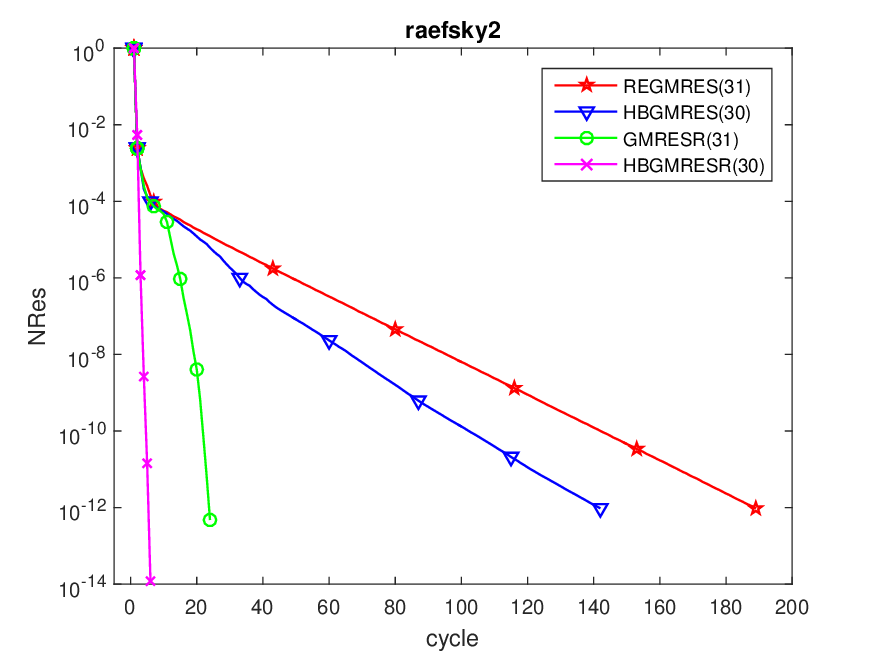}}
&  \hspace{-0.5 cm}
\resizebox*{0.48\textwidth}{0.240\textheight}{\includegraphics{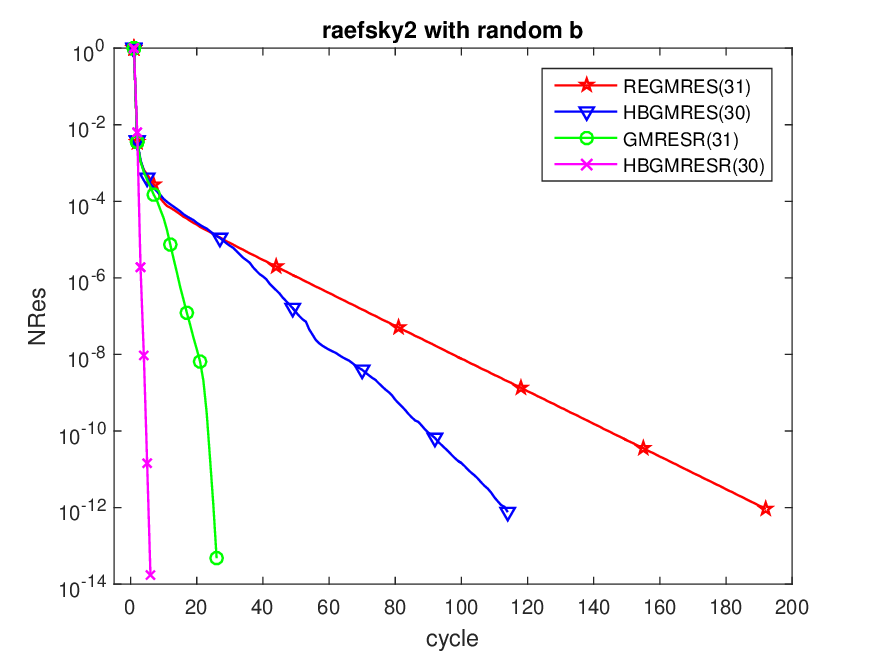}} \\
\hspace{-0.3 cm}
\resizebox*{0.48\textwidth}{0.240\textheight}{\includegraphics{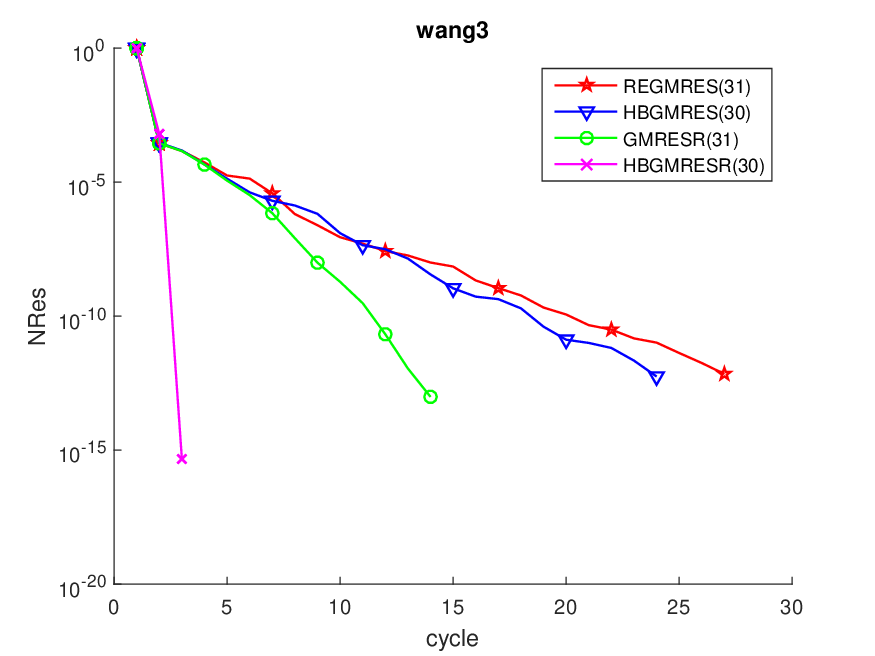}}
&  \hspace{-0.5 cm}
\resizebox*{0.48\textwidth}{0.240\textheight}{\includegraphics{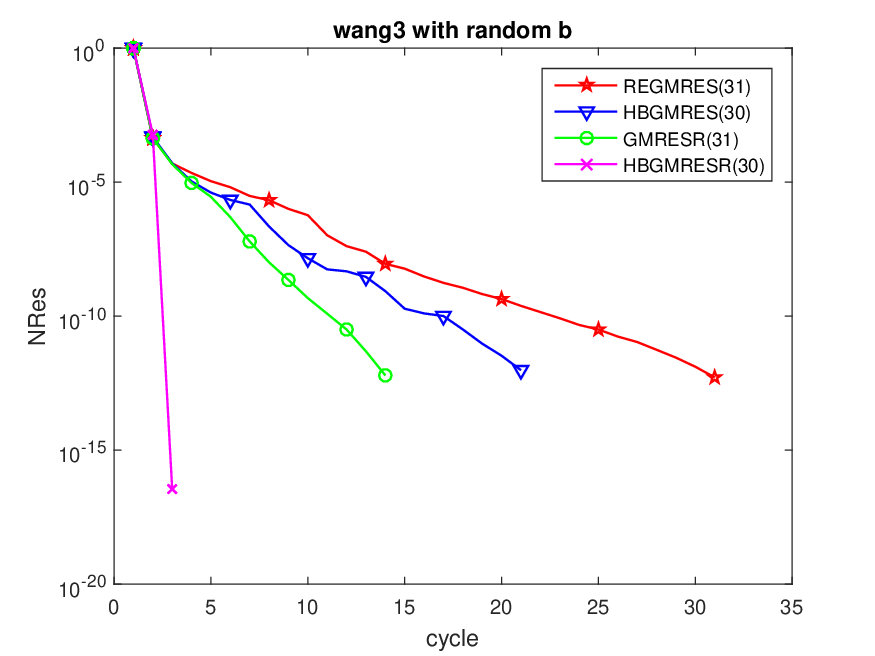}}
\end{tabular}\par
}\vspace{-0.15 cm}
\caption{\small 
    NRes {\em vs.} cycle for
   {\tt ns3da}, {\tt raefsky1}, {\tt raefsky2}, and {\tt wang3}
   with selective reorthogonalization. {\em Left:\/} original $b$; {\em Right:\/} random $b$.
   }
\label{fig:3rd4}
\end{figure}
\vspace{2mm}

\newpage
\begin{figure}
{\centering
\begin{tabular}{cc}
\hspace{-0.3 cm}
\resizebox*{0.48\textwidth}{0.240\textheight}{\includegraphics{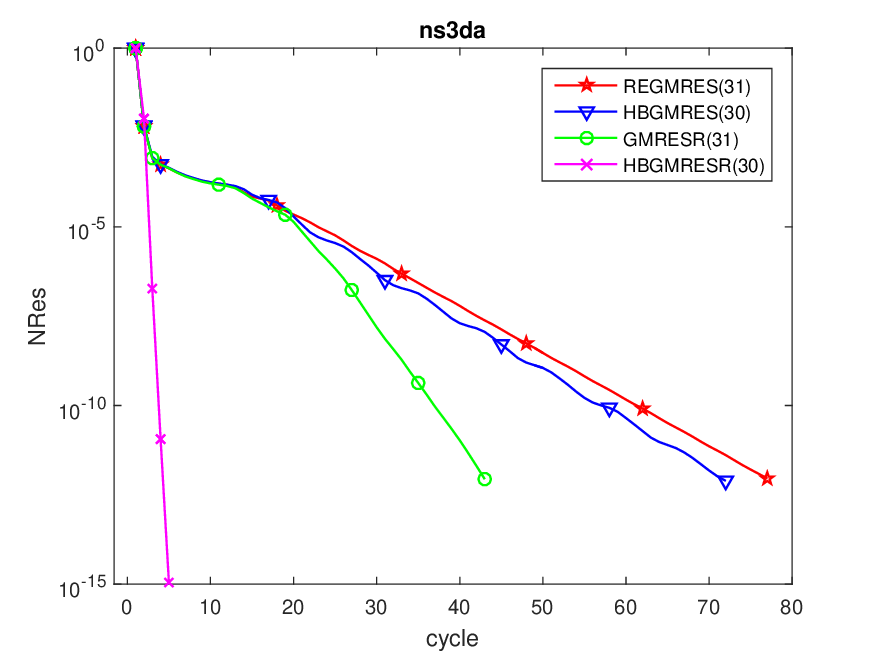}}
&  \hspace{-0.5 cm}
\resizebox*{0.48\textwidth}{0.240\textheight}{\includegraphics{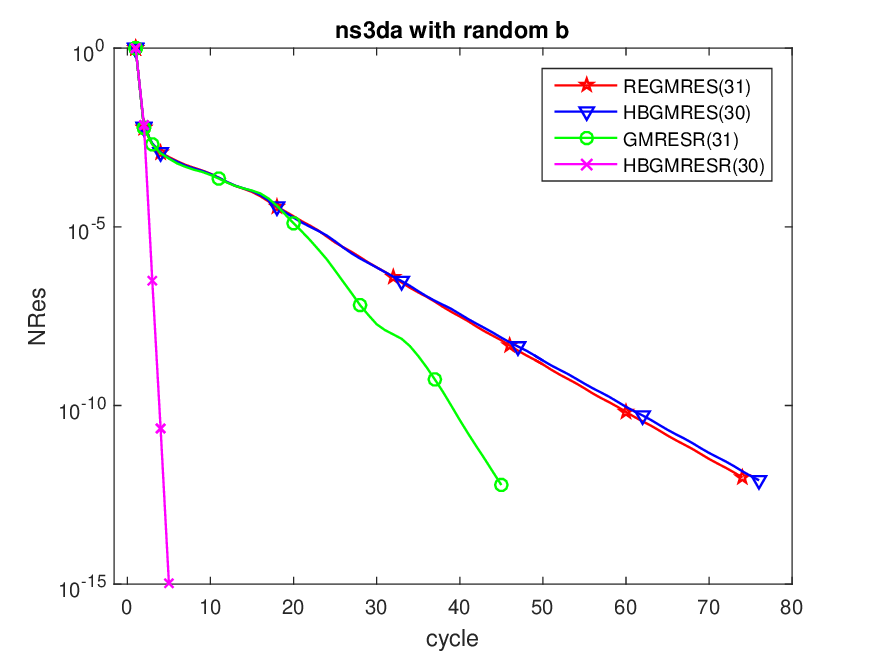}} \\
\hspace{-0.3 cm}
\resizebox*{0.48\textwidth}{0.240\textheight}{\includegraphics{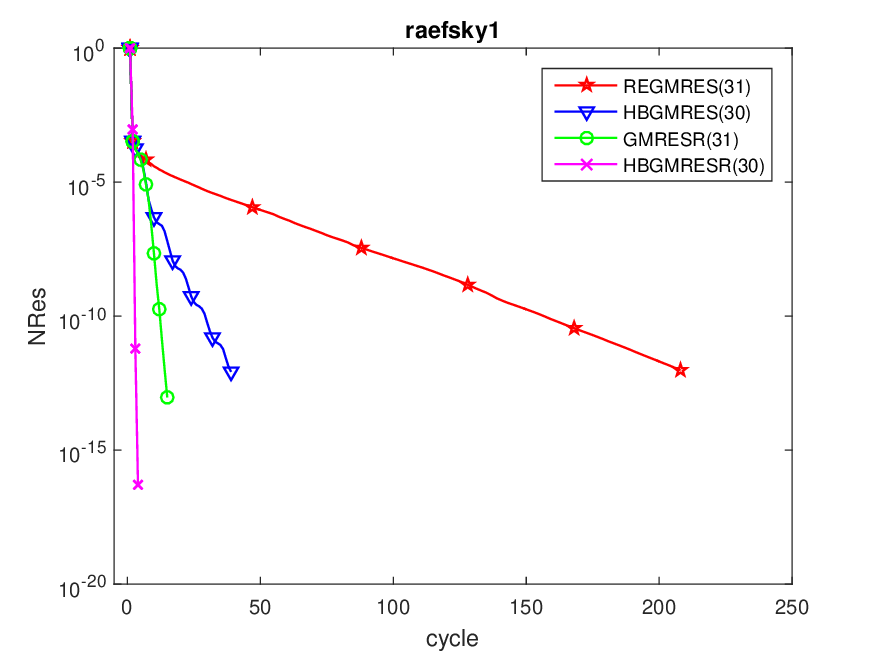}}
&  \hspace{-0.5 cm}
\resizebox*{0.48\textwidth}{0.240\textheight}{\includegraphics{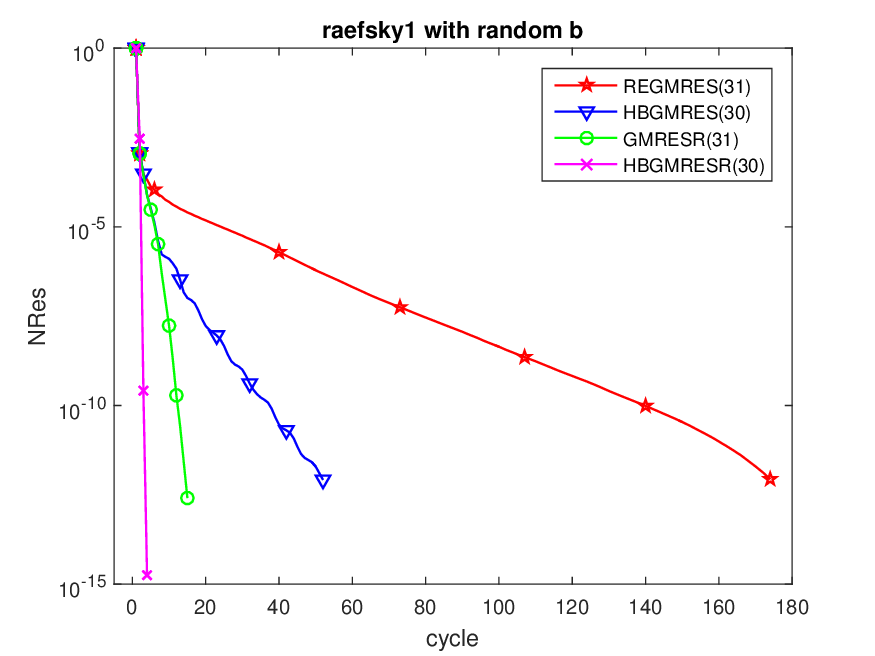}} \\
\hspace{-0.3 cm}
\resizebox*{0.48\textwidth}{0.240\textheight}{\includegraphics{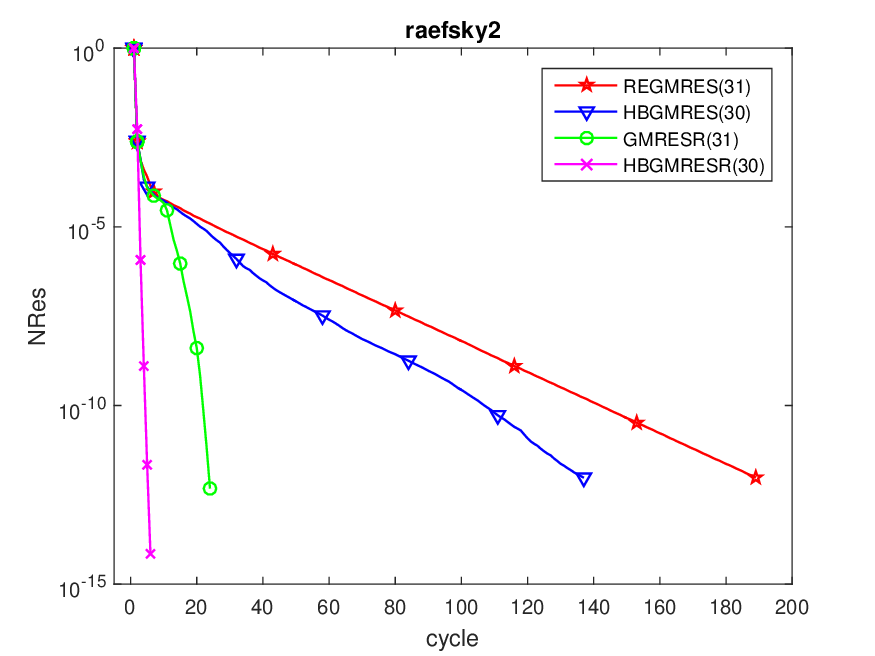}}
&  \hspace{-0.5 cm}
\resizebox*{0.48\textwidth}{0.240\textheight}{\includegraphics{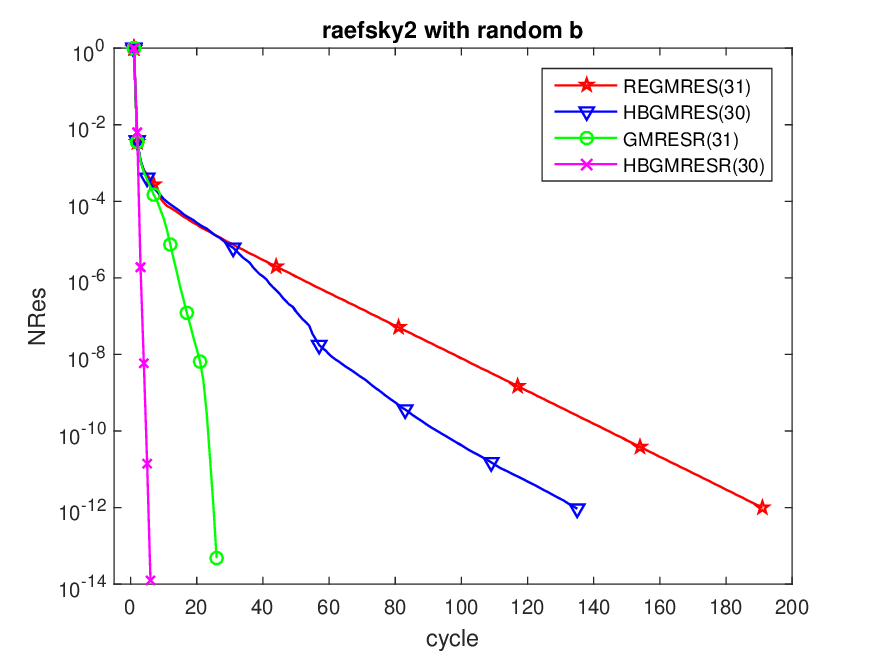}} \\
\hspace{-0.3 cm}
\resizebox*{0.48\textwidth}{0.240\textheight}{\includegraphics{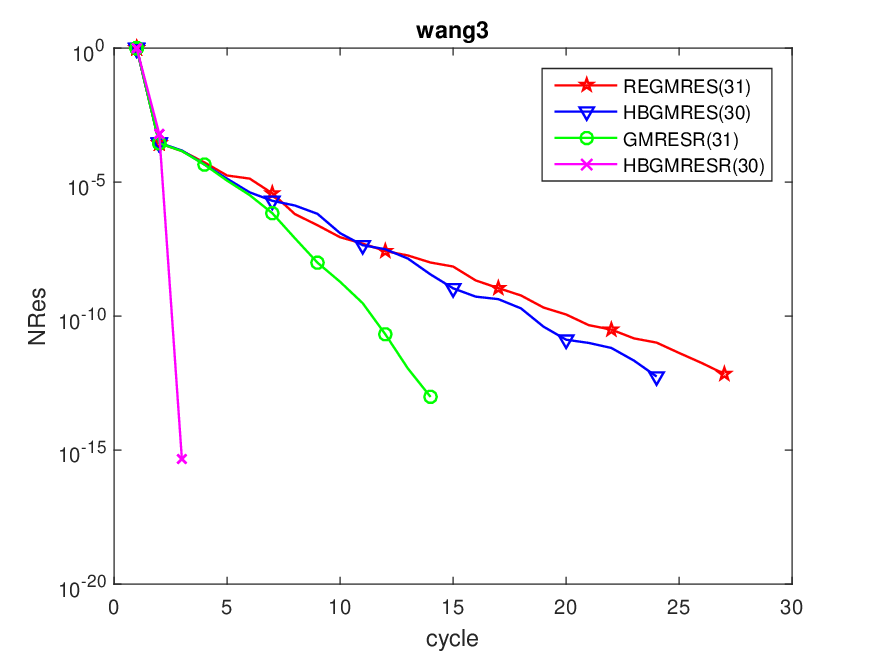}}
&  \hspace{-0.5 cm}
\resizebox*{0.48\textwidth}{0.240\textheight}{\includegraphics{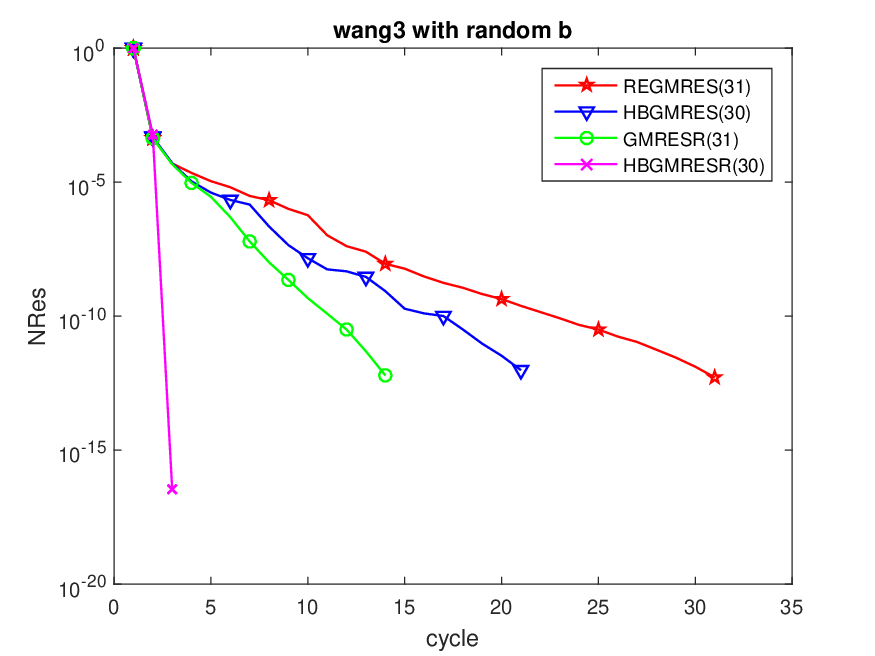}}
\end{tabular}\par
}\vspace{-0.15 cm}
\caption{\small 
    NRes {\em vs.} cycle for
   {\tt ns3da}, {\tt raefsky1}, {\tt raefsky2}, and {\tt wang3}
   with always reorthogonalization. {\em Left:\/} original $b$; {\em Right:\/} random $b$.
   }
\label{fig:3rd4w}
\end{figure}
\vspace{2mm}

\newpage
\begin{figure}
{\centering
\begin{tabular}{cc}
\hspace{-0.3 cm}
\resizebox*{0.48\textwidth}{0.240\textheight}{\includegraphics{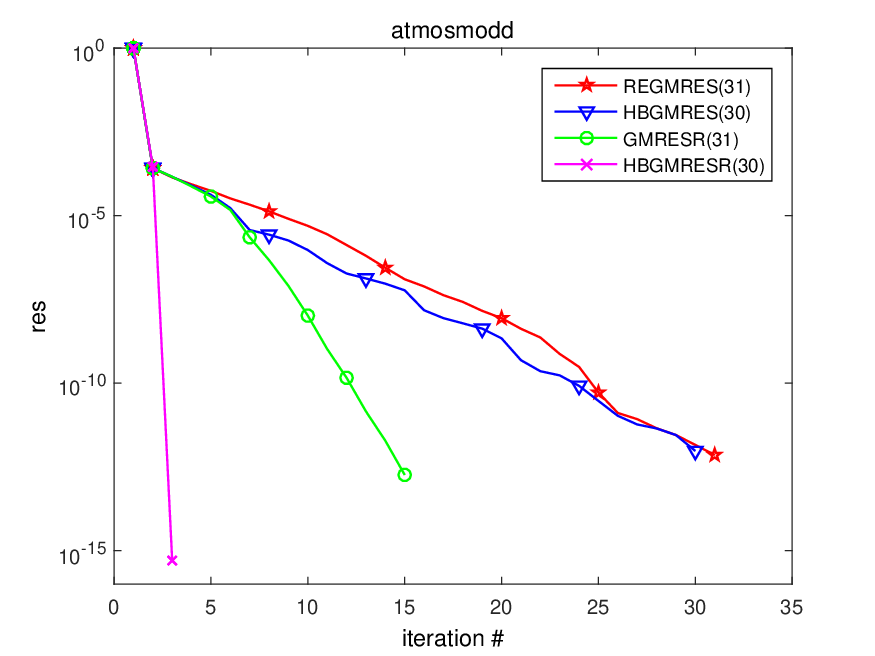}}
&  \hspace{-0.5 cm}
\resizebox*{0.48\textwidth}{0.240\textheight}{\includegraphics{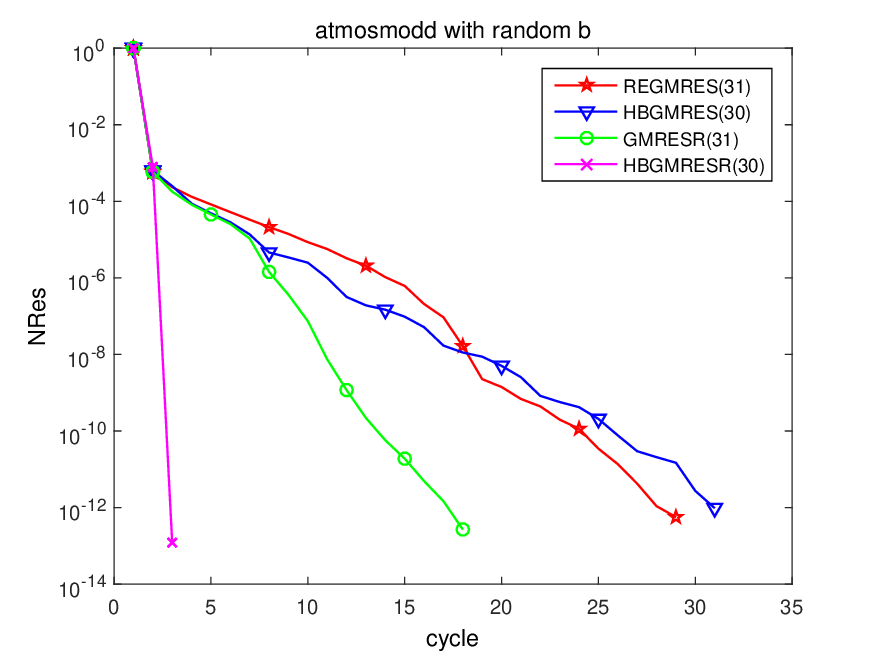}} \\
\hspace{-0.3 cm}
\resizebox*{0.48\textwidth}{0.240\textheight}{\includegraphics{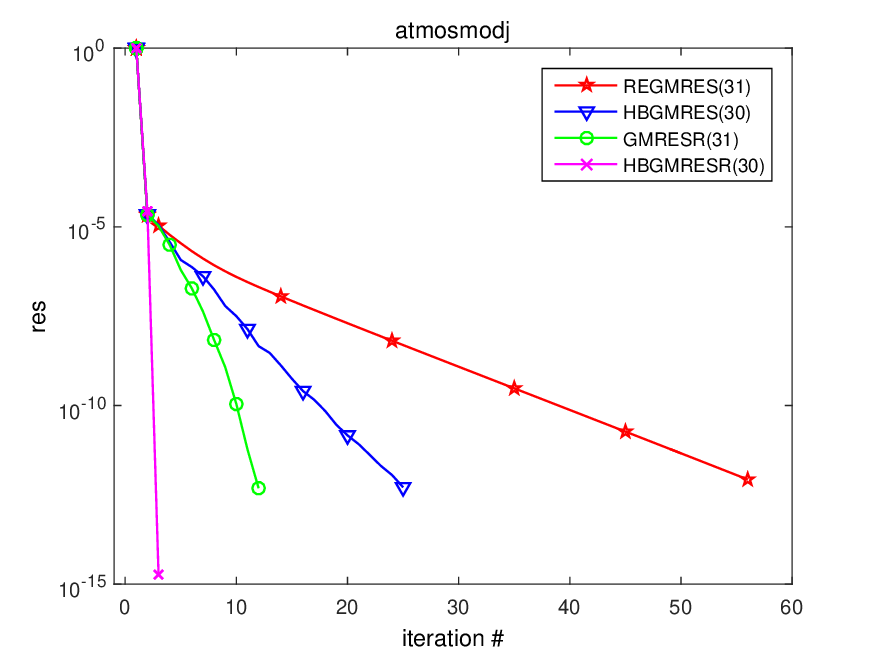}}
&  \hspace{-0.5 cm}
\resizebox*{0.48\textwidth}{0.240\textheight}{\includegraphics{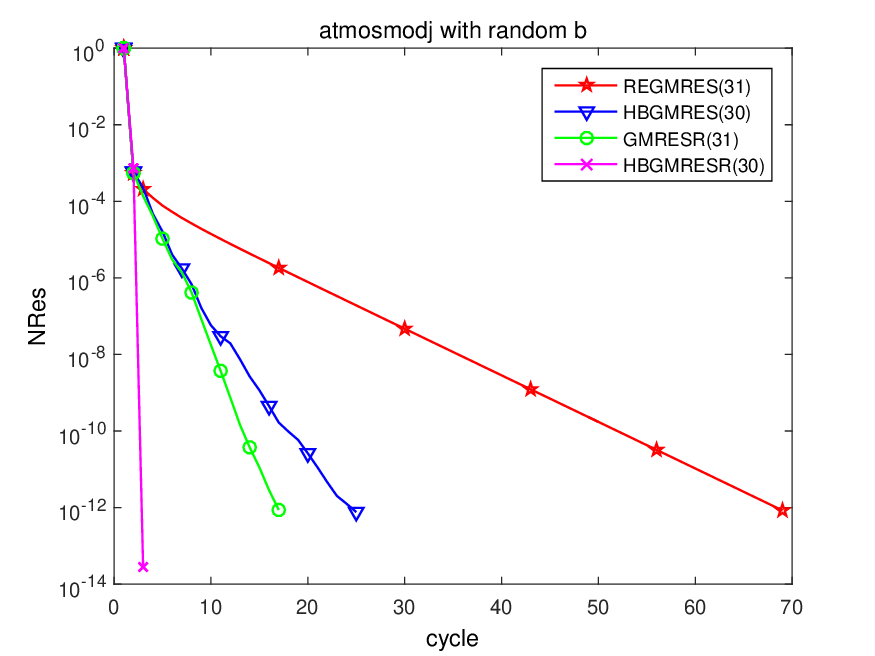}} \\
\hspace{-0.3 cm}
\resizebox*{0.48\textwidth}{0.240\textheight}{\includegraphics{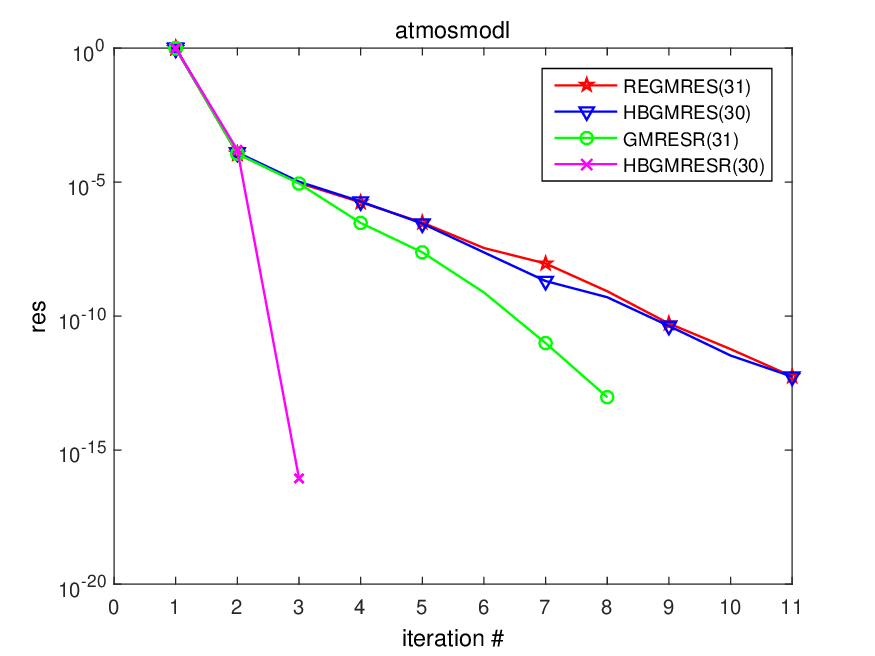}}
&  \hspace{-0.5 cm}
\resizebox*{0.48\textwidth}{0.240\textheight}{\includegraphics{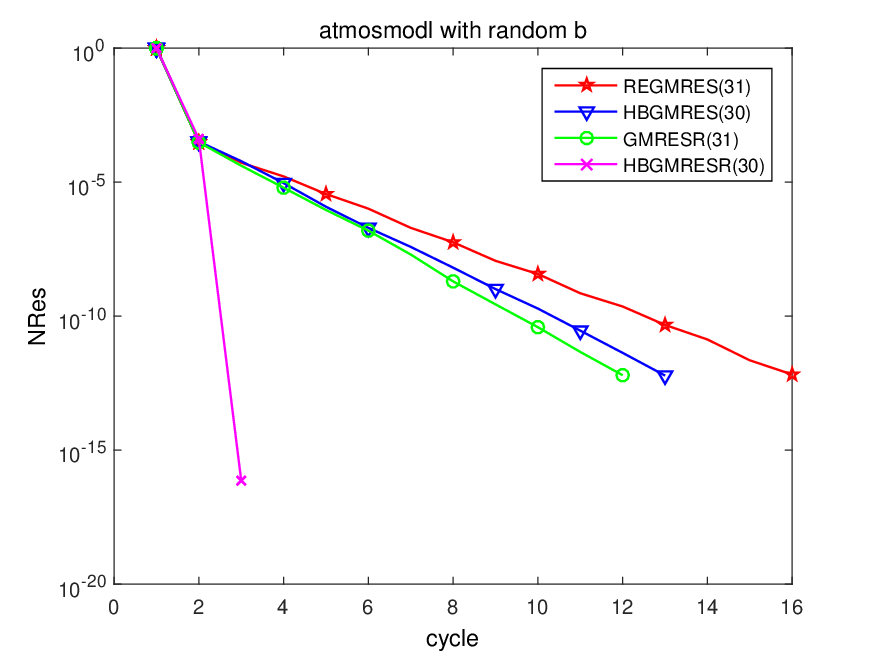}} \\
\hspace{-0.3 cm}
\resizebox*{0.48\textwidth}{0.240\textheight}{\includegraphics{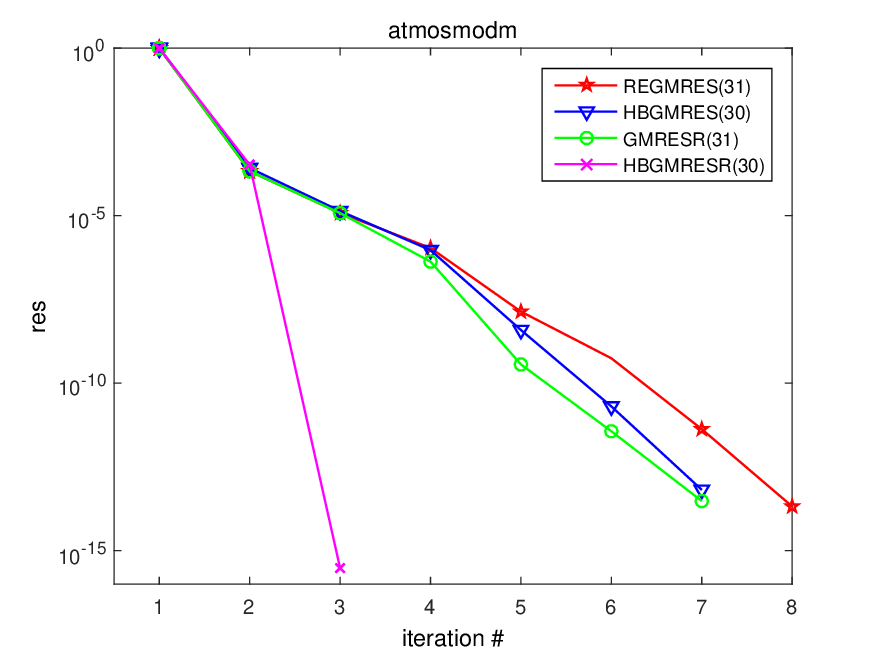}}
&  \hspace{-0.5 cm}
\resizebox*{0.48\textwidth}{0.240\textheight}{\includegraphics{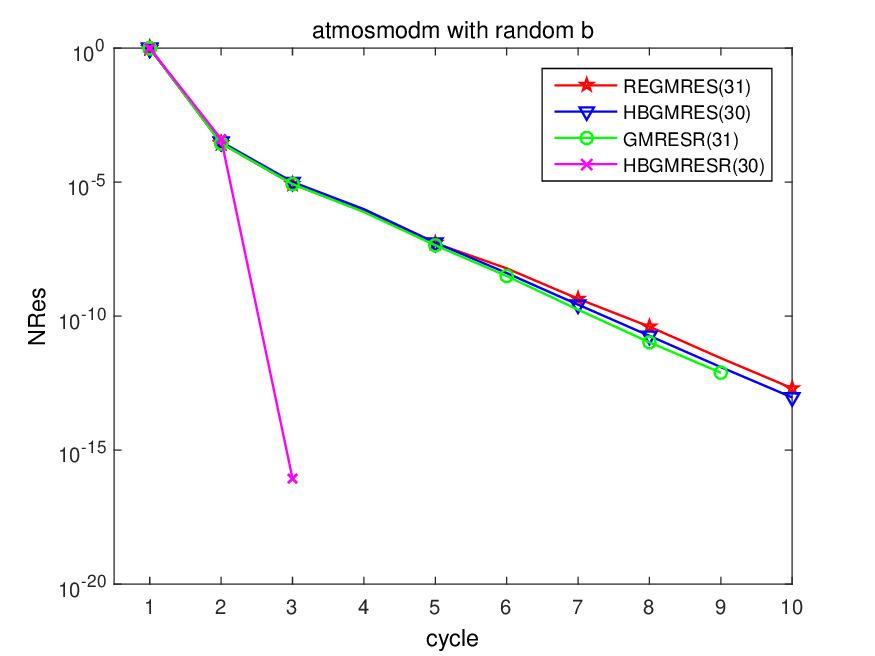}}
\end{tabular}\par
}\vspace{-0.15 cm}
\caption{\small 
    NRes {\em vs.} cycle for
   {\tt atmosmodd}, {\tt atmosmodj}, {\tt atmosmodl}, and {\tt atmosmodm}
   with selective reorthogonalization. {\em Left:\/} original $b$; {\em Right:\/} random $b$.
   }
\label{fig:4th4}
\end{figure}

\newpage
\begin{figure}
{\centering
\begin{tabular}{cc}
\hspace{-0.3 cm}
\resizebox*{0.48\textwidth}{0.240\textheight}{\includegraphics{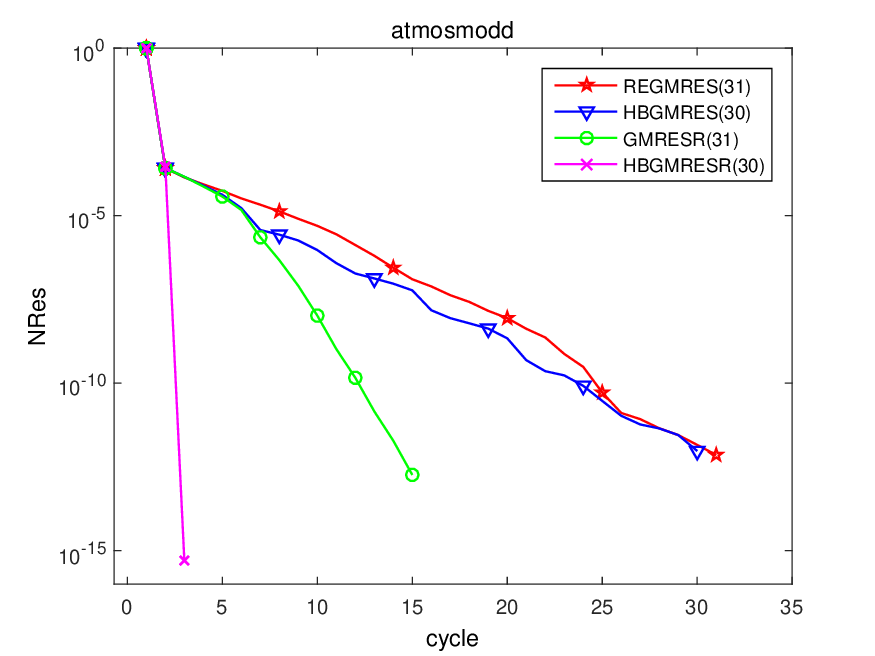}}
&  \hspace{-0.5 cm}
\resizebox*{0.48\textwidth}{0.240\textheight}{\includegraphics{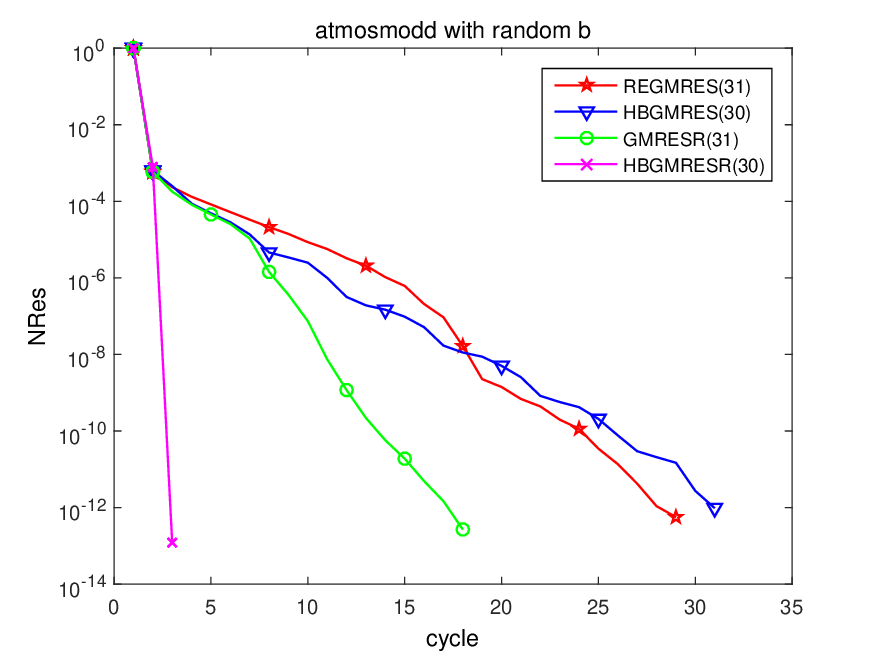}} \\
\hspace{-0.3 cm}
\resizebox*{0.48\textwidth}{0.240\textheight}{\includegraphics{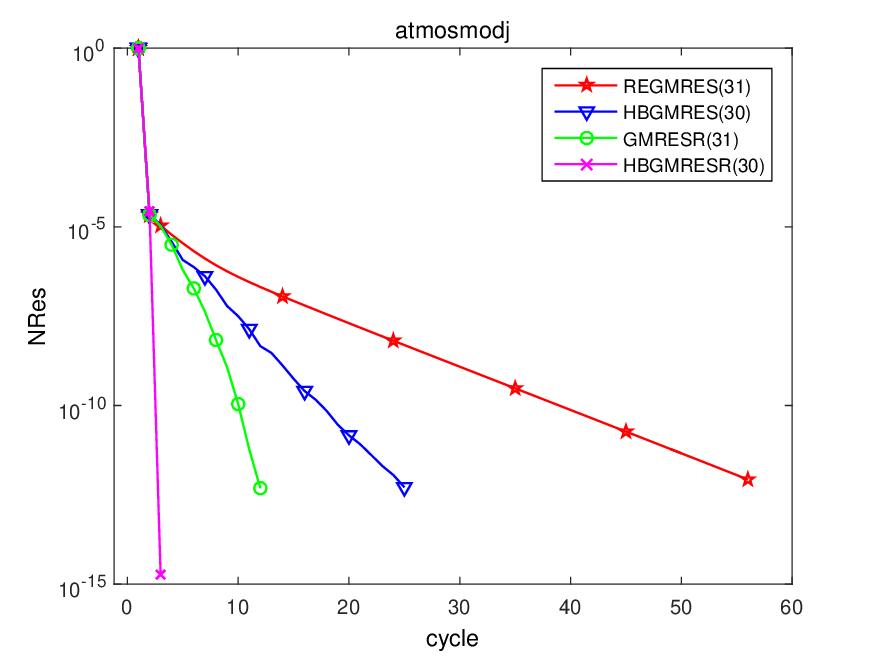}}
&  \hspace{-0.5 cm}
\resizebox*{0.48\textwidth}{0.240\textheight}{\includegraphics{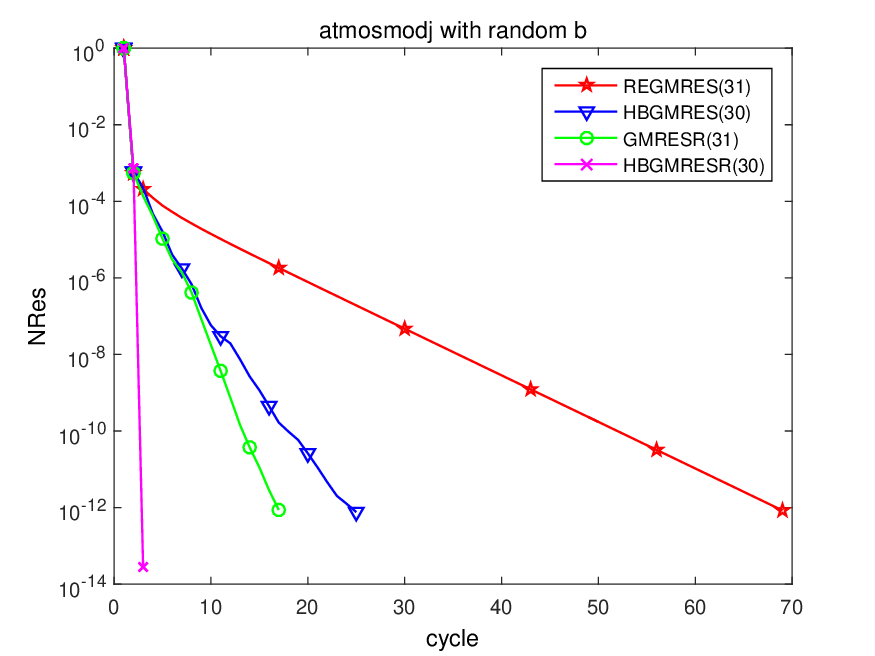}} \\
\hspace{-0.3 cm}
\resizebox*{0.48\textwidth}{0.240\textheight}{\includegraphics{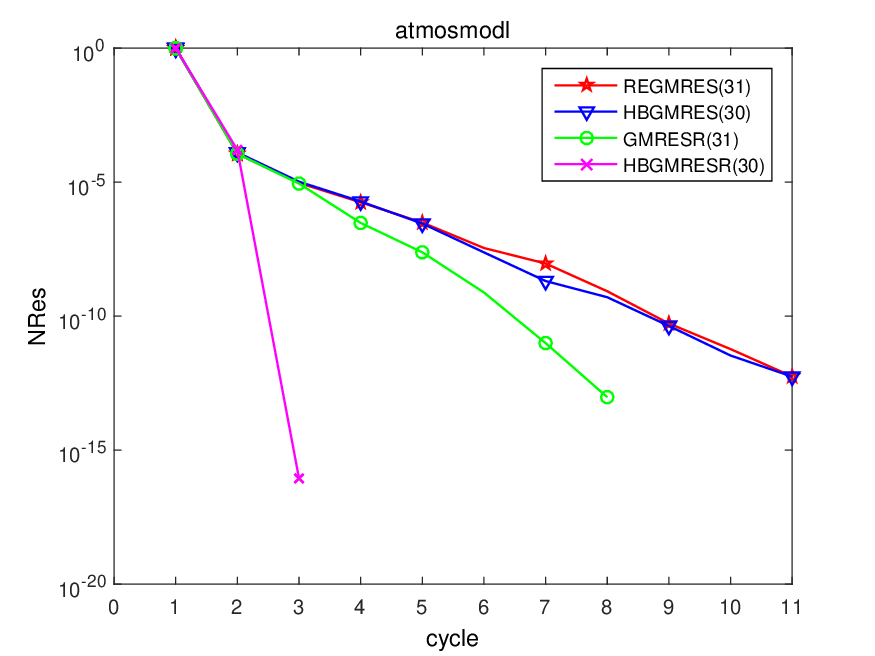}}
&  \hspace{-0.5 cm}
\resizebox*{0.48\textwidth}{0.240\textheight}{\includegraphics{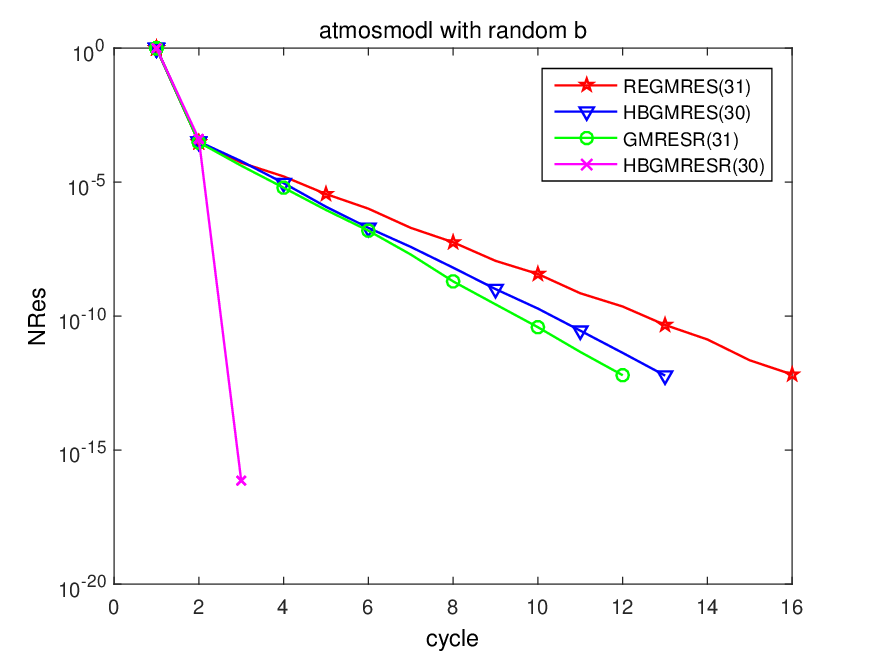}} \\
\hspace{-0.3 cm}
\resizebox*{0.48\textwidth}{0.240\textheight}{\includegraphics{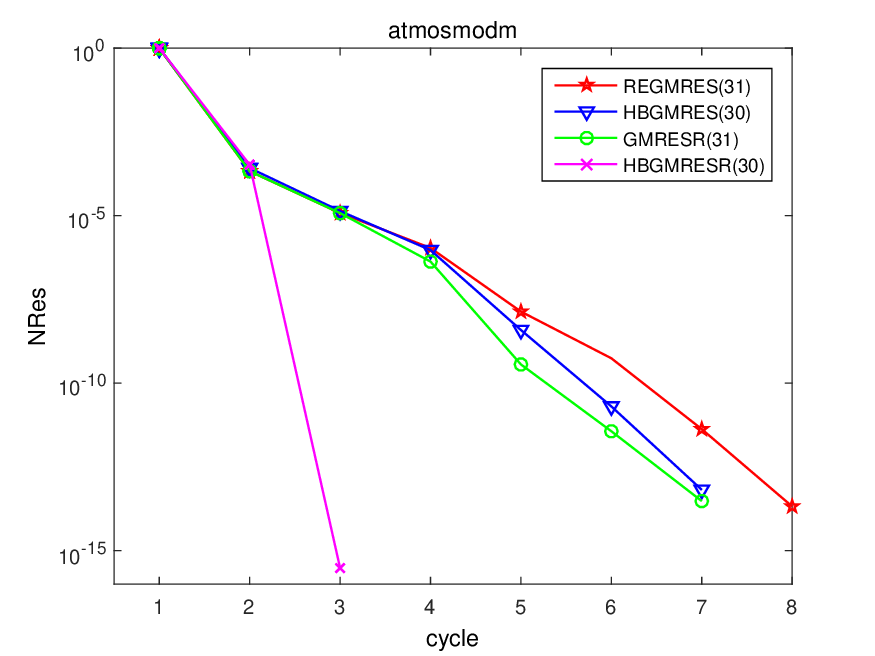}}
&  \hspace{-0.5 cm}
\resizebox*{0.48\textwidth}{0.240\textheight}{\includegraphics{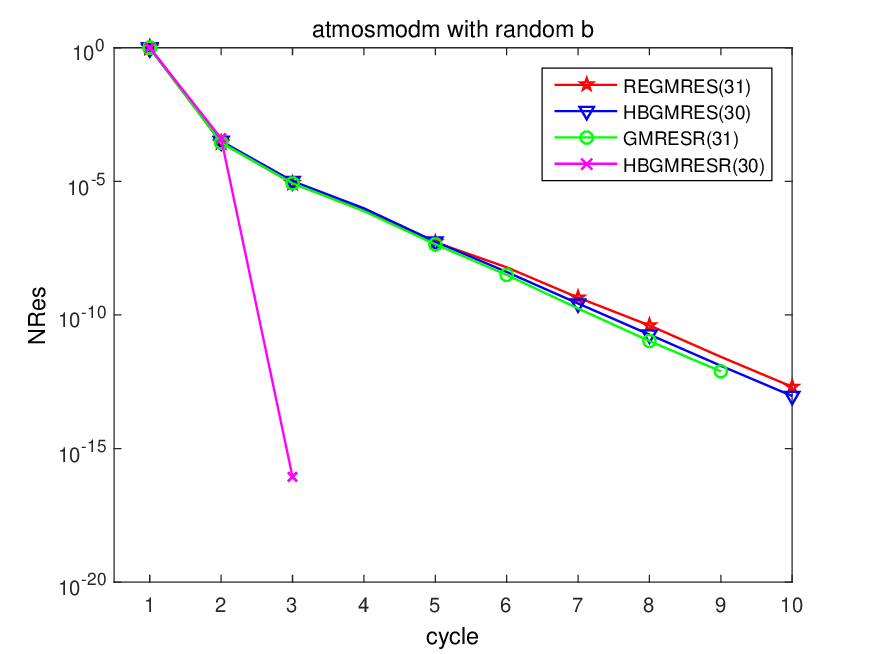}}
\end{tabular}\par
}\vspace{-0.15 cm}
\caption{\small 
    NRes {\em vs.} cycle for
   {\tt atmosmodd}, {\tt atmosmodj}, {\tt atmosmodl}, and {\tt atmosmodm}
   with selective reorthogonalization. {\em Left:\/} original $b$; {\em Right:\/} random $b$.
   }
\label{fig:4th4w}
\end{figure}

\vspace{2mm}
\newpage
\begin{figure}
{\centering
\begin{tabular}{cc}
\hspace{-0.3 cm}
\resizebox*{0.48\textwidth}{0.240\textheight}{\includegraphics{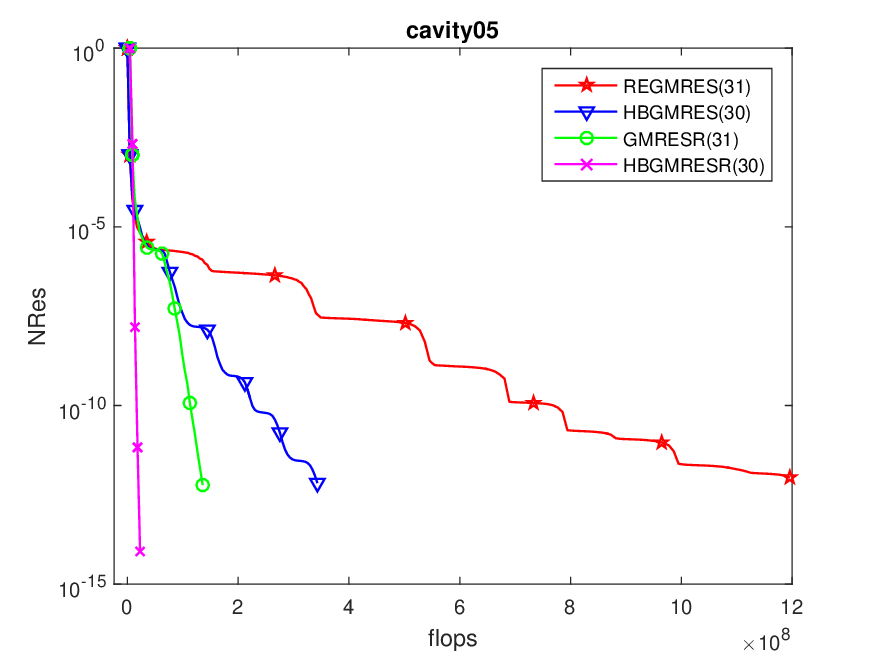}}
&  \hspace{-0.5 cm}
\resizebox*{0.48\textwidth}{0.240\textheight}{\includegraphics{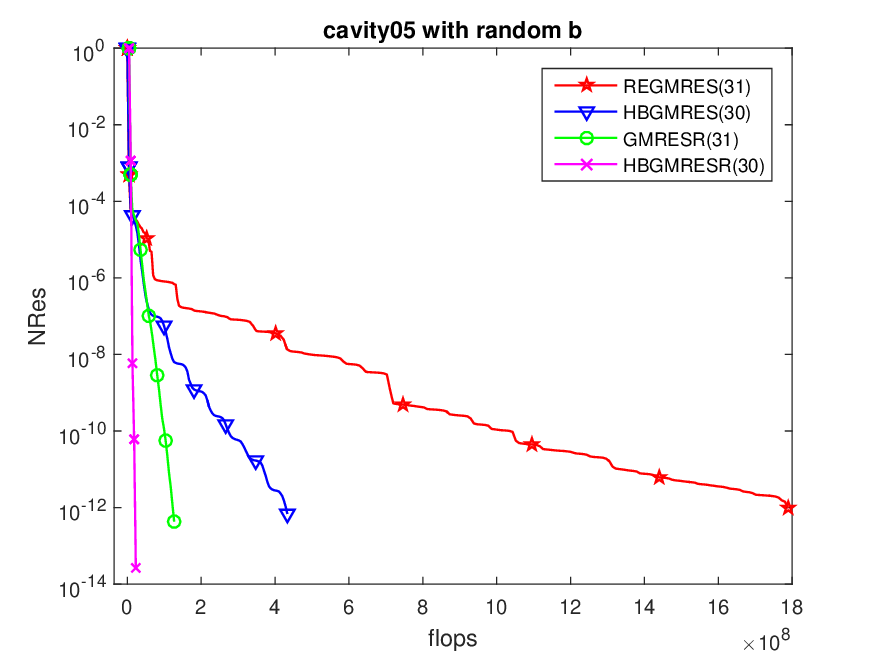}} \\
\hspace{-0.3 cm}
\resizebox*{0.48\textwidth}{0.240\textheight}{\includegraphics{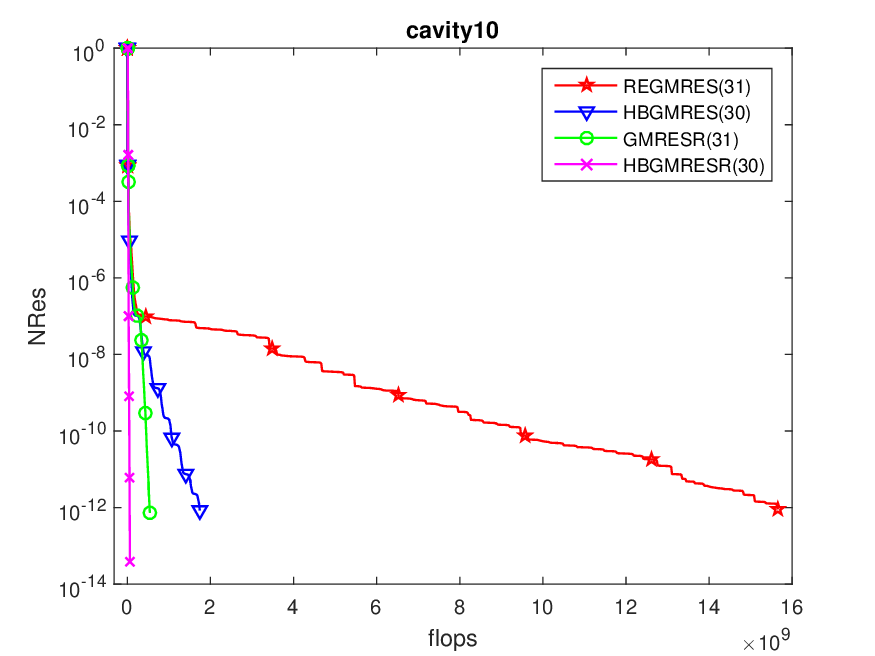}}
&  \hspace{-0.5 cm}
\resizebox*{0.48\textwidth}{0.240\textheight}{\includegraphics{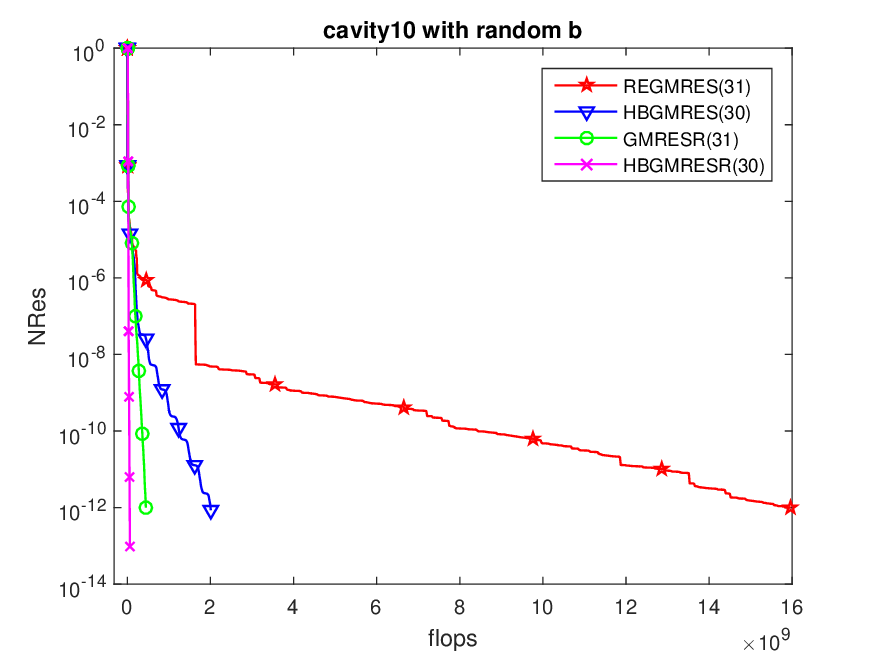}} \\
\hspace{-0.3 cm}
\resizebox*{0.48\textwidth}{0.240\textheight}{\includegraphics{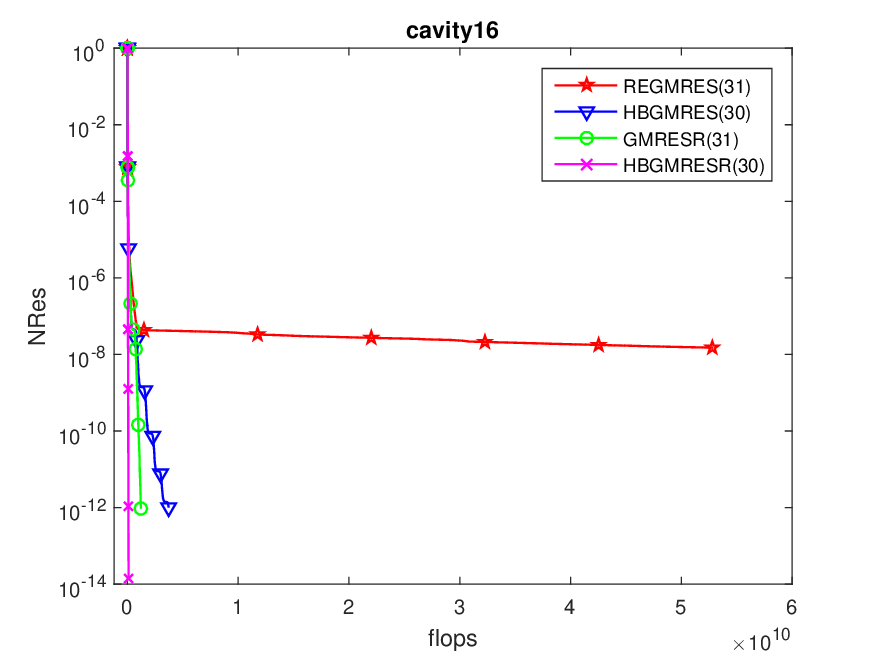}}
&  \hspace{-0.5 cm}
\resizebox*{0.48\textwidth}{0.240\textheight}{\includegraphics{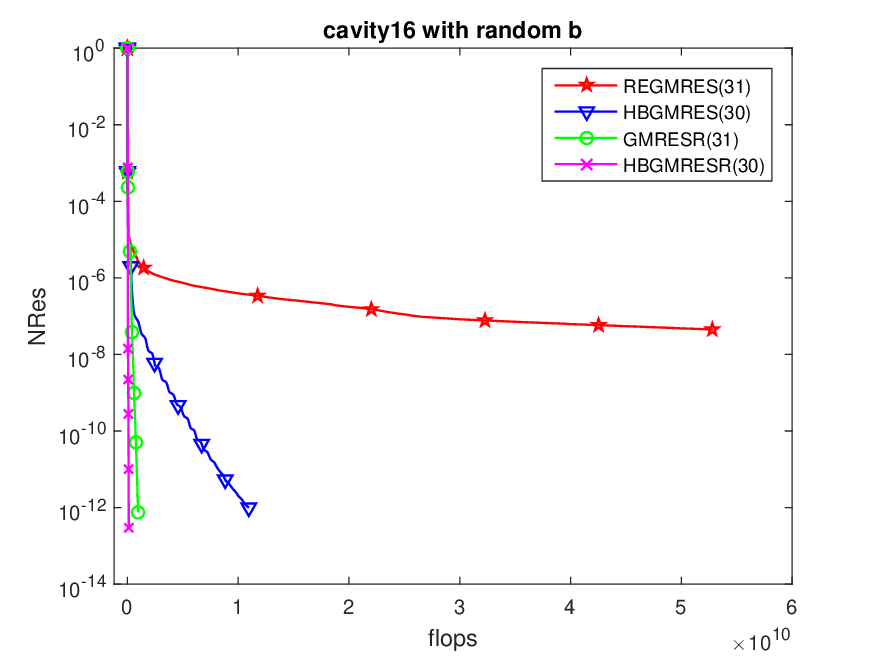}} \\
\hspace{-0.3 cm}
\resizebox*{0.48\textwidth}{0.240\textheight}{\includegraphics{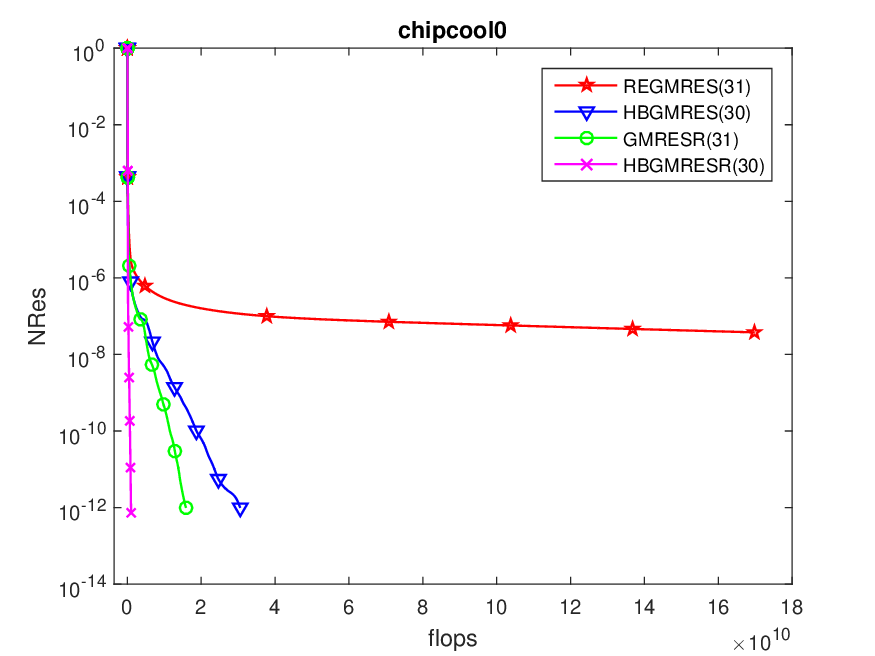}}
&  \hspace{-0.5 cm}
\resizebox*{0.48\textwidth}{0.240\textheight}{\includegraphics{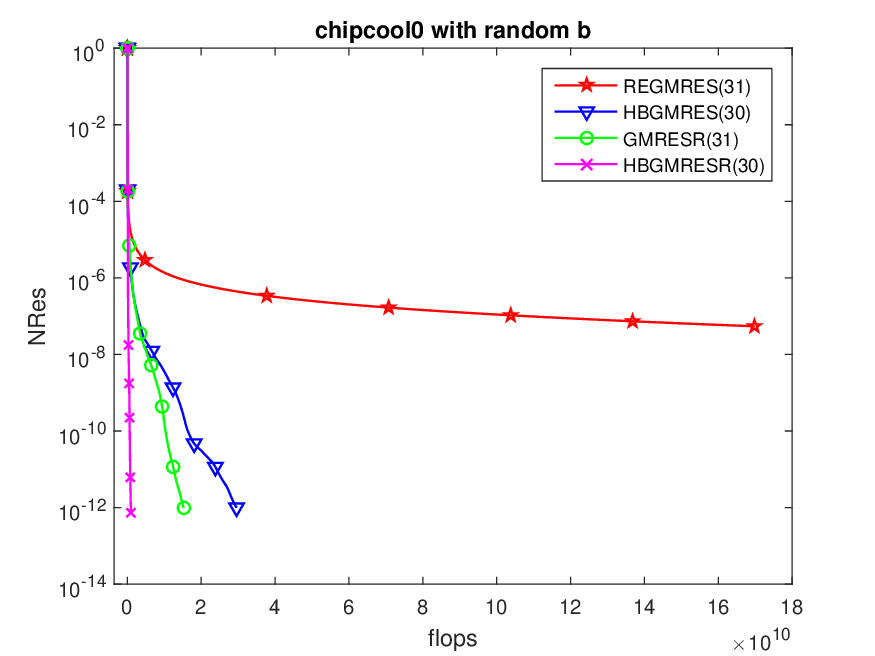}}
\end{tabular}\par
}\vspace{-0.15 cm}
\caption{\small 
    NRes {\em vs.} cycle for
   {\tt cavity05}, {\tt cavity10}, {\tt cavity16}, and {\tt chipcool0}
   with selective reorthogonalization. {\em Left:\/} original $b$; {\em Right:\/} random $b$.
   }
\label{fig:1st4flops}
\end{figure}
\vspace{2mm}
\newpage
\begin{figure}
{\centering
\begin{tabular}{cc}
\hspace{-0.3 cm}
\resizebox*{0.48\textwidth}{0.240\textheight}{\includegraphics{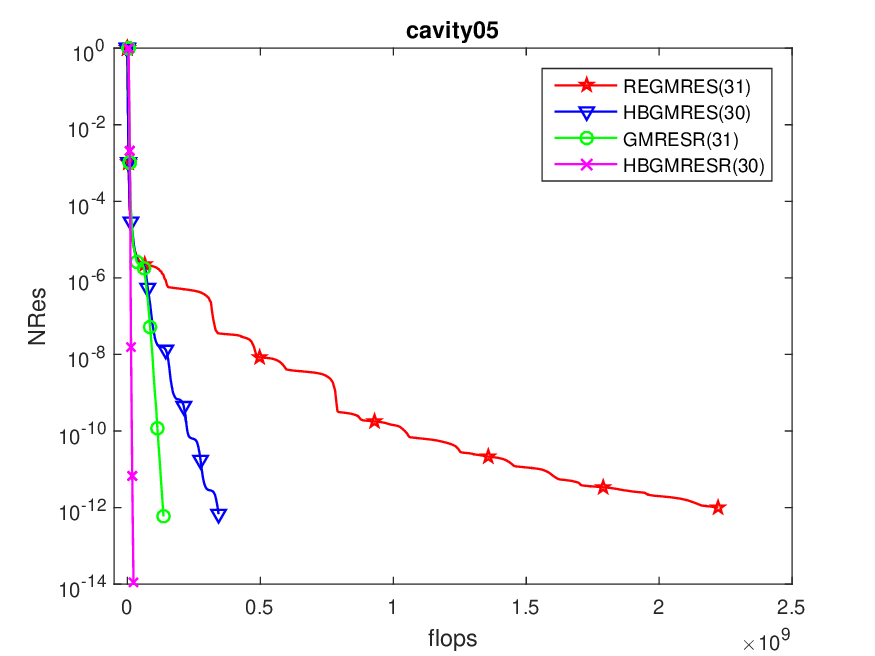}}
&  \hspace{-0.5 cm}
\resizebox*{0.48\textwidth}{0.240\textheight}{\includegraphics{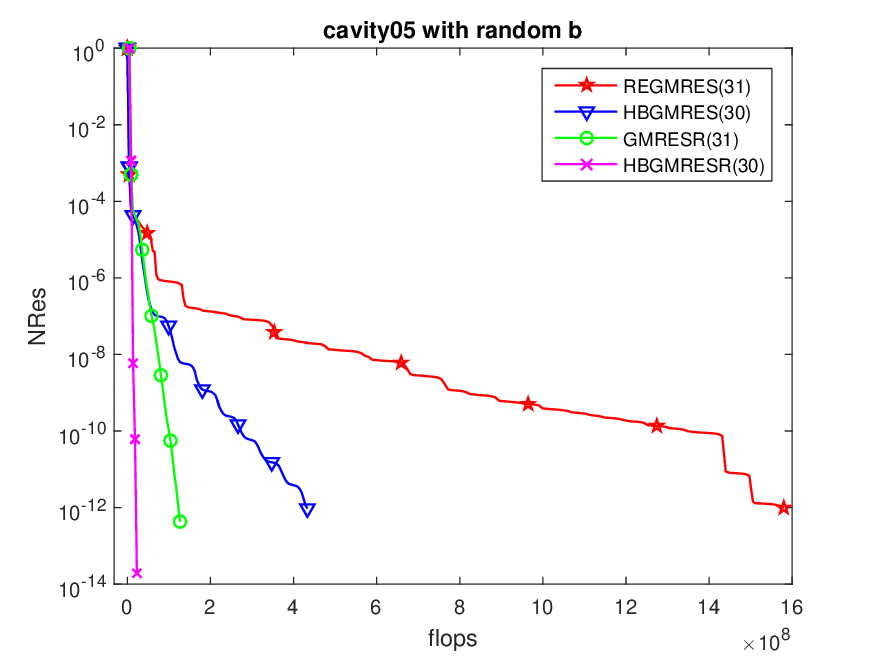}} \\
\hspace{-0.3 cm}
\resizebox*{0.48\textwidth}{0.240\textheight}{\includegraphics{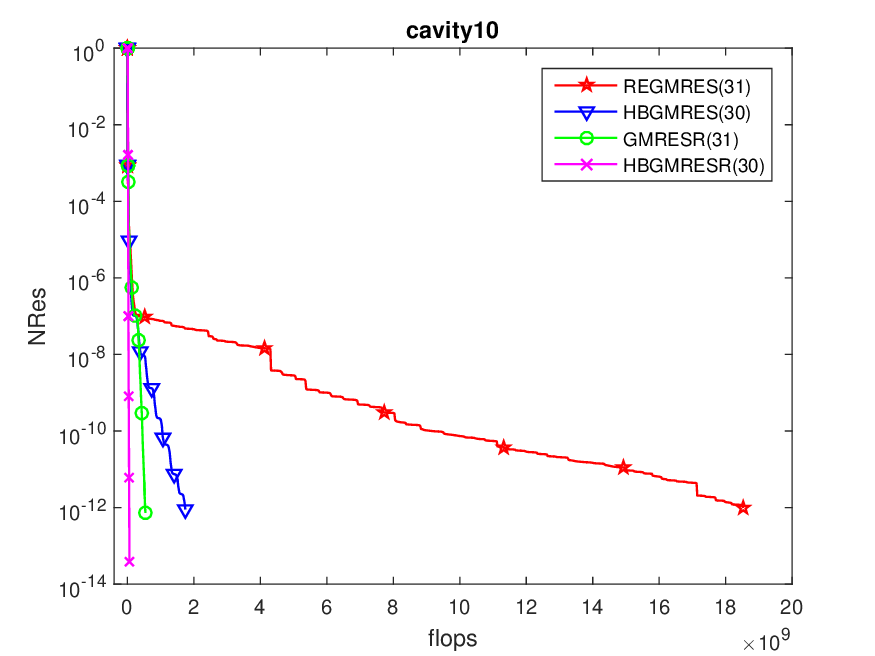}}
&  \hspace{-0.5 cm}
\resizebox*{0.48\textwidth}{0.240\textheight}{\includegraphics{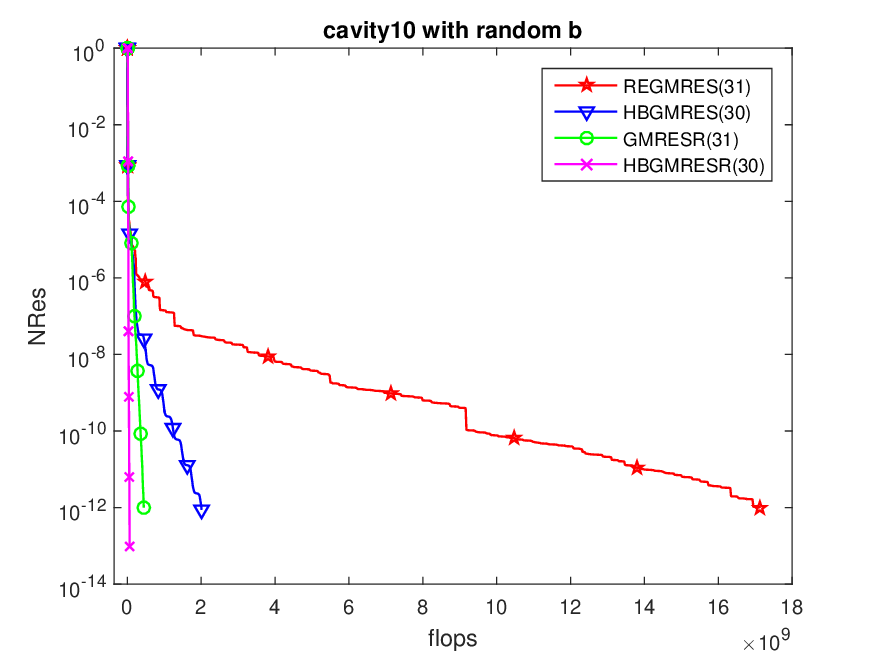}} \\
\hspace{-0.3 cm}
\resizebox*{0.48\textwidth}{0.240\textheight}{\includegraphics{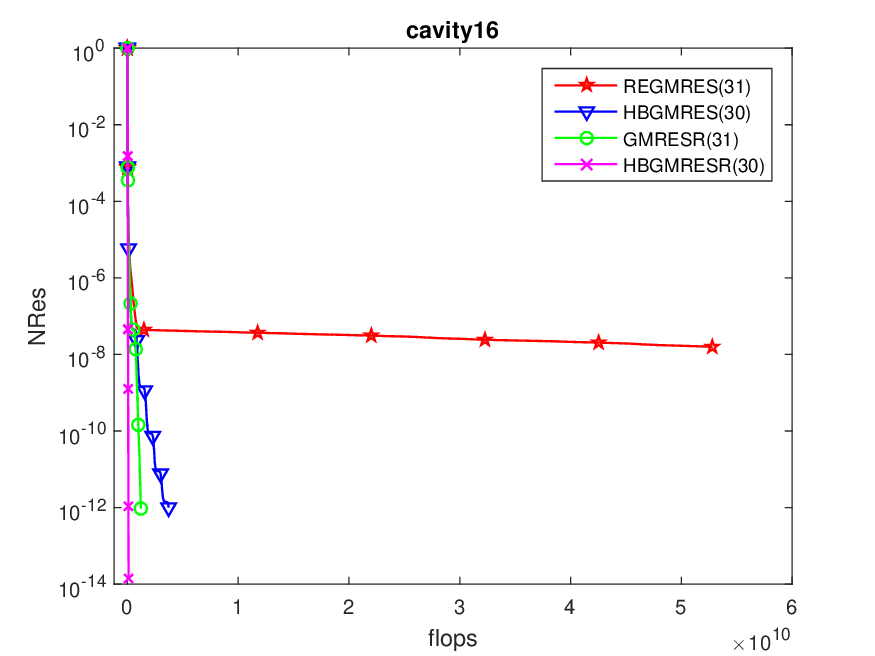}}
&  \hspace{-0.5 cm}
\resizebox*{0.48\textwidth}{0.240\textheight}{\includegraphics{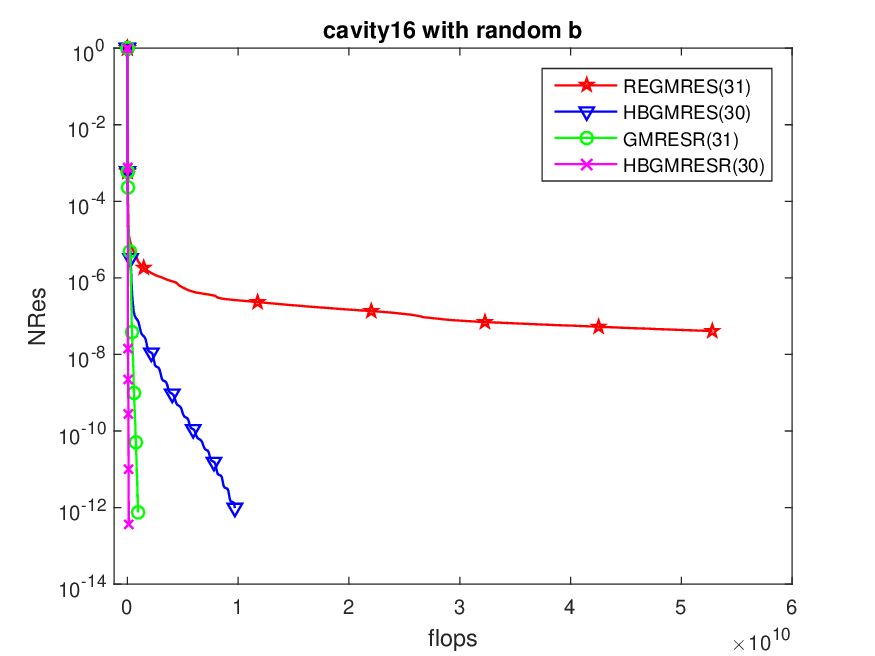}} \\
\hspace{-0.3 cm}
\resizebox*{0.48\textwidth}{0.240\textheight}{\includegraphics{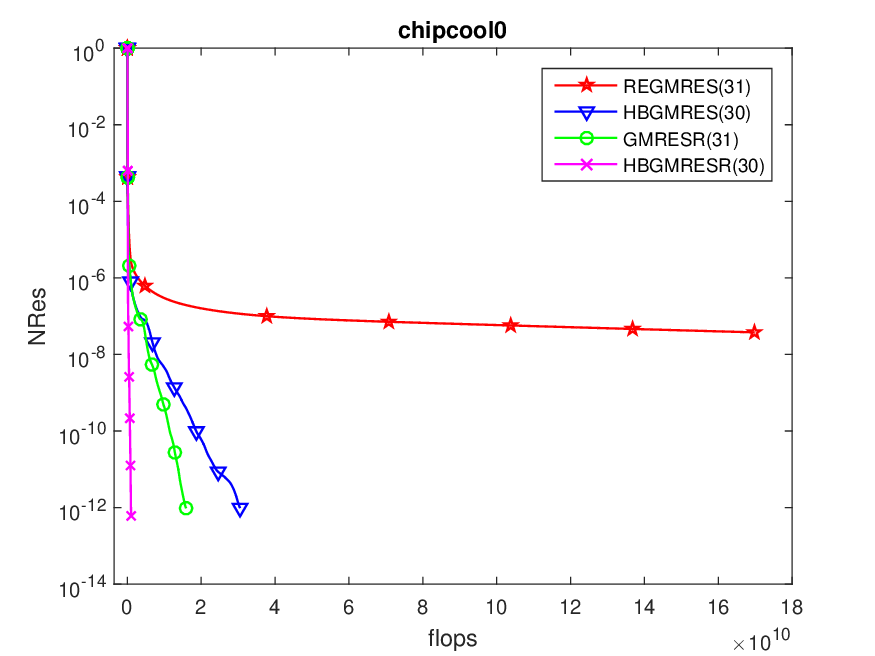}}
&  \hspace{-0.5 cm}
\resizebox*{0.48\textwidth}{0.240\textheight}{\includegraphics{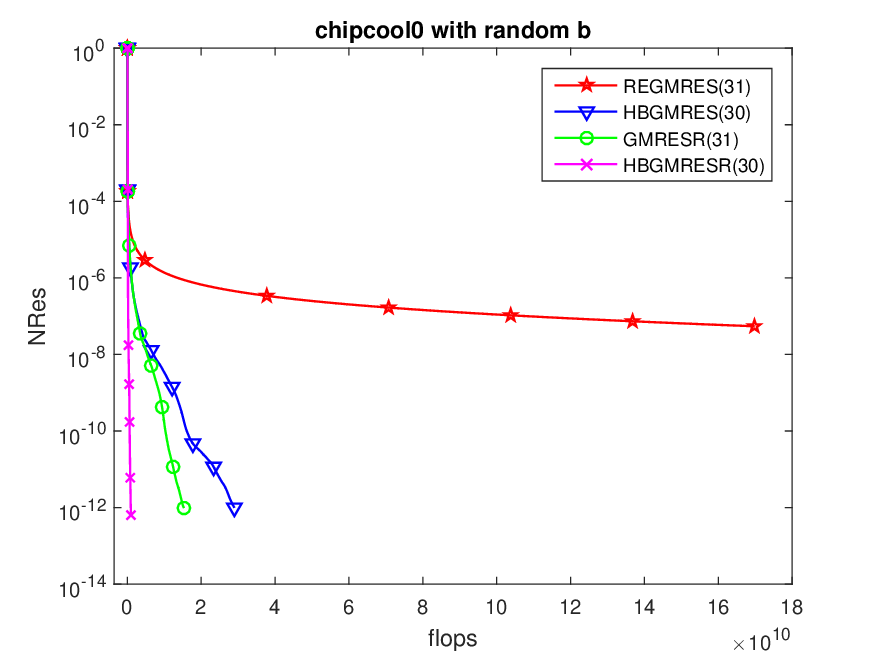}}
\end{tabular}\par
}\vspace{-0.15 cm}
\caption{\small 
    NRes {\em vs.} cycle for
   {\tt cavity05}, {\tt cavity10}, {\tt cavity16}, and {\tt chipcool0}
   with always reorthogonalization. {\em Left:\/} original $b$; {\em Right:\/} random $b$.
   }
\label{fig:1st4wflops}
\end{figure}
\vspace{2mm}

\newpage
\begin{figure}
{\centering
\begin{tabular}{cc}
\hspace{-0.3 cm}
\resizebox*{0.48\textwidth}{0.240\textheight}{\includegraphics{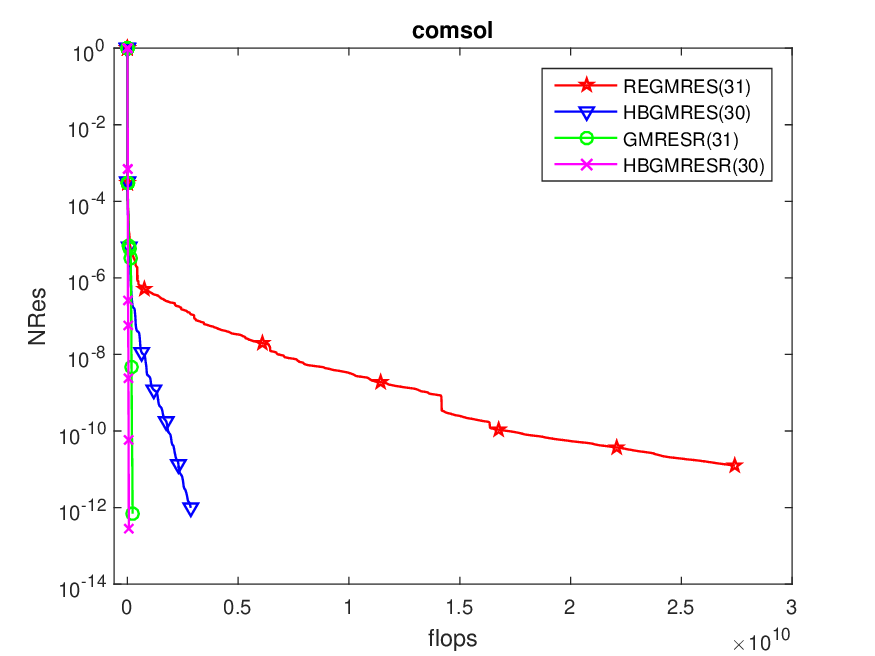}}
&  \hspace{-0.5 cm}
\resizebox*{0.48\textwidth}{0.240\textheight}{\includegraphics{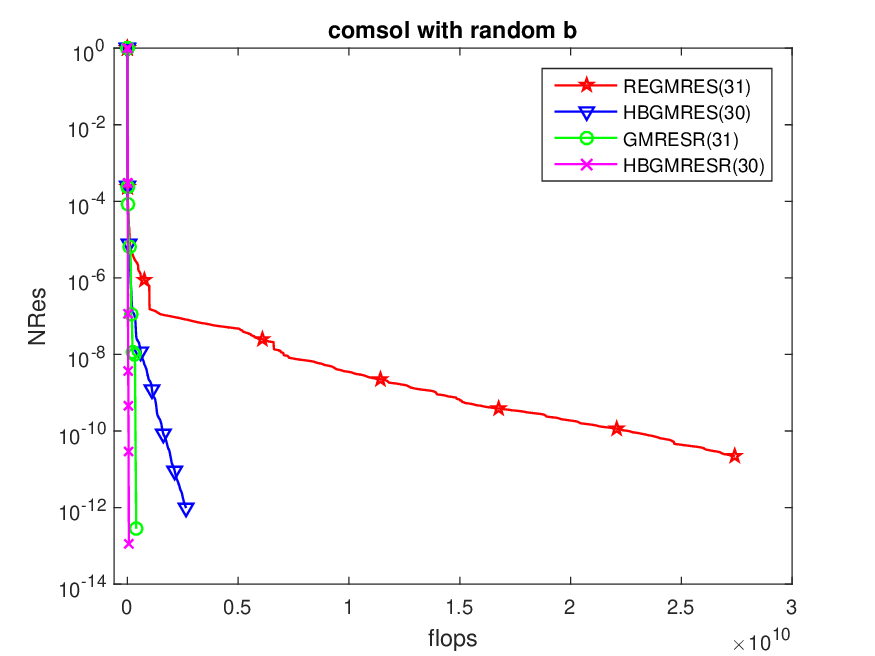}} \\
\hspace{-0.3 cm}
\resizebox*{0.48\textwidth}{0.240\textheight}{\includegraphics{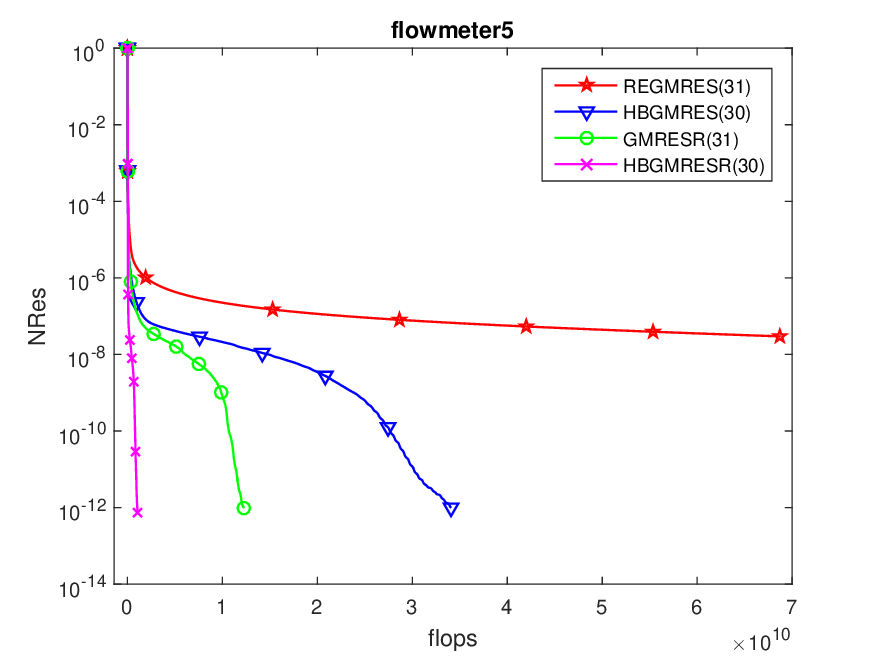}}
&  \hspace{-0.5 cm}
\resizebox*{0.48\textwidth}{0.240\textheight}{\includegraphics{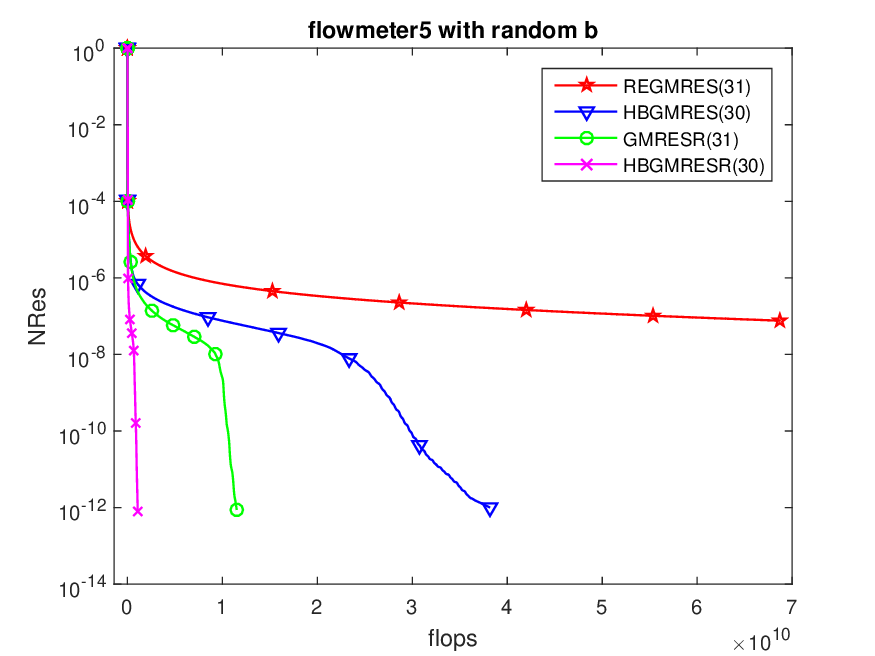}} \\
\hspace{-0.3 cm}
\resizebox*{0.48\textwidth}{0.240\textheight}{\includegraphics{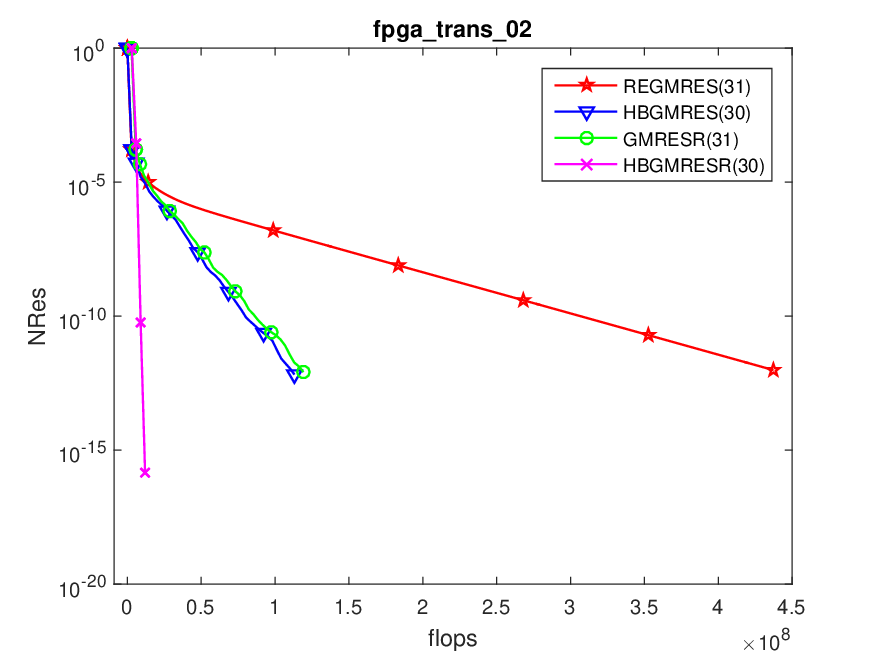}}
&  \hspace{-0.5 cm}
\resizebox*{0.48\textwidth}{0.240\textheight}{\includegraphics{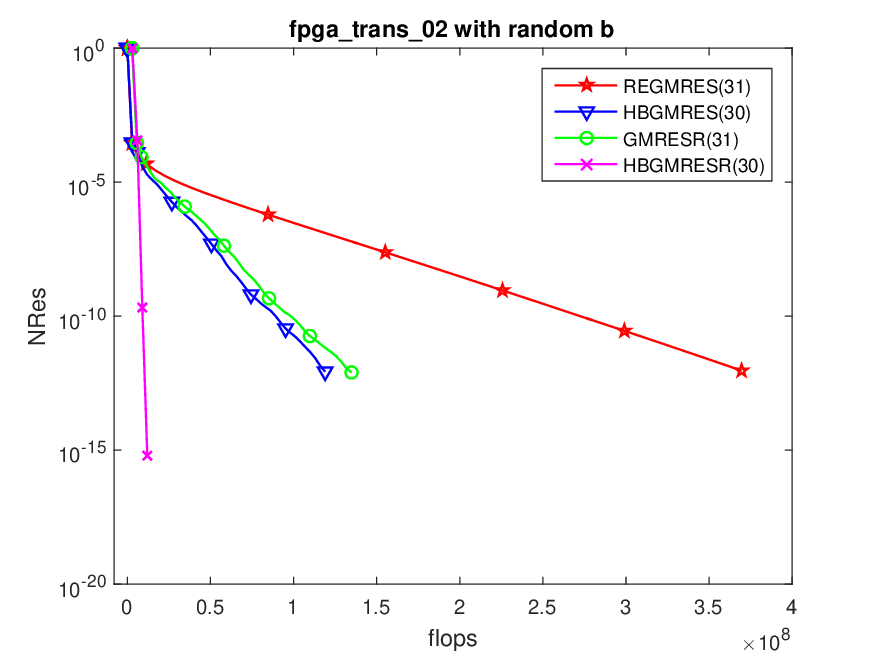}} \\
\hspace{-0.3 cm}
\resizebox*{0.48\textwidth}{0.240\textheight}{\includegraphics{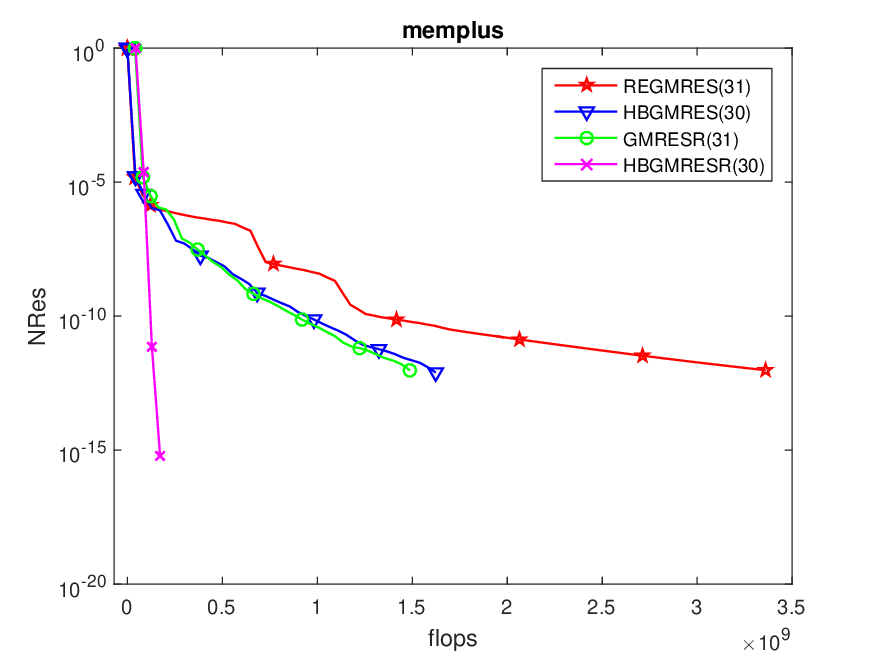}}
&  \hspace{-0.5 cm}
\resizebox*{0.48\textwidth}{0.240\textheight}{\includegraphics{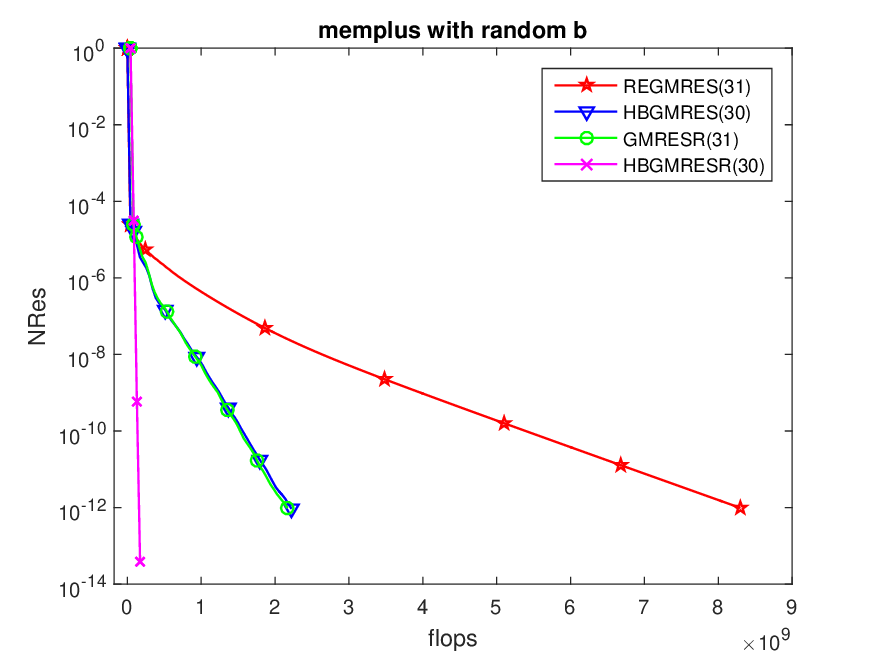}}
\end{tabular}\par
}\vspace{-0.15 cm}
\caption{\small 
    NRes {\em vs.} cycle for
   {\tt comsol}, {\tt flowmeter5}, {\tt fpga\_trans\_02}, and {\tt memplus}
   with selective reorthogonalization. {\em Left:\/} original $b$; {\em Right:\/} random $b$.
   }
\label{fig:2nd4flops}
\end{figure}
\vspace{2mm}

\newpage
\begin{figure}
{\centering
\begin{tabular}{cc}
\hspace{-0.3 cm}
\resizebox*{0.48\textwidth}{0.240\textheight}{\includegraphics{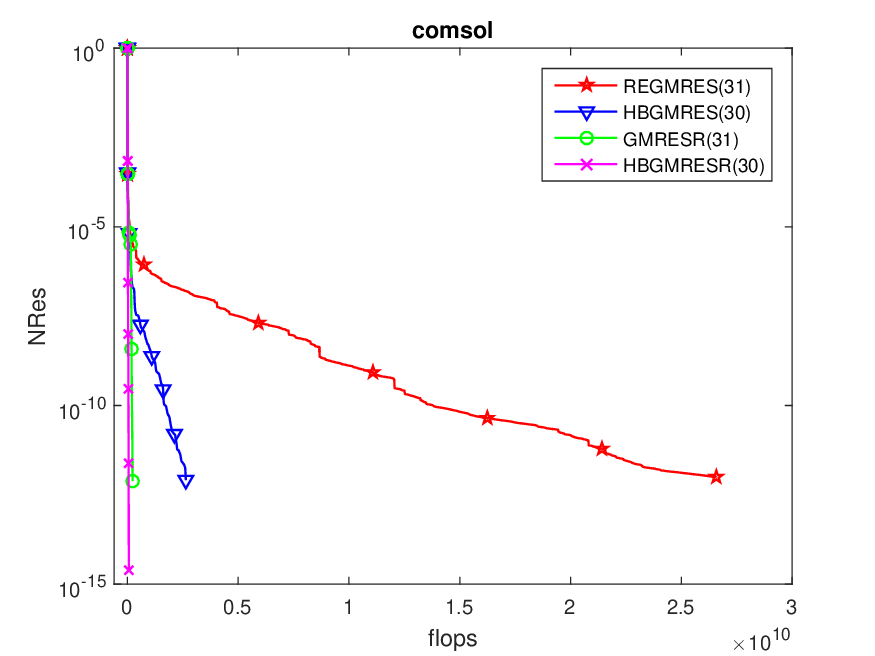}}
&  \hspace{-0.5 cm}
\resizebox*{0.48\textwidth}{0.240\textheight}{\includegraphics{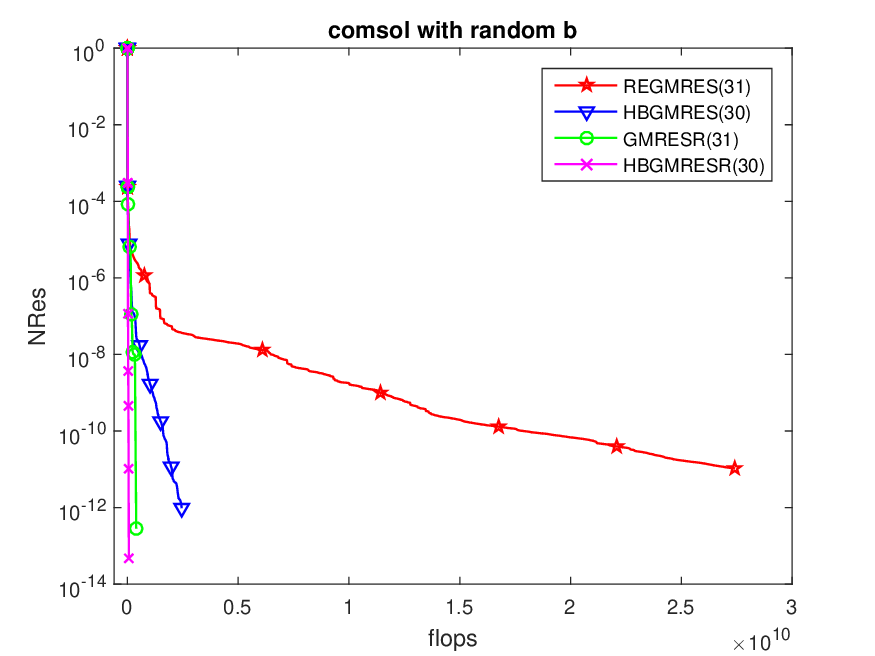}} \\
\hspace{-0.3 cm}
\resizebox*{0.48\textwidth}{0.240\textheight}{\includegraphics{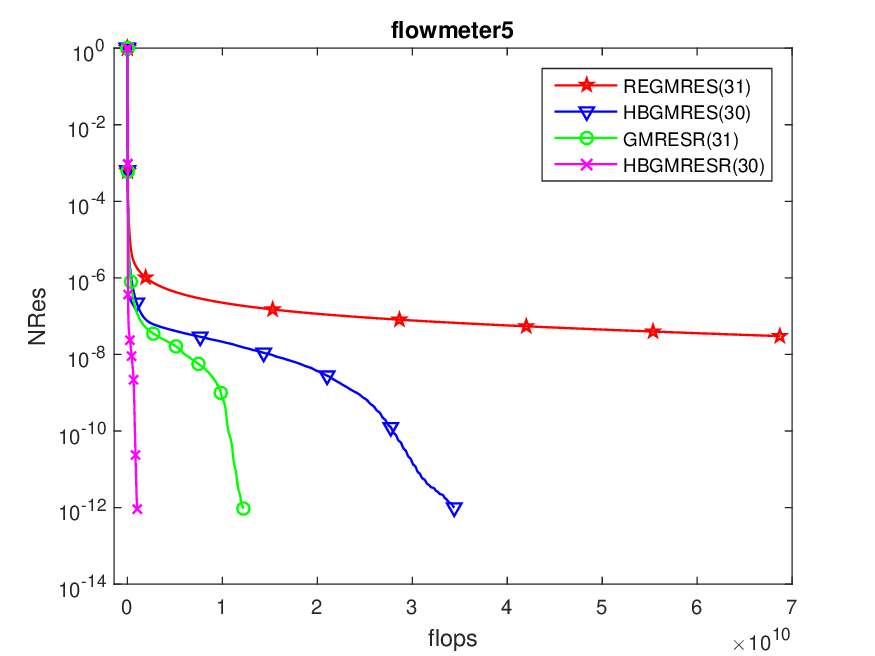}}
&  \hspace{-0.5 cm}
\resizebox*{0.48\textwidth}{0.240\textheight}{\includegraphics{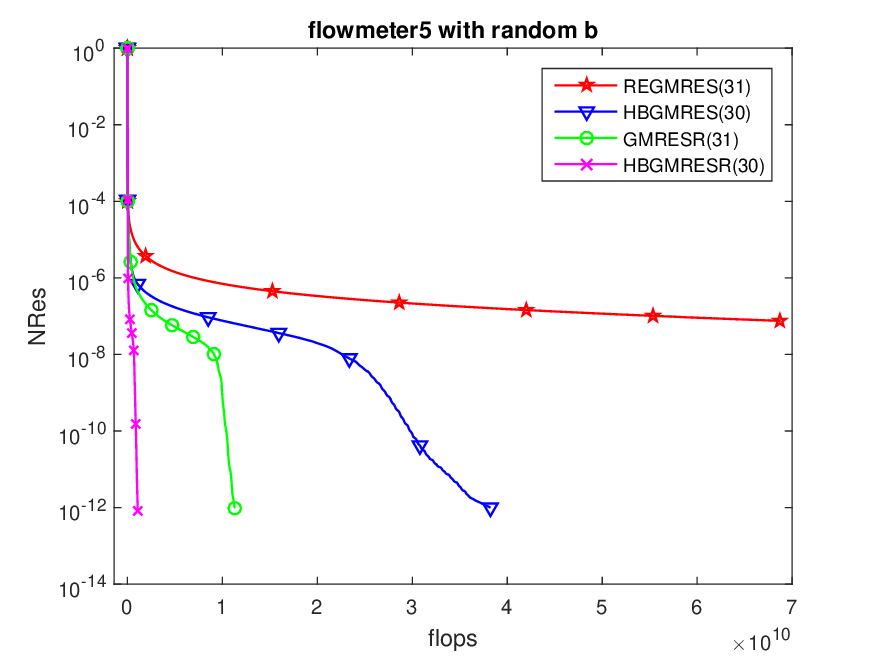}} \\
\hspace{-0.3 cm}
\resizebox*{0.48\textwidth}{0.240\textheight}{\includegraphics{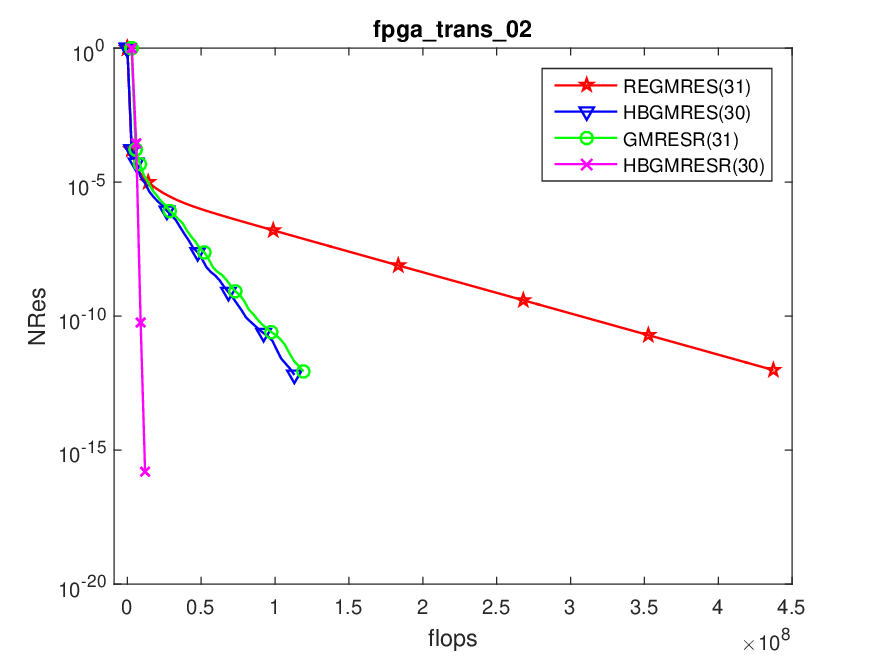}}
&  \hspace{-0.5 cm}
\resizebox*{0.48\textwidth}{0.240\textheight}{\includegraphics{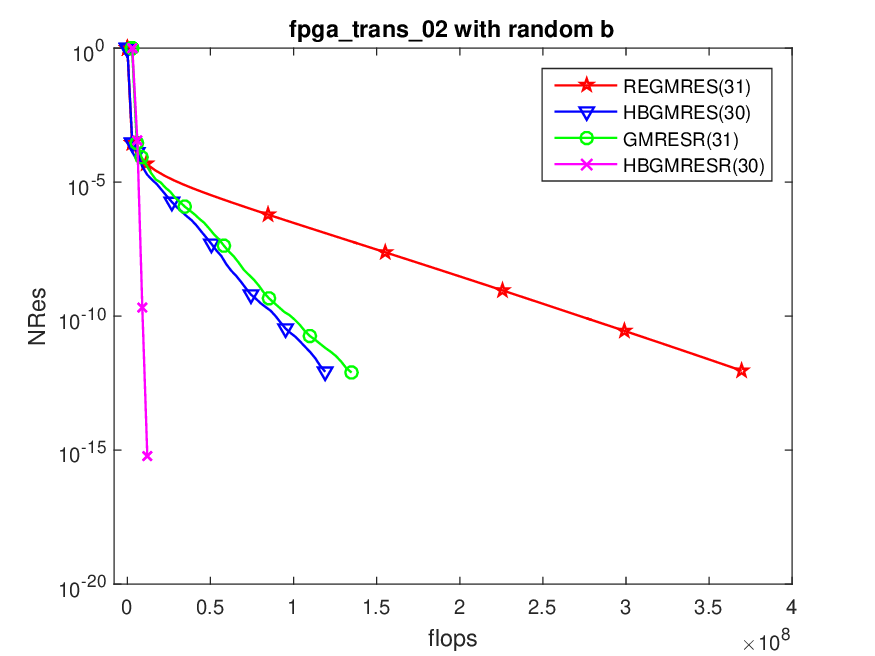}} \\
\hspace{-0.3 cm}
\resizebox*{0.48\textwidth}{0.240\textheight}{\includegraphics{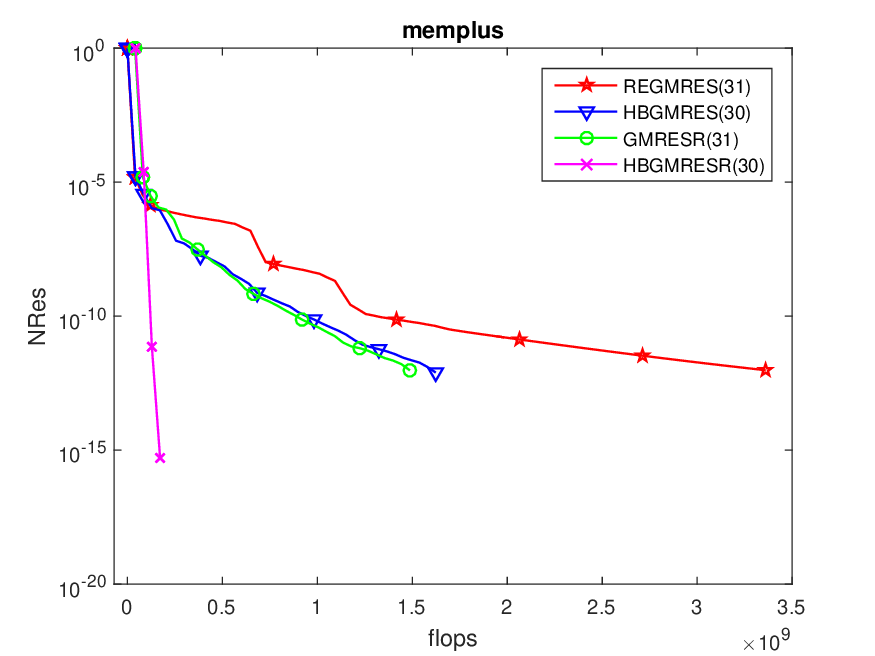}}
&  \hspace{-0.5 cm}
\resizebox*{0.48\textwidth}{0.240\textheight}{\includegraphics{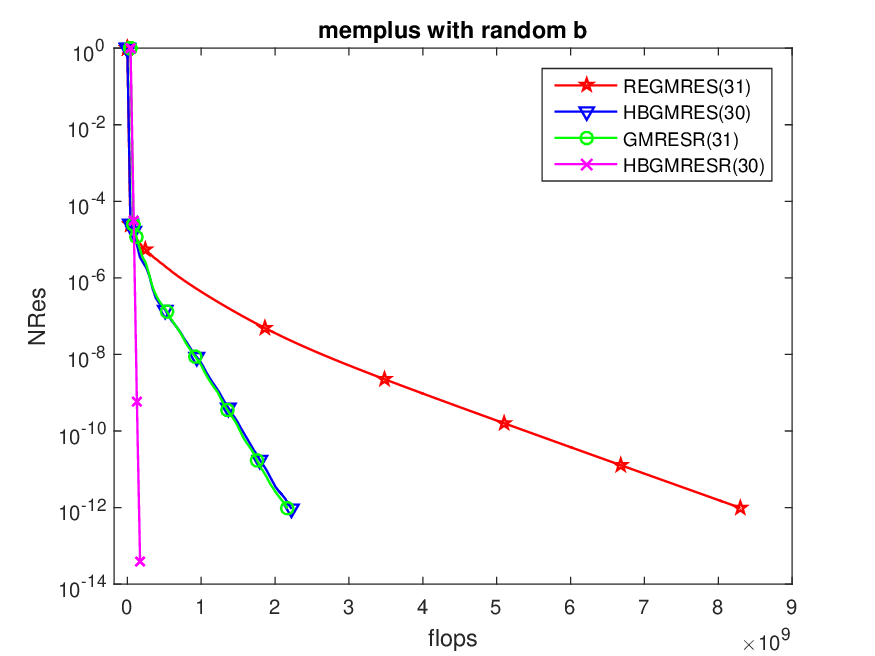}}
\end{tabular}\par
}\vspace{-0.15 cm}
\caption{\small 
    NRes {\em vs.} cycle for
   {\tt comsol}, {\tt flowmeter5}, {\tt fpga\_trans\_02}, and {\tt memplus}
   with always reorthogonalization. {\em Left:\/} original $b$; {\em Right:\/} random $b$.
   }
\label{fig:2nd4wflops}
\end{figure}

\newpage
\begin{figure}
{\centering
\begin{tabular}{cc}
\hspace{-0.3 cm}
\resizebox*{0.48\textwidth}{0.240\textheight}{\includegraphics{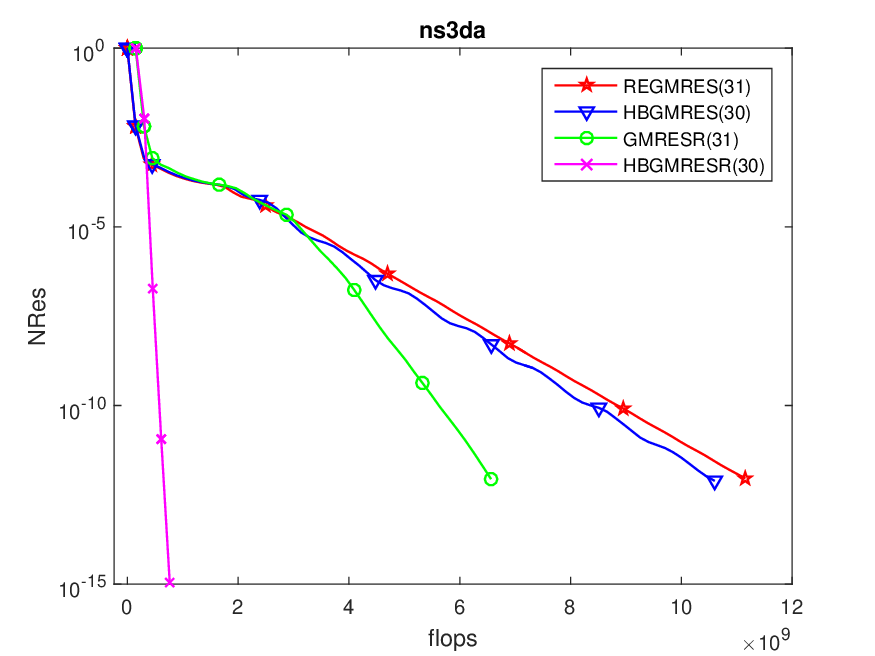}}
&  \hspace{-0.5 cm}
\resizebox*{0.48\textwidth}{0.240\textheight}{\includegraphics{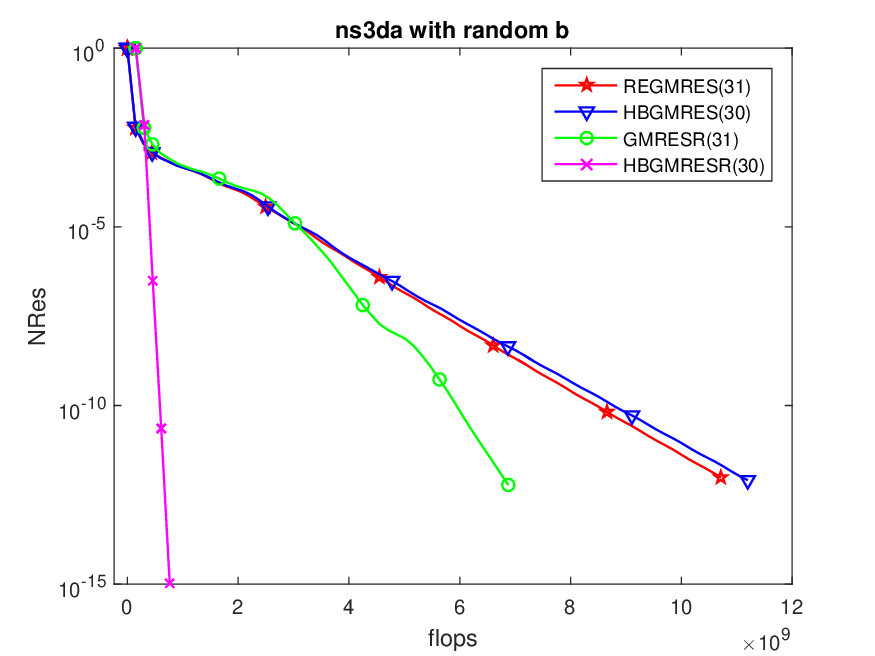}} \\
\hspace{-0.3 cm}
\resizebox*{0.48\textwidth}{0.240\textheight}{\includegraphics{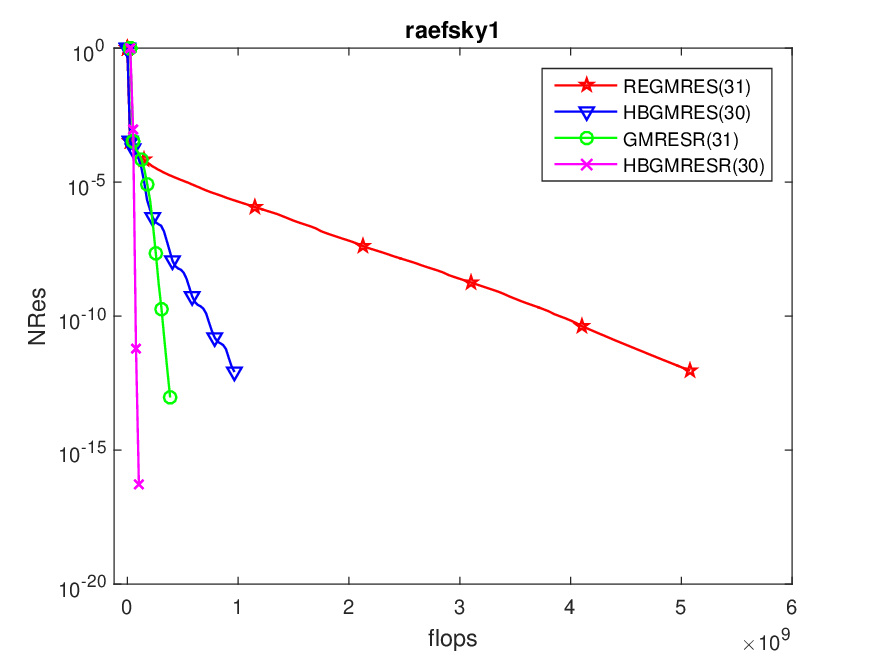}}
&  \hspace{-0.5 cm}
\resizebox*{0.48\textwidth}{0.240\textheight}{\includegraphics{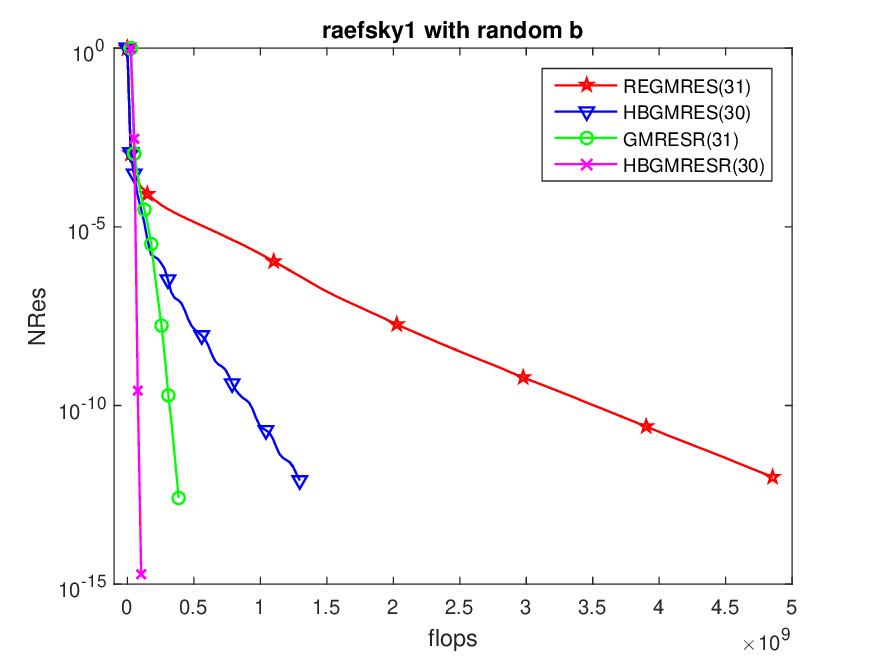}} \\
\hspace{-0.3 cm}
\resizebox*{0.48\textwidth}{0.240\textheight}{\includegraphics{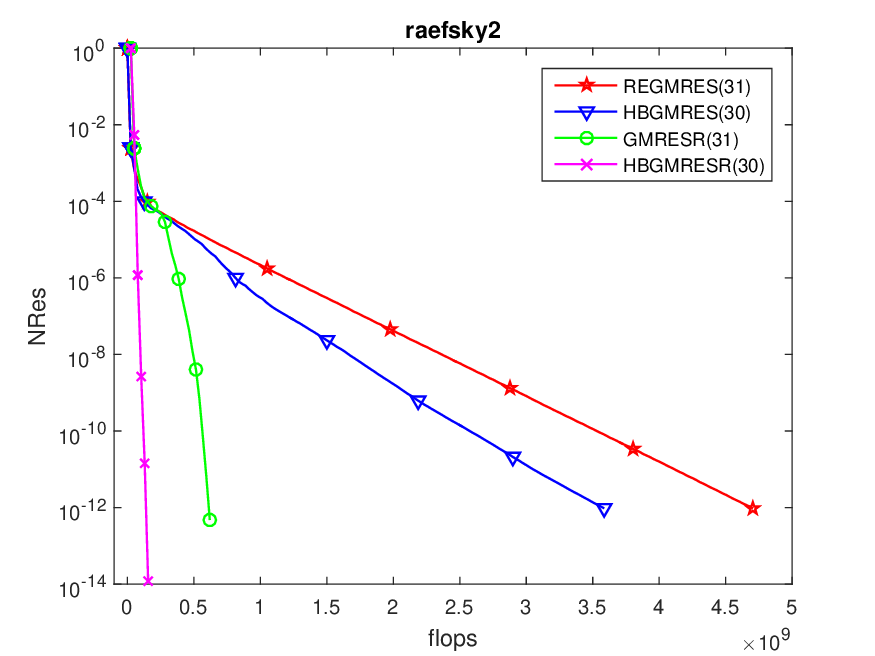}}
&  \hspace{-0.5 cm}
\resizebox*{0.48\textwidth}{0.240\textheight}{\includegraphics{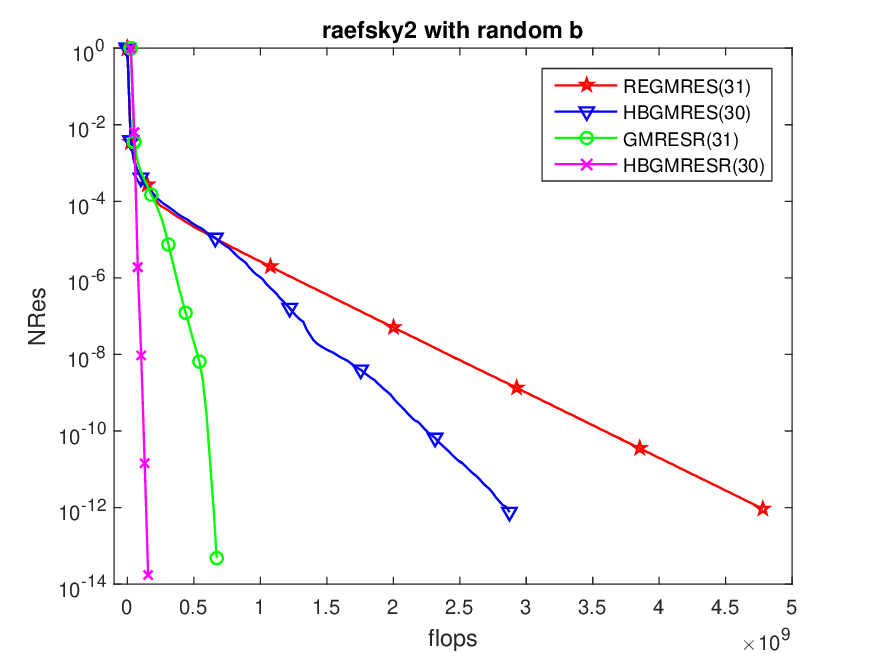}} \\
\hspace{-0.3 cm}
\resizebox*{0.48\textwidth}{0.240\textheight}{\includegraphics{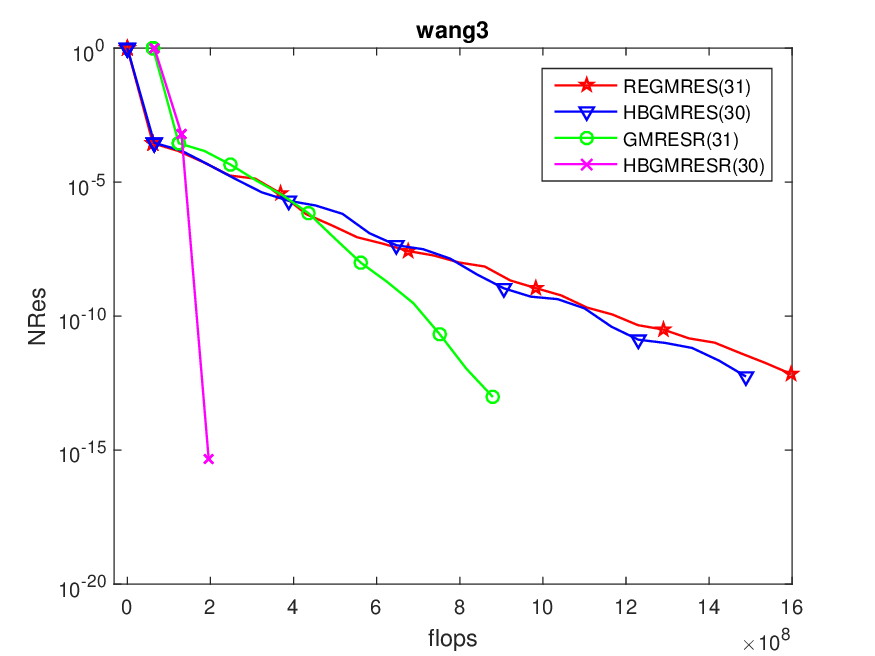}}
&  \hspace{-0.5 cm}
\resizebox*{0.48\textwidth}{0.240\textheight}{\includegraphics{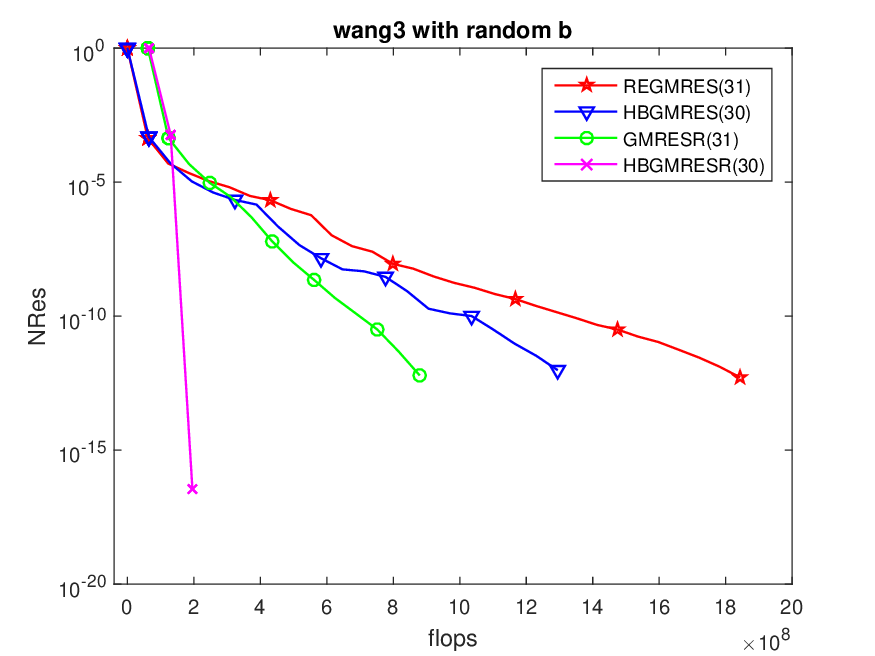}}
\end{tabular}\par
}\vspace{-0.15 cm}
\caption{\small 
    NRes {\em vs.} cycle for
   {\tt ns3da}, {\tt raefsky1}, {\tt raefsky2}, and {\tt wang3}
   with selective reorthogonalization. {\em Left:\/} original $b$; {\em Right:\/} random $b$.
   }
\label{fig:3rd4flops}
\end{figure}
\vspace{2mm}

\newpage
\begin{figure}
{\centering
\begin{tabular}{cc}
\hspace{-0.3 cm}
\resizebox*{0.48\textwidth}{0.240\textheight}{\includegraphics{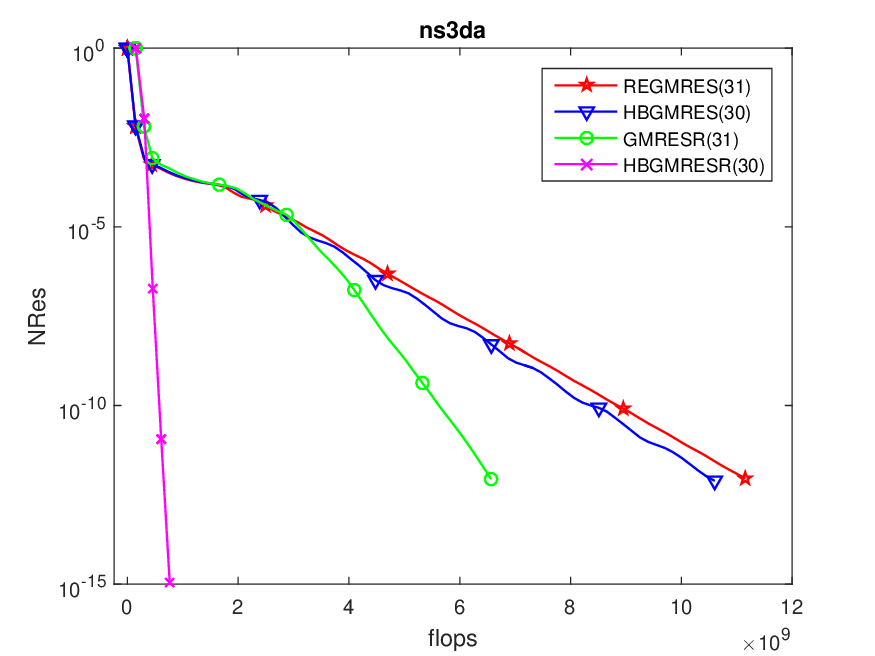}}
&  \hspace{-0.5 cm}
\resizebox*{0.48\textwidth}{0.240\textheight}{\includegraphics{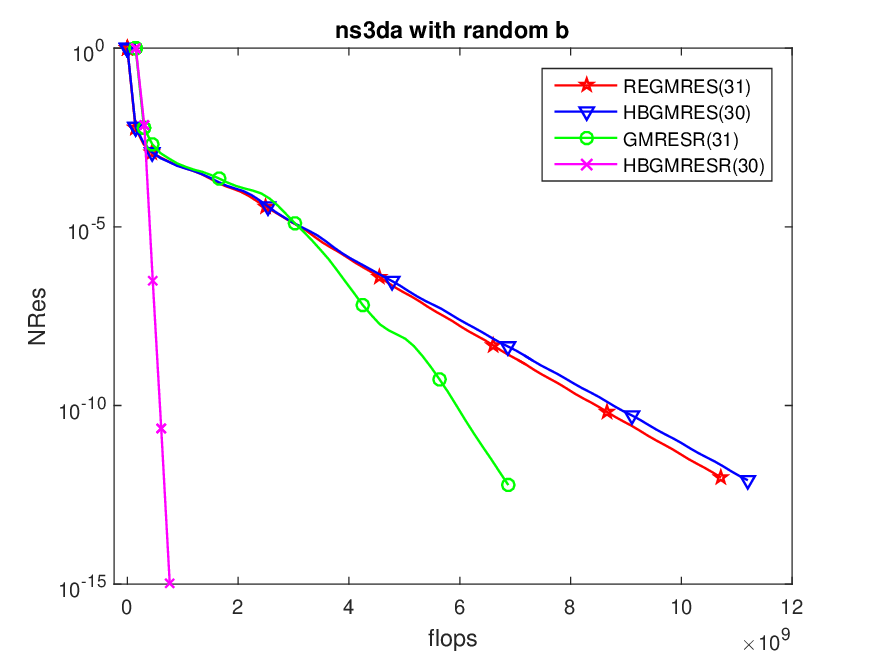}} \\
\hspace{-0.3 cm}
\resizebox*{0.48\textwidth}{0.240\textheight}{\includegraphics{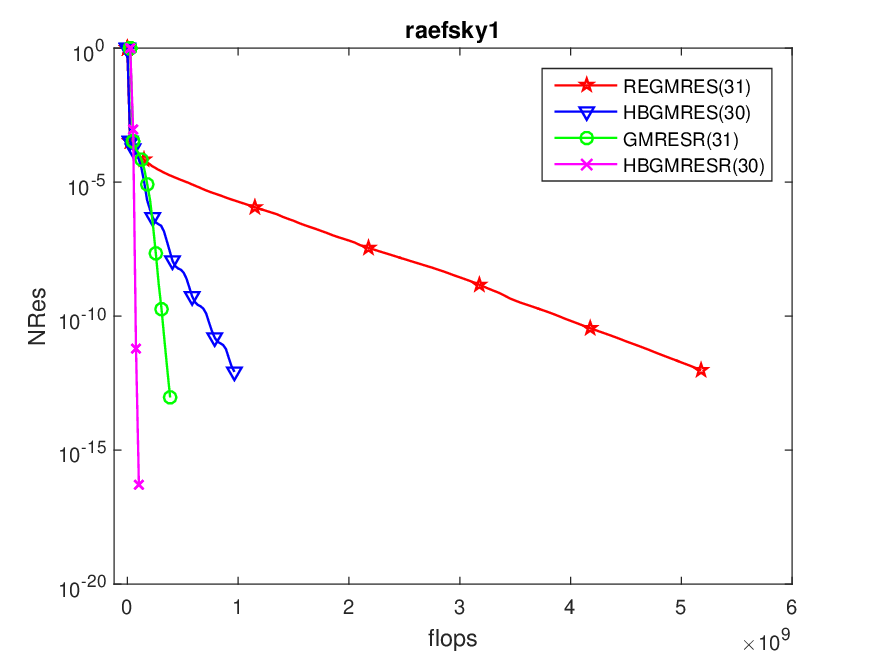}}
&  \hspace{-0.5 cm}
\resizebox*{0.48\textwidth}{0.240\textheight}{\includegraphics{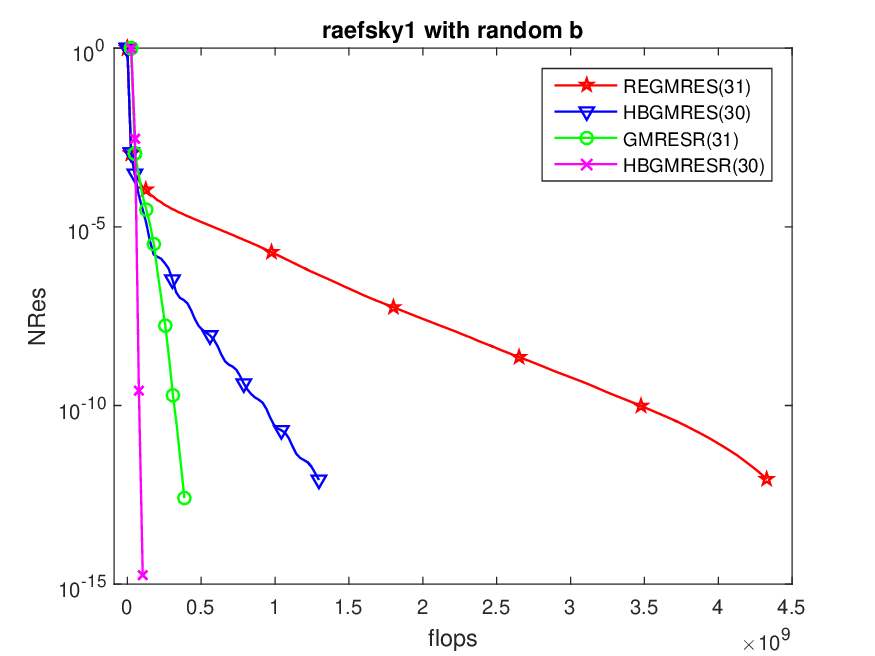}} \\
\hspace{-0.3 cm}
\resizebox*{0.48\textwidth}{0.240\textheight}{\includegraphics{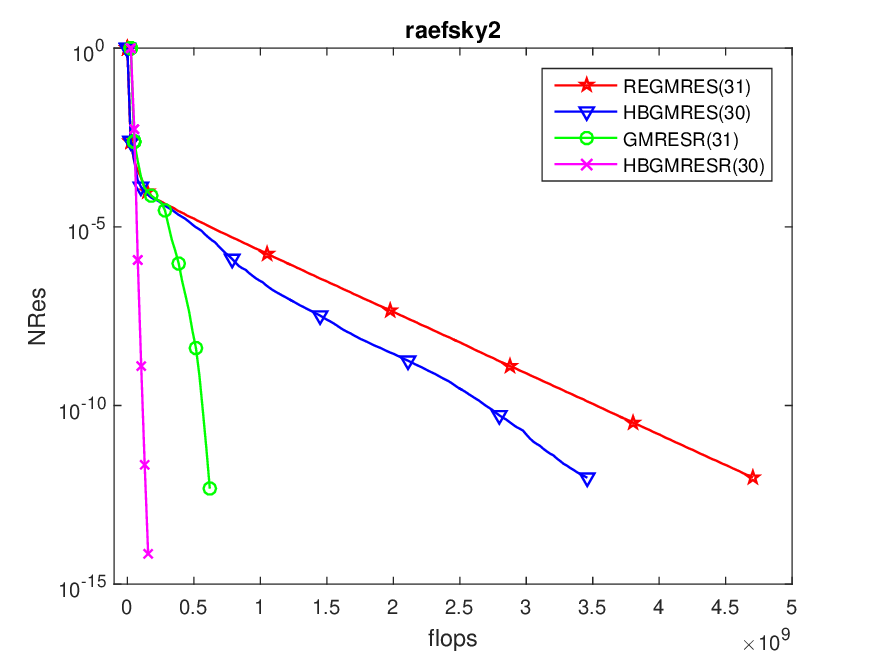}}
&  \hspace{-0.5 cm}
\resizebox*{0.48\textwidth}{0.240\textheight}{\includegraphics{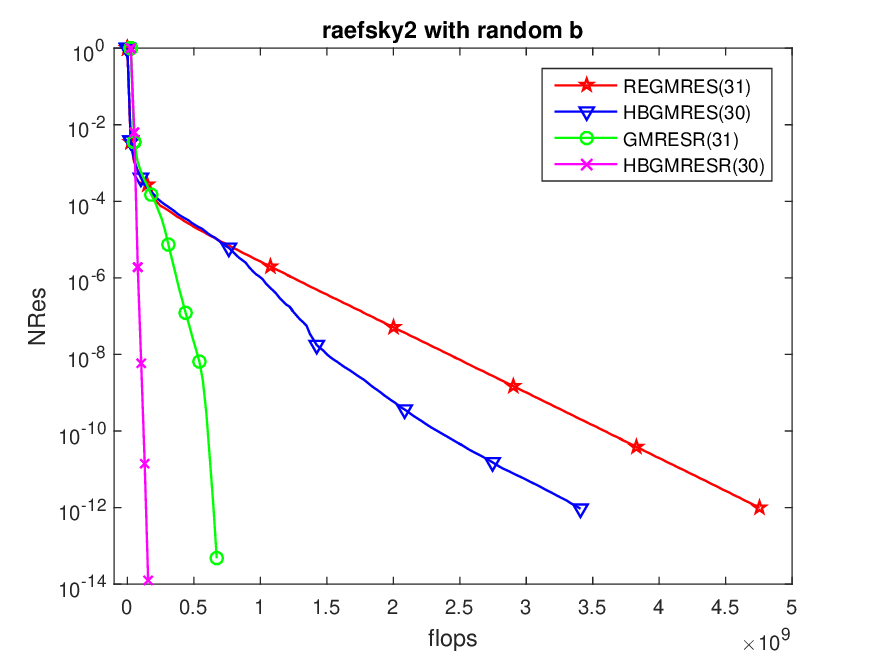}} \\
\hspace{-0.3 cm}
\resizebox*{0.48\textwidth}{0.240\textheight}{\includegraphics{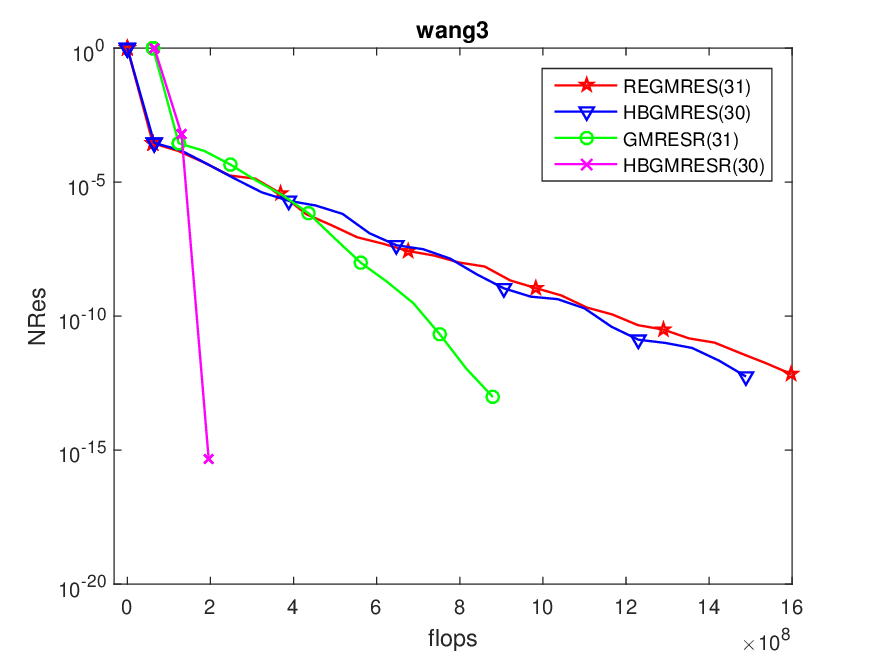}}
&  \hspace{-0.5 cm}
\resizebox*{0.48\textwidth}{0.240\textheight}{\includegraphics{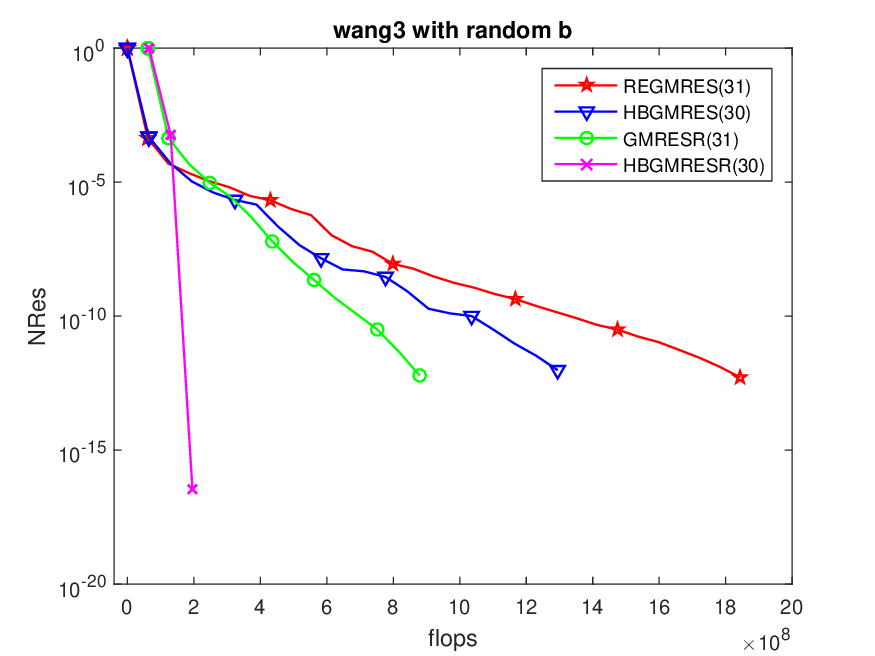}}
\end{tabular}\par
}\vspace{-0.15 cm}
\caption{\small 
    NRes {\em vs.} cycle for
   {\tt ns3da}, {\tt raefsky1}, {\tt raefsky2}, and {\tt wang3}
   with always reorthogonalization. {\em Left:\/} original $b$; {\em Right:\/} random $b$.
   }
\label{fig:3rd4wflops}
\end{figure}
\vspace{2mm}

\newpage
\begin{figure}
{\centering
\begin{tabular}{cc}
\hspace{-0.3 cm}
\resizebox*{0.48\textwidth}{0.240\textheight}{\includegraphics{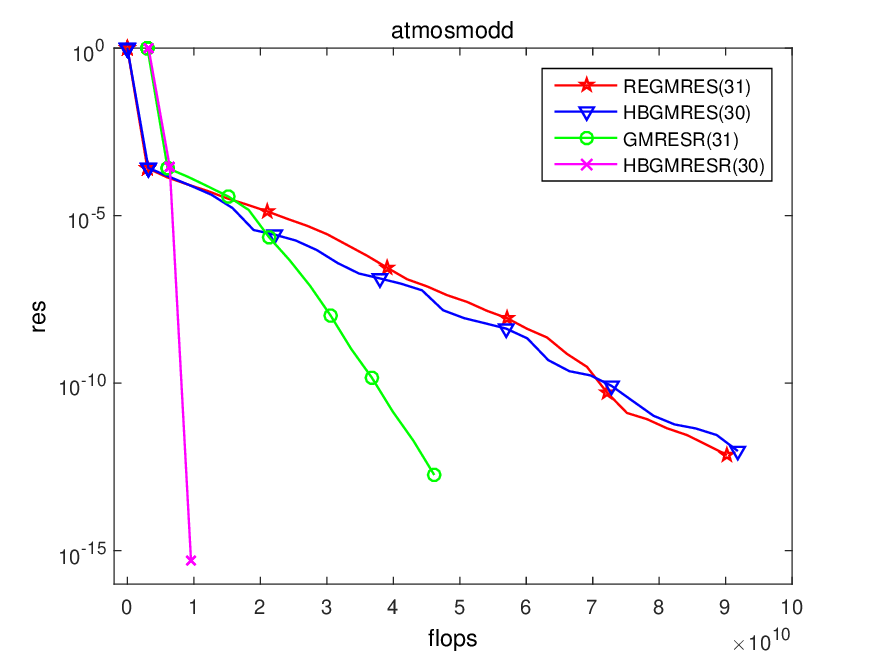}}
&  \hspace{-0.5 cm}
\resizebox*{0.48\textwidth}{0.240\textheight}{\includegraphics{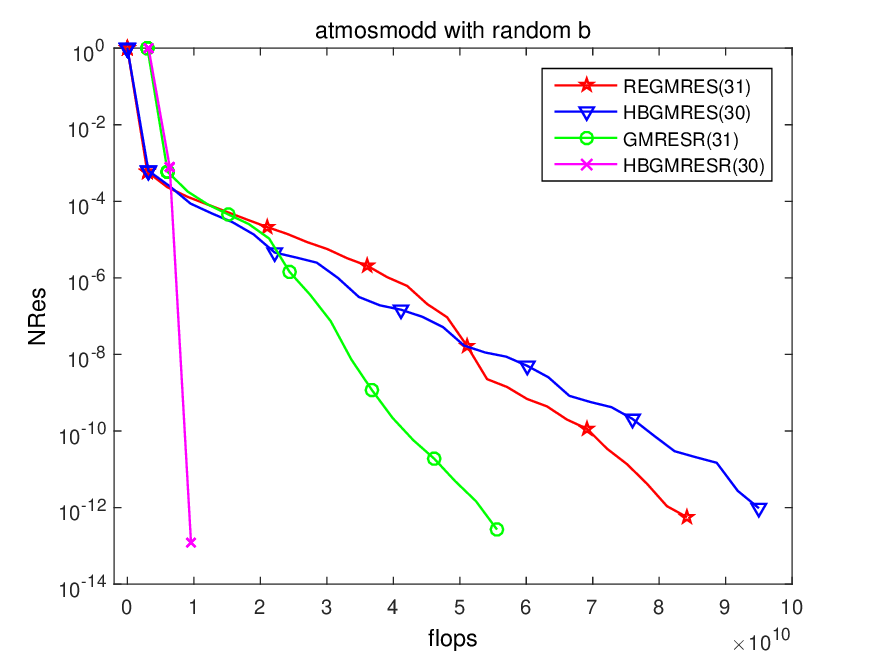}} \\
\hspace{-0.3 cm}
\resizebox*{0.48\textwidth}{0.240\textheight}{\includegraphics{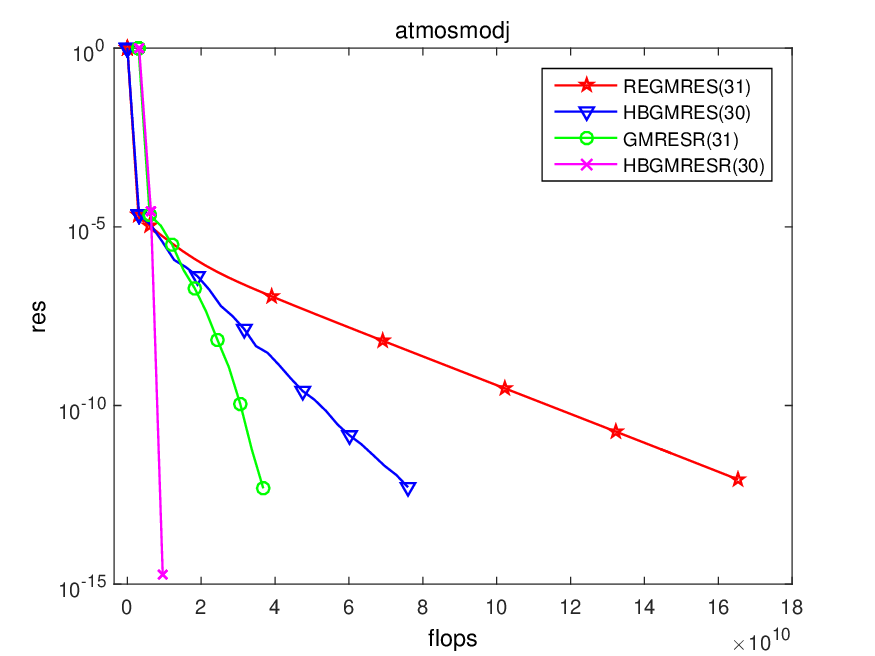}}
&  \hspace{-0.5 cm}
\resizebox*{0.48\textwidth}{0.240\textheight}{\includegraphics{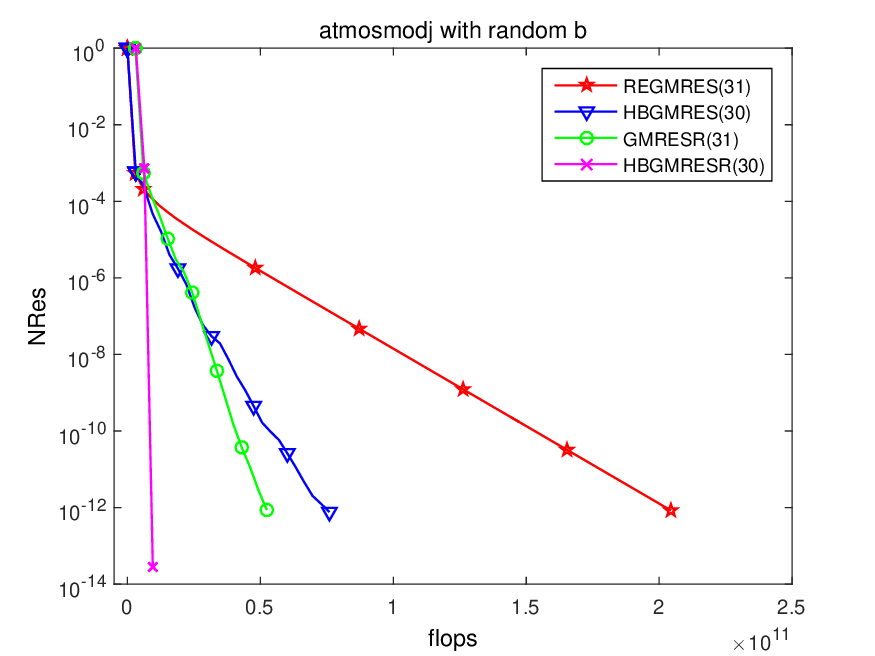}} \\
\hspace{-0.3 cm}
\resizebox*{0.48\textwidth}{0.240\textheight}{\includegraphics{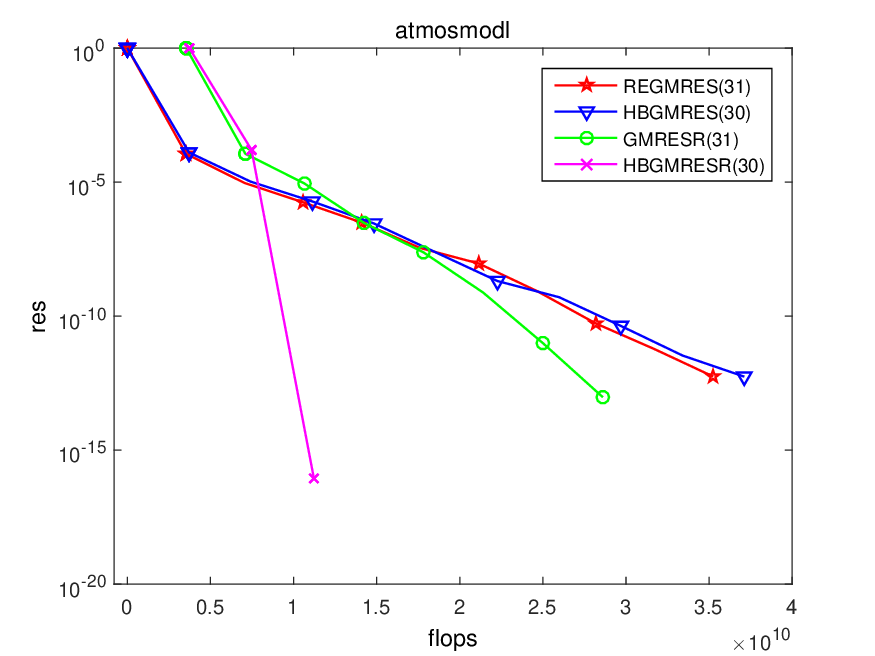}}
&  \hspace{-0.5 cm}
\resizebox*{0.48\textwidth}{0.240\textheight}{\includegraphics{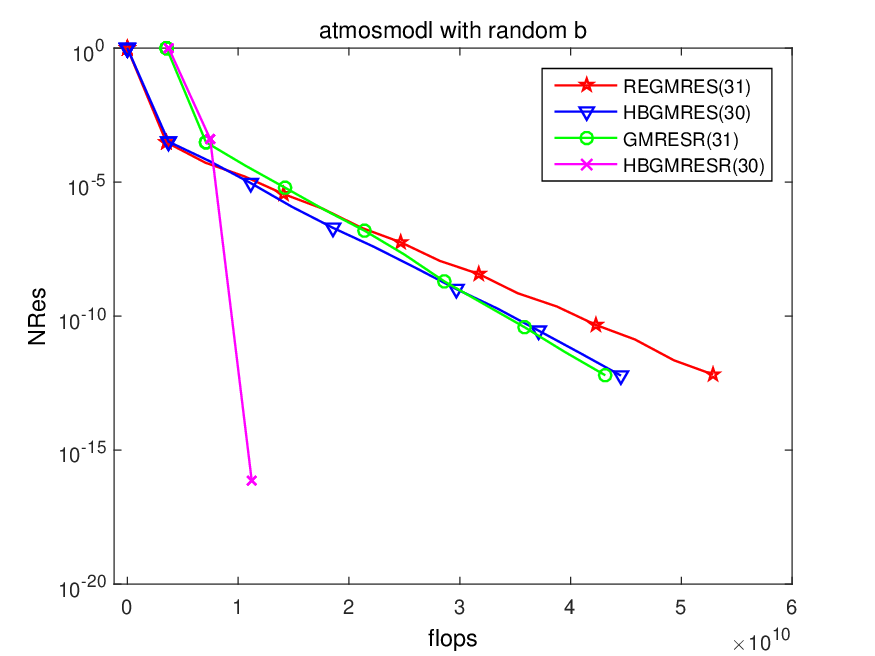}} \\
\hspace{-0.3 cm}
\resizebox*{0.48\textwidth}{0.240\textheight}{\includegraphics{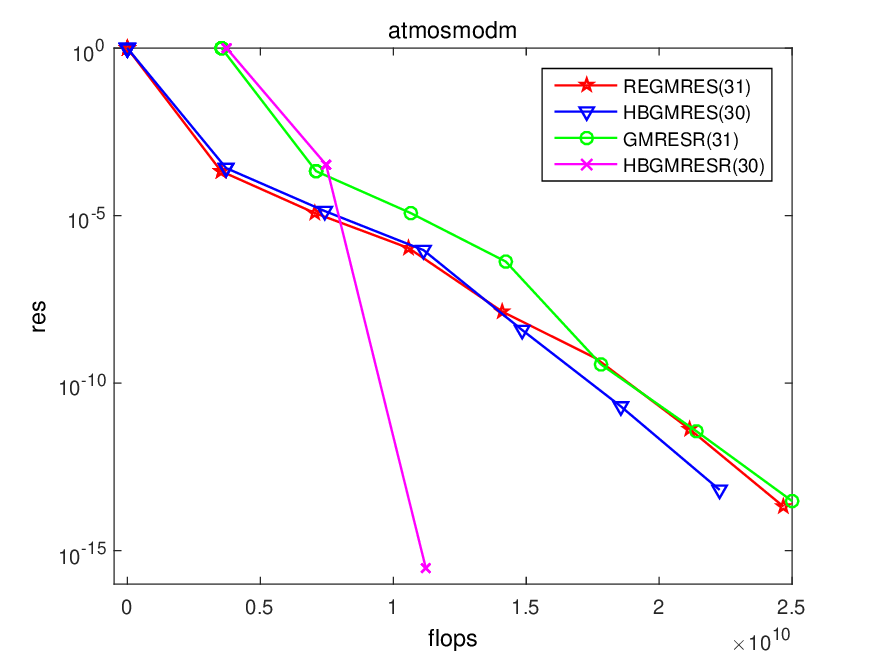}}
&  \hspace{-0.5 cm}
\resizebox*{0.48\textwidth}{0.240\textheight}{\includegraphics{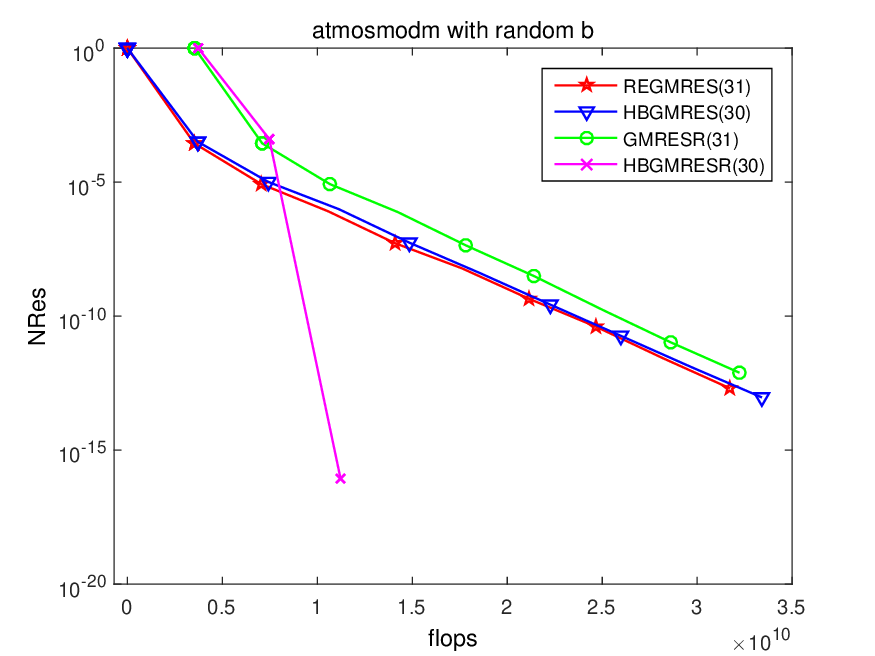}}
\end{tabular}\par
}\vspace{-0.15 cm}
\caption{\small 
    NRes {\em vs.} cycle for
   {\tt atmosmodd}, {\tt atmosmodj}, {\tt atmosmodl}, and {\tt atmosmodm}
   with selective reorthogonalization. {\em Left:\/} original $b$; {\em Right:\/} random $b$.
   }
\label{fig:4th4flops}
\end{figure}

\newpage
\begin{figure}
{\centering
\begin{tabular}{cc}
\hspace{-0.3 cm}
\resizebox*{0.48\textwidth}{0.240\textheight}{\includegraphics{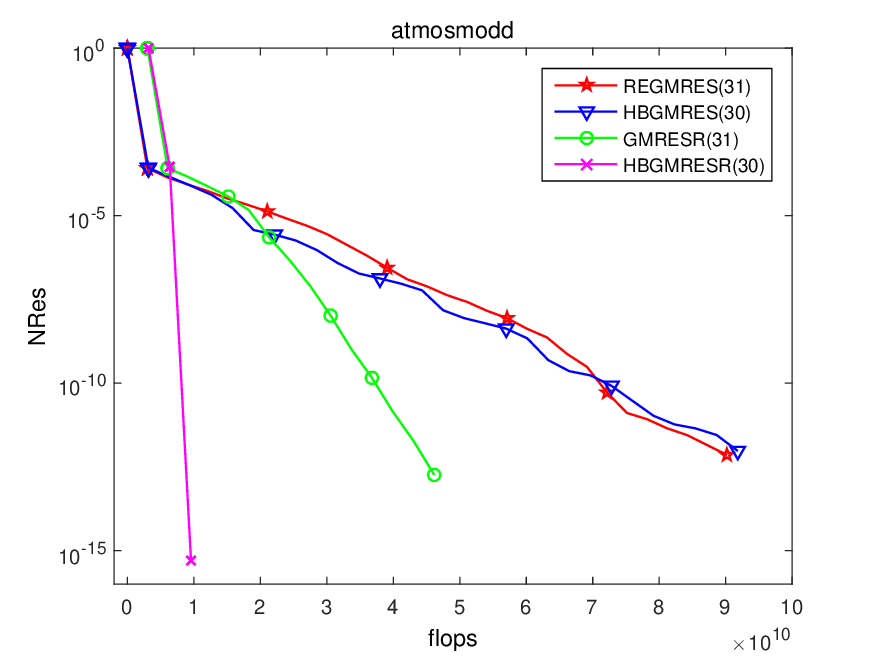}}
&  \hspace{-0.5 cm}
\resizebox*{0.48\textwidth}{0.240\textheight}{\includegraphics{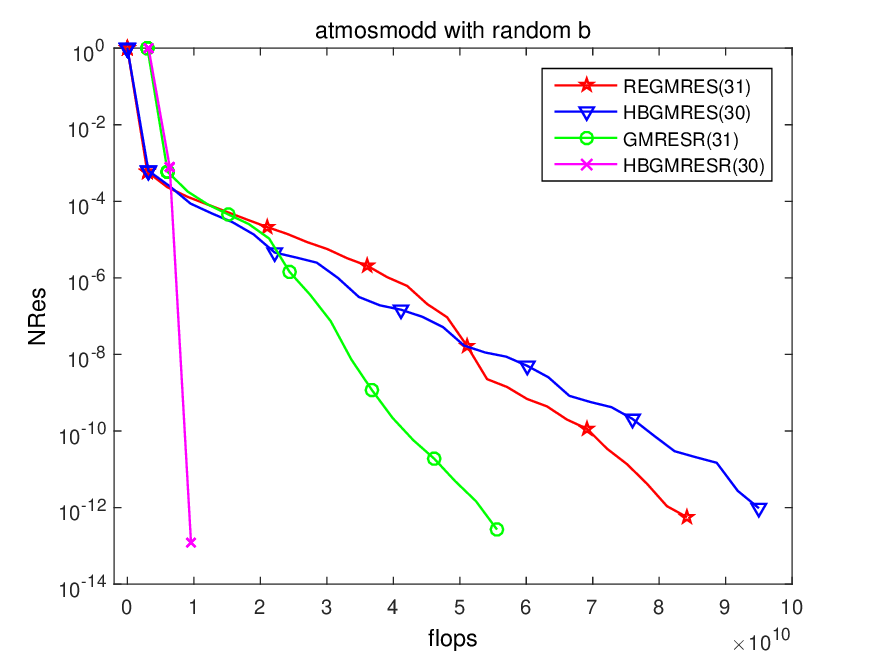}} \\
\hspace{-0.3 cm}
\resizebox*{0.48\textwidth}{0.240\textheight}{\includegraphics{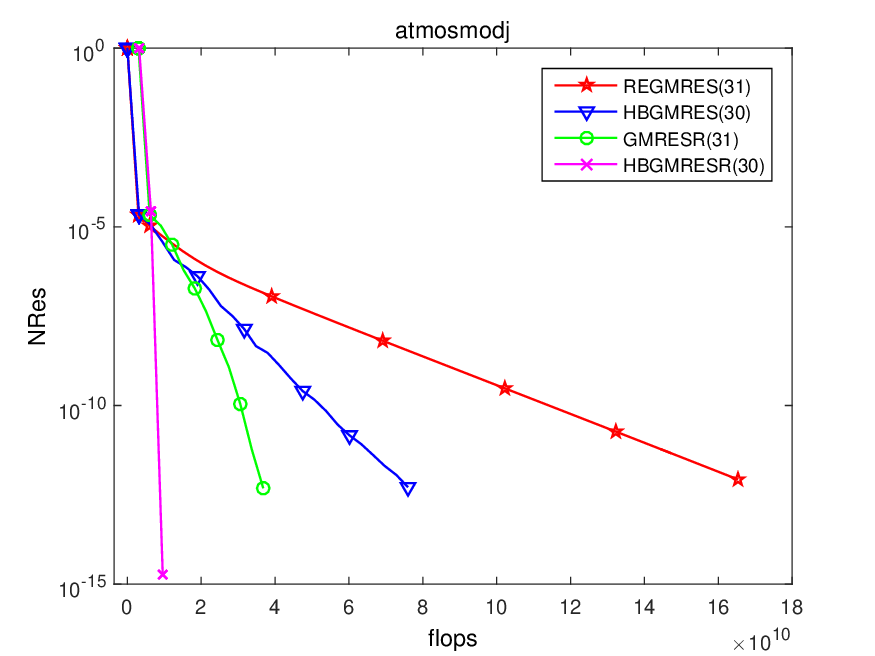}}
&  \hspace{-0.5 cm}
\resizebox*{0.48\textwidth}{0.240\textheight}{\includegraphics{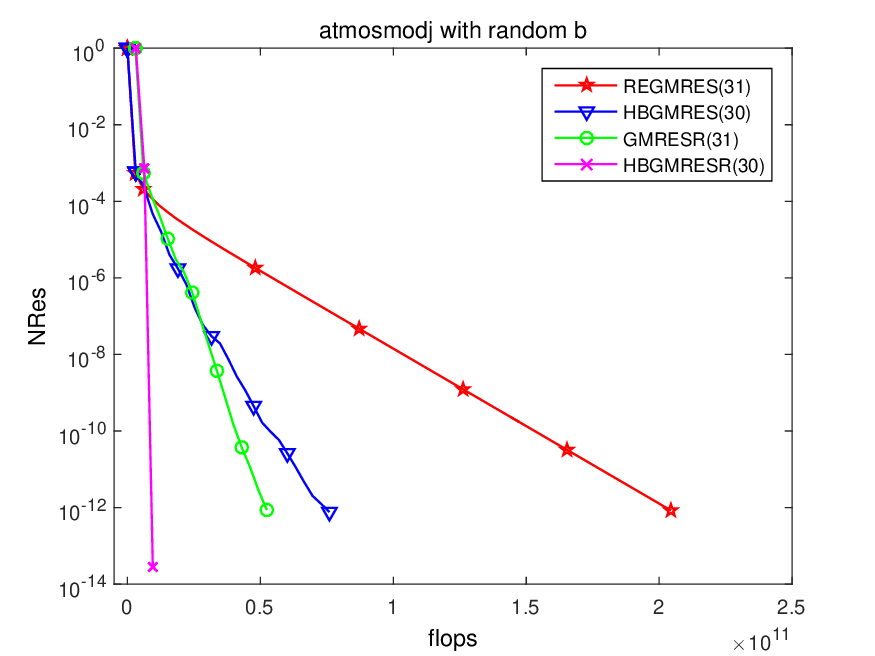}} \\
\hspace{-0.3 cm}
\resizebox*{0.48\textwidth}{0.240\textheight}{\includegraphics{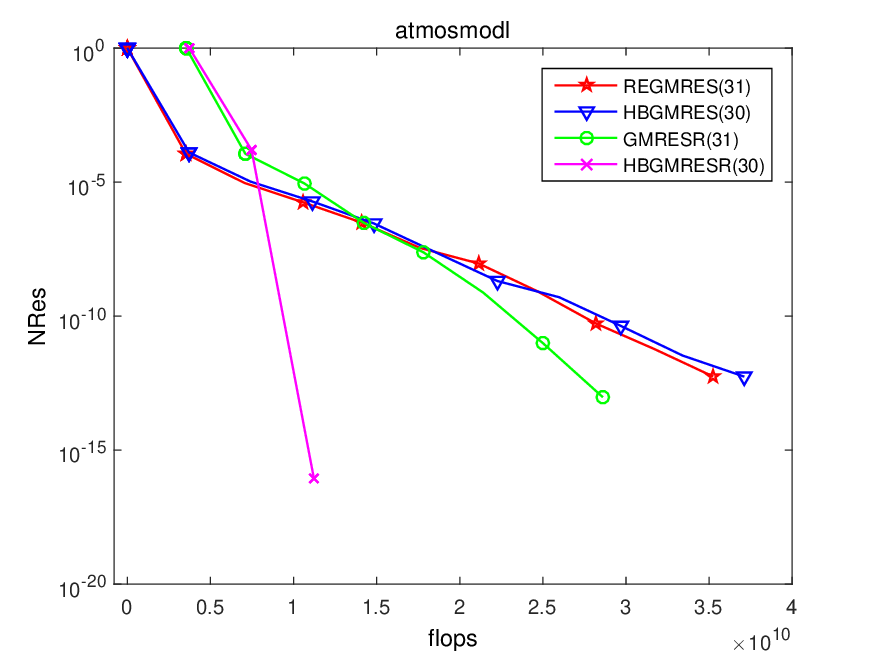}}
&  \hspace{-0.5 cm}
\resizebox*{0.48\textwidth}{0.240\textheight}{\includegraphics{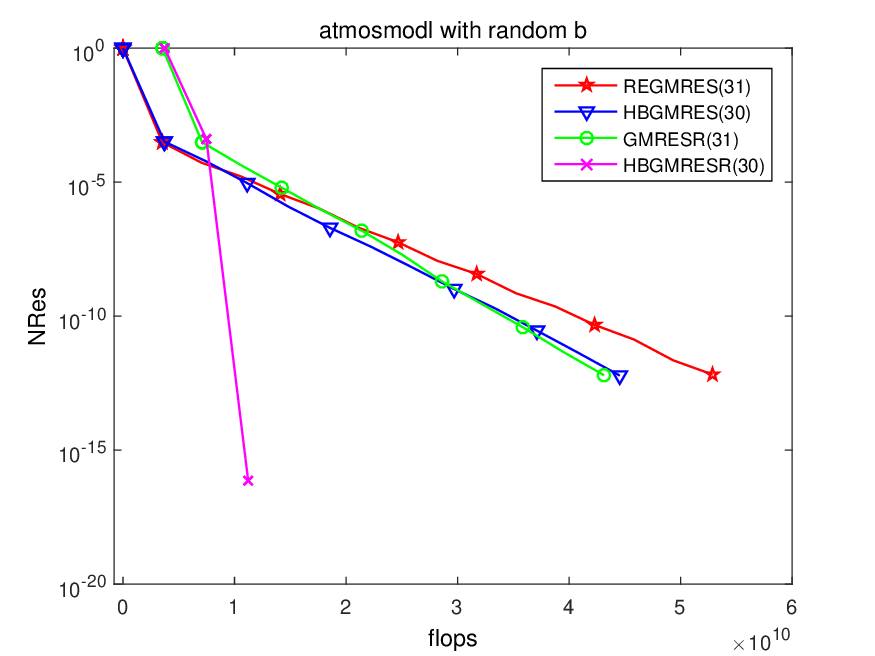}} \\
\hspace{-0.3 cm}
\resizebox*{0.48\textwidth}{0.240\textheight}{\includegraphics{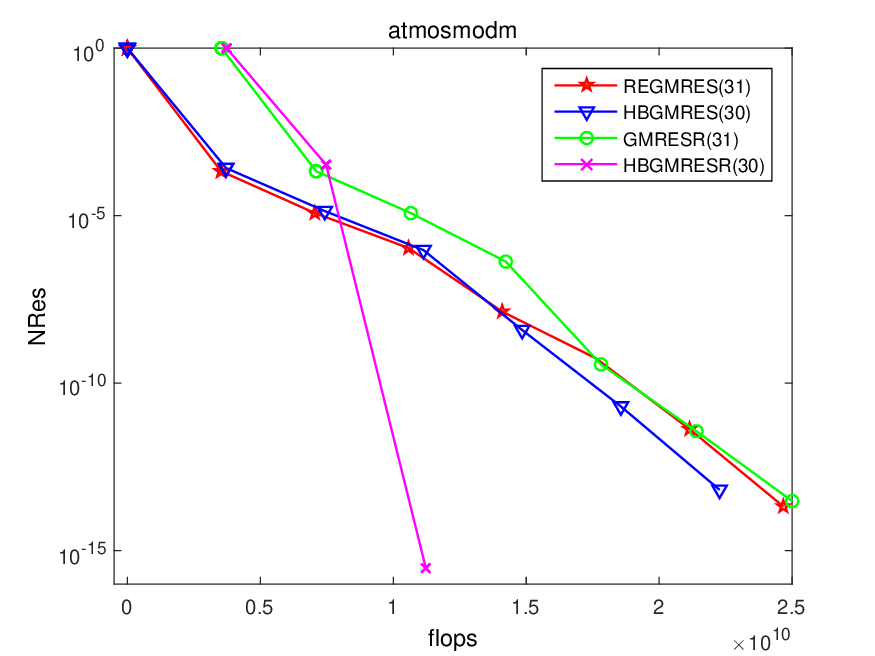}}
&  \hspace{-0.5 cm}
\resizebox*{0.48\textwidth}{0.240\textheight}{\includegraphics{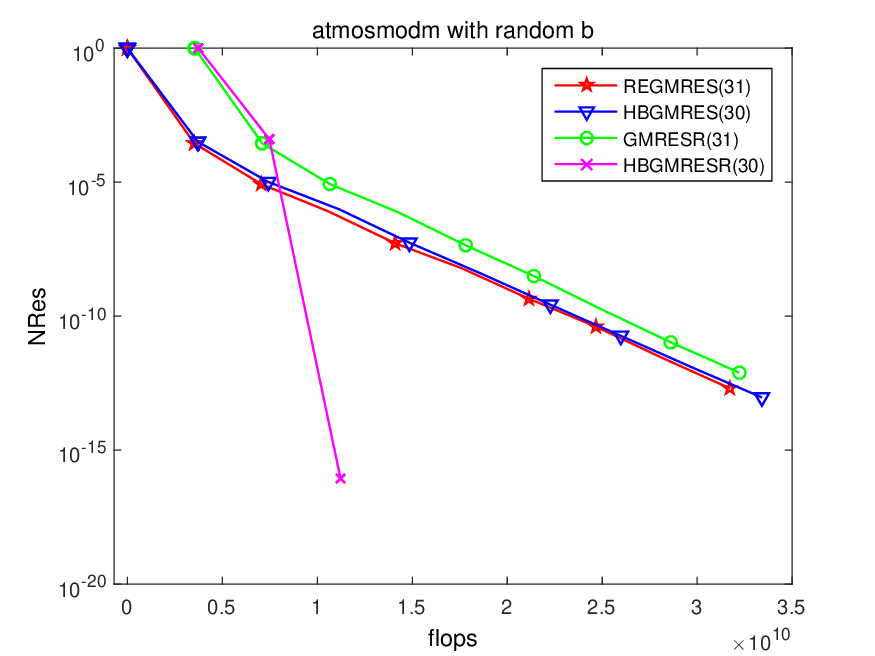}}
\end{tabular}\par
}\vspace{-0.15 cm}
\caption{\small 
    NRes {\em vs.} cycle for
   {\tt atmosmodd}, {\tt atmosmodj}, {\tt atmosmodl}, and {\tt atmosmodm}
   with selective reorthogonalization. {\em Left:\/} original $b$; {\em Right:\/} random $b$.
   }
\label{fig:4th4wflops}
\end{figure}

\bibliographystyle{amsplain}

\end{document}